# Symmetrisation of a class of two-sample tests by mutually considering depth ranks including functional spaces

**Felix Gnettner**[1], **Claudia Kirch**[1] and **Alicia Nieto-Reyes**[2]

[1] *Institut für Mathematische Stochastik, Fakultät für Mathematik, Otto-von-Guericke-Universität Magdeburg, Germany,*
*e-mail:* felix.gnettner@ovgu.de; claudia.kirch@ovgu.de

[2] *Departamento de Matemáticas, Estadística y Computación, Facultad de Ciencias, Universidad de Cantabria, Spain,*
*e-mail:* alicia.nieto@unican.es

**Abstract:** Statistical depth functions provide measures of the outlyingness, or centrality, of the elements of a space with respect to a distribution. It is a nonparametric concept applicable to spaces of any dimension, for instance, multivariate and functional. Liu and Singh (1993) presented a multivariate two-sample test based on depth-ranks. We dedicate this paper to improving the power of the associated test statistic and incorporating its applicability to functional data. In doing so, we obtain a more natural test statistic that is symmetric in both samples. We derive the null asymptotic of the proposed test statistic, also proving the validity of the testing procedure for functional data. Finally, the finite sample performance of the test for functional data is illustrated by means of a simulation study and a real data analysis on annual temperature curves of ocean drifters is executed.

## 1. Introduction

Two-sample testing aims at answering the question whether two independent observed datasets have the same underlying distribution (null hypothesis) or whether the distributions differ (alternative hypothesis). This question plays a role in all applications where data from two independent groups are compared. In univariate statistics, the two commonly used tests for this problem are the two-sample $t$-test for normal distributions and the nonparametric Wilcoxon rank sum test, which is suited for heavy tailed data or data containing outliers. While Hotelling's $T^2$-test provides a suitable generalisation of the two-sample-$t$-test to multivariate data, there is not a unique natural generalisation of the Wilcoxon rank sum test. Generally, two-sample tests can be based on any reasonable measure of discrepancy between the distributions, for example [65] proposed a methodology based on the so-called ball divergence while [41, 91] developed kernel based two-sample tests for deviations in the mean based on the maximum mean discrepancy.

Jurečková and Kalina [46] discussed several multivariate rank based test including a test first proposed by Liu and Singh [52] and later further studied by [13, 71, 79, 93]. The statistic involved in the Liu and Singh test is known as the Liu-Singh-statistic. This statistic is based on the notion of a statistical depth function [94], which provides ranks of the elements of a sample with respect to itself and other given sample(s). The Liu-Singh-statistic detects both location and scale differences in principle but is particularly powerful against pure scale differences. A different test also based on depth ranks was proposed and numerically investigated in Chenouri and Small [16]. Statistical depth functions are known for being robust measures that detect characteristics of the underlying distribution such as location, scale, bias, skewness and kurtosis [12, 51], and, in fact, statistical depth functions can characterise the underlying distribution [19]. Thus, a depth based two-sample test is potentially powerful against any alternative, or, at least, a wide variety. Furthermore, they have the potential of being more robust (depending on the depths being used).

Control charts based on depth-based ranks as in [52] have been discussed in [49] and applied to an aviation dataset in [14]. As more and more high-dimensional data are observed in practice, which can often be regarded as random curves, there has been an increased interest in statistical methodology for functional data in the last decades. For instance, two-sample tests for functional data have been considered by several authors. Horváth, Kokoszka and Reeder [43] proposed a test for differences in the mean function for $L^2(0,1)$-valued functional data. Dette, Kokot and Aue [24] developed a test for relevant differences in the mean function for observations in the Banach space of continuous functions. Munko et al. [58] proposed a bootstrap based multisample ANOVA procedure for functional data. Two-sample tests for differences in the covariance operators have been introduced by [5, 25, 36]. Wynne and Duncan [87] proposed a kernel two-sample test for functional data that is based on the maximum mean discrepancy. Jiménez-Gamero and Franco-Pereira [44] proposed a $k$-sample test that allows for $k$ to grow. Some of these papers also consider change point statistics as generalisations of their two-sample tests. Further tests and related estimators for functional data have been developed by [4, 6, 76] for changes in the mean and by [77, 81] for changes in the covariance.

### 1.1. Main contributions

Below, we summarise the main contributions made in this paper:

**Symmetrisation of a class of two-sample tests** A very important property when testing for the inequality of the underlying distributions of two samples is the symmetry of the procedure with respect to the two samples, i.e. the outcome of the test should not depend on the choice of label given to the two samples.

- In Section 2.2, we prove that the originally proposed Liu-Singh-test does not have this property to the extent that the power of the test significantly depends on the labelling with only one labelling having asymptotic power one.
- In Section 2.4 we derive asymptotic results for the joint tuple obtained from both orderings under the null hypothesis, which allows us to propose symmetric, asymptotic size $\alpha$ tests based on this joint tuple.

**Improvements on the size and power of the tests** Due to the particular construction of the rejection regions of the proposed test statistics, we tend to get a better small sample size behaviour. This is achieved by taking more information from the proofs of the null asymptotics into account. Additionally, as indicated by an illustrative example in Section 2.2, the proposed symmetrisation leads to a potentially large increase in power, in particular in those situations where the original statistic behaved poorly. A different (symmetrised) modification of the original statistic was already mentioned in the original paper by Liu and Singh [52] and later investigated for the two- and multisample case [15, 53] as well as for the change point situation [66, 68, 67]. But this modification was constructed to detect differences in the scale of the distributions and, unlike our proposal, has no or only very little power against location shifts (see, for instance, [52], Section 4).

**Extensions to functional spaces** In Section 3, we generalise the depth-based test for the nonparametric two-sample problem as proposed by Liu and Singh [52] to functional data. This is done by using depth functions for Hilbert and Banach space valued observations, see also [61]. Unlike for some other two-sample tests for functional data, no dimension reduction techniques and no or only weak moment assumptions are required depending on the depth used. The proposed tests – after the modifications already mentioned above – can simultaneously detect deviations in the mean as well as in the covariance operator and even have some power against different alternatives such as differences in higher-order moments or different shapes.

**Corrections of previous results** Finally, we fill a gap in the proof of the null asymptotics as provided by Zuo and He [93] and present a counterexample to their asymptotics result under the alternative (see Remark 2.4).

### *1.2. Organisation of the material*

In Section 2, different variations of the Liu-Singh-statistic are considered and their asymptotics under the null hypothesis are investigated under reasonable assumptions. The problem with the asymmetry is addressed and an illustrative counterexample to [93], Theorem 1 under the alternative, is presented. Section 3 is dedicated to verify the assumptions on the data depth functions used with our proposed statistics such that the asymptotics results from Section 2 hold true.

Section 4 contains a summary of our simulation study as well as an analysis of temperature curves of ocean drifters in order to distinguish between periods of El Niño and La Niña (in different intensities). For the proofs, an extensive simulation study and further analysis of the real dataset, we refer the reader to the appendix.

## 2. Variations of Liu-Singh-statistics

Wilcoxon's rank sum statistic measures whether the ranks corresponding to one sample are too small or too large compared to what could be expected if the distribution of both samples were the same. Statistical depth functions provide a centre-outward ordering and as such allow to measure how deep one sample lies within the other. This motivated Liu and Singh [52] to use statistical depth functions to generalise Wilcoxon's rank sum test, as described in the following.

### 2.1. The Liu-Singh-statistic and an initial modification

Consider two samples $X_1, ..., X_m$ (referred to as $X$-sample) and $Y_1, ..., Y_n$ ($Y$-sample) such that $X; X_1, ..., X_m \sim P$ iid and independent of $Y; Y_1, ..., Y_n \sim Q$ iid both with values in $\Omega$, where we want to test the null hypothesis $P = Q$ against the alternative $P \neq Q$. For a depth function $D : \Omega \times \mathcal{P}_\Omega \to \mathbb{R}_{\geqslant 0}$, where $\mathcal{P}_\Omega$ denotes the collection of all probability measures on $\Omega$ with respect to an appropriate $\sigma$-algebra, the expression

$$\widetilde{R}(y, P) = \int \mathbb{1}_{\{D(x,P) \leqslant D(y,P)\}} \, \mathrm{d}P(x). \tag{1}$$

can be regarded as a measure for the relative deepness of a point $y$ with respect to the probability measure $P$. The expected value is given by $\mathbb{E}\widetilde{R}(Y, P) = \mathbb{P}(D(X, P) \leqslant D(Y, P))$, which equals 1/2 under the null hypothesis of $P = Q$ if $D(X, P)$ is continuous (as under Assumption $\mathrm{Asm}_1$ below in Subsection 2.3.1). With $\widehat{P}_m$ denoting the empirical measure based on $X_1, ..., X_m$, the corresponding empirical version

$$\widetilde{R}(Y_j, \widehat{P}_m) = \frac{1}{m} \sum_{i=1}^{m} \mathbb{1}_{\{D(X_i, \widehat{P}_m) \leqslant D(Y_j, \widehat{P}_m)\}}$$

can be seen as a generalisation of the rank of observation $Y_j$ with respect to the $X$-sample. Unlike for $\widetilde{R}(Y, P)$, the expected value of the empirical version is not equal to 1/2, not even for continuous $D(X, P)$, because equality of $D(X, P)$ and $D(Y, Q)$ can occur with positive probability under the empirical measure. For this bias to be asymptotically negligible we need an additional regularity condition, which can be difficult to evaluate for specific depth functions (see Remark 2.2 below in Subsection 2.3.1). Therefore, we propose the following modification

$$R(y, P) = \int \mathbb{1}_{\{D(x,P) < D(y,P)\}} \, \mathrm{d}P(x) + \frac{1}{2} \int \mathbb{1}_{\{D(x,P) = D(y,P)\}} \, \mathrm{d}P(x), \tag{2}$$

$$R(Y_j, \widehat{P}_m) = \frac{1}{m}\sum_{i=1}^{m} \mathbb{1}_{\{D(X_i,\widehat{P}_m)<D(Y_j,\widehat{P}_m)\}} + \frac{1}{2}\frac{1}{m}\sum_{i=1}^{m} \mathbb{1}_{\{D(X_i,\widehat{P}_m)=D(Y_j,\widehat{P}_m)\}}. \tag{3}$$

For the theoretical quantity it clearly holds $R(y,P) = \widetilde{R}(y,P)$ a.e. for continuous distributions $P$.

Similarly to the Wilcoxon statistic we sum up all these *generalised ranks* $R(Y_j, \widehat{P}_m)$ leading to the following version of the Liu-Singh-statistic, to which we refer as $\mathcal{LS}$-statistic in the rest of the paper.

$$\mathcal{LS}(\widehat{P}_m, \widehat{Q}_n) = \int R(y, \widehat{P}_m)\, \mathrm{d}\widehat{Q}_n(y) = \frac{1}{n}\sum_{j=1}^{n} R(Y_j, \widehat{P}_m),$$

where $\widehat{Q}_n$ is the empirical measure based on $Y_1, ..., Y_n$. This, in turn, is the empirical version of the expected value of $R(y,P)$ with respect to the underlying distribution of the $Y$-sample $Q$, given by

$$\begin{aligned} \mathcal{LS}(P,Q) &= \int R(y,P)\, \mathrm{d}Q(y) \\ &= \mathbb{P}(D(X,P) < D(Y,P)) + \frac{1}{2}\mathbb{P}(D(X,P) = D(Y,P)). \end{aligned} \tag{4}$$

Roughly speaking, this measure-based index indicates whether or not the random variable $Y$ has a similar deepness with respect to $P$, as $X$ has. Or, in other words, it measures the *outlyingness* of the distribution $Q$ with respect to the distribution $P$ in terms of the depth. Thus, both large and small values of $\mathcal{LS}(P,Q)$ indicate a departure from the null hypothesis of equality of the distributions.

The original Liu-Singh-statistic introduced in [52] was constructed in the same way but using $\widetilde{R}$ instead of $R$. In the following sections, we refer to that statistic as $\widetilde{\mathcal{LS}}$-statistic.

## 2.2. Drawbacks of the $\mathcal{LS}$-statistics: an illustrative example

One major drawback of the $\mathcal{LS}$-statistic is the fact that it is not symmetric in the two samples in the sense that exchanging the labels of the samples may lead to a different test decision. Any two-sample test checking for inequality of the two underlying distributions should have such a symmetry property to guarantee a meaningful statistical interpretation. Otherwise, two practitioners may come to very different conclusions when applying the same test statistics to the same data but with exchanged labels.

Furthermore, we illustrate in this section that the combined information from the two $\mathcal{LS}$-statistics obtained by both labellings, i.e. the $\mathcal{LS}$-tuple $(\mathcal{LS}(\widehat{P}_m, \widehat{Q}_n), \mathcal{LS}(\widehat{Q}_n, \widehat{P}_m))$, is much greater than each of the two statistics separately. In addition, this is later confirmed by our simulation results in Sec-tion C.1. See also Figure 2, in Subsection 2.4, to get an impression of the rejection regions (and corresponding empirical size) for different statistics based on the $\mathcal{LS}$-tuple.

To discuss this point in detail, we use the following illustrative example: Let one sample follow a $U(0, 0.5)$ distribution and the other a $U(0, 1)$ distribution.

Computing the depth of a $U(0, 0.5)$ sample with respect to a $U(0, 1)$ sample leads to very similar values to those obtained when computing the depth of a $U(0, 1)$ sample with respect to a $U(0, 1)$ sample: In the former case, all mass lies in the same direction from the centre point as opposed to equally split to both directions in the latter. However, this direction is not taken into account by depth functions which only constitute a centre-outward ordering. Therefore, very little power can be expected from a test that is based on this notion. Indeed, this is confirmed by Theorem 2.1 (a) below for the Tukey depth, which is defined as

$$D_T(y, P) = \min(P((-\infty, y]), P([y, \infty))) \tag{5}$$

in the univariate case, cf. (15) in Section 3 below. The original rank definition in one dimension takes both the centre-outward ordering and the direction from the centre into account, so the Wilcoxon test does not suffer from this problem.

On the other hand, with exchanged labels the depths of the $U(0, 1)$ sample compared to the $U(0, 0.5)$ samples is considered, such that about half of the samples of the $U(0, 1)$ sample will have depth value 0. This is much more than can be expected under the null hypothesis, resulting in a good power behaviour. This is confirmed by Theorem 2.1 (b) for the Tukey depth.

**Theorem 2.1.** *Consider two independent samples of iid random variables $X, X_1, ..., X_m \sim U(0, 1)$ with empirical distribution function (EDF) $\widehat{P}_m$ and $Y, Y_1, ..., Y_n \sim U(0, 1/2)$ with EDF $\widehat{Q}_n$ and $m/(m + n) \to \tau \in (0, 1)$. Then, it holds for the $\mathcal{LS}$-statistic with respect to the Tukey depth, $\mathcal{LS}_T$, that*

$$\text{(a)} \quad \sqrt{\frac{12 \cdot m \cdot n}{m + n}} \left( \mathcal{LS}_T(\widehat{P}_m, \widehat{Q}_n) - \frac{1}{2} \right) \xrightarrow{\mathcal{D}} N(0, 5 - 4\tau),$$

$$\text{(b)} \quad \sqrt{\frac{12 \cdot m \cdot n}{m + n}} \left( \mathcal{LS}_T(\widehat{Q}_n, \widehat{P}_m) - \frac{1}{4} \right) \xrightarrow{\mathcal{D}} N\left(0, \frac{5}{4} - \tau\right).$$

*Both (a) and (b) remain true when substituting $\mathcal{LS}_T$ by $\widetilde{\mathcal{LS}}_T$, the original statistic with respect to the Tukey depth.*

The limit distributions under the null hypothesis (i.e. when both samples follow a $U(0, 1)$ distribution) is given by (see Theorem 2.3 below)

$$\sqrt{\frac{12 \cdot m \cdot n}{m + n}} \left( \mathcal{LS}_T(\widehat{P}_m, \widehat{Q}_n) - \frac{1}{2} \right) \xrightarrow{\mathcal{D}} N(0, 1) \text{ and, equally,}$$

$$\sqrt{\frac{12 \cdot m \cdot n}{m + n}} \left( \mathcal{LS}_T(\widehat{Q}_n, \widehat{P}_m) - \frac{1}{2} \right) \xrightarrow{\mathcal{D}} N(0, 1).$$

Consequently, for the labelling as in (a) the limit distribution under this alternative coincedes with the null distribution except for a larger variance. Therefore, the power of the corresponding test will remain bounded away from one due to the same expectation but with a power larger than the level due to the larger variance. On the other hand, for the ordering as in (b) the two limit distributions (under the alternative and the null hypothesis) contract around different expectations (1/2 under the null hypothesis and 1/4 under the alternative) so that the test will have asymptotic power one.

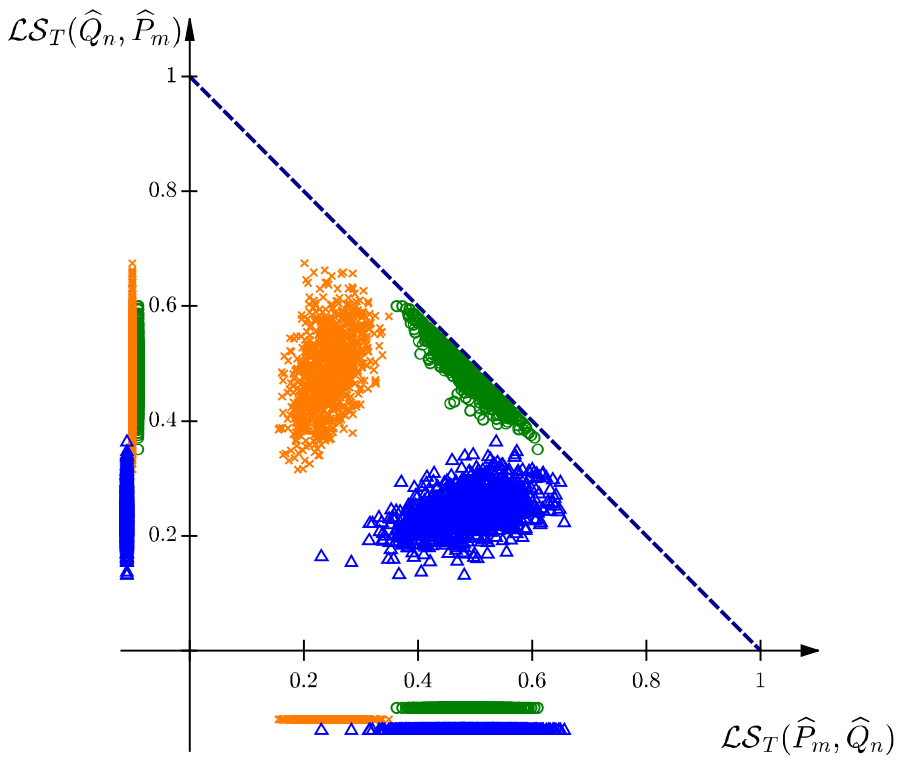

FIG 1. *For 1000 trials, we plotted for $m, n = 100$ the $\mathcal{LS}$-tuple $(\mathcal{LS}(\widehat{P}_m, \widehat{Q}_n), \mathcal{LS}(\widehat{Q}_n, \widehat{P}_m))$ and the corresponding projections. Green circles: $P, Q = U(0, 1)$. Orange crosses: $P = U(0, 0.5), Q = U(0, 1)$. Blue triangles: $P = U(0, 1), Q = U(0, 0.5)$.*

A corresponding empirical illustration, based on the univariate Tukey depth, can be found in Figure 1 for $m = n = 100$. In the figure, we have a scatter plot, with 1000 repetitions, of the $\mathcal{LS}$-tuple $(\mathcal{LS}_T(\widehat{P}_m, \widehat{Q}_n), \mathcal{LS}_T(\widehat{Q}_n, \widehat{P}_m))$ with $\mathcal{LS}_T(\widehat{P}_m, \widehat{Q}_n)$ corresponding to the $x$-axis and $\mathcal{LS}_T(\widehat{Q}_n, \widehat{P}_m)$ to the $y$-axis. In addition, $\mathcal{LS}_T(\widehat{P}_m, \widehat{Q}_n)$ is represented below the $x$-axis and $\mathcal{LS}_T(\widehat{Q}_n, \widehat{P}_m)$ in the left of the $y$-axis. In the figure, the green circles correspond to the null hypothesis, being all random variables drawn independently from $U(0, 1)$ distributions. Meanwhile, the orange crosses correspond to the scenario in Theorem 2.1, with the $X$ random variables being drawn independently from a $U(0, 0.5)$ distribution. The dark blue dash-dotted secondary diagonal in Figure 1 divides the unit square into two triangles and the $\mathcal{LS}$-tuple usually takes values in the below triangle, see Section C.4 for details. As suggested by the above considerations regarding Theorem 2.1 the orange and green points of the projection onto the $y$-axis, $\mathcal{LS}(\widehat{Q}_n, \widehat{P}_m)$, are centred around 1/2 with the variability of the orange crosses larger than that of the green circles. On the other hand, when the labels are switched as in $\mathcal{LS}(\widehat{P}_m, \widehat{Q}_n)$ (projection onto the $x$-axis), then the orange crosses are well separated from the green circles (as indicated by the different asymptotic expectations). The blue triangles in Figure 1 are a result of switching the roles of $U(0, 0.5)$ and $U(0, 1)$, where a behaviour in an analogous fashion to the orange crosses can be observed. The clearly separated point clouds in Figure 1 indicate that tests based on both entries of the $\mathcal{LS}$-tuple can be constructed to have higher power and at the same time not suffer from asymmetry problems. From our empirical observations, we conjecture that the location of the points of the $\mathcal{LS}$-tuple in the unit square gives an indication on the type of difference between the underlying distributions (see Section C.4).

A different way to symmetrise the $\mathcal{LS}$-test has already been pointed out by Liu and Singh [52], and later discussed by [15, 53] but has less power against the important alternative of mean differences (compared to a joint consideration of both orderings based on the scatter plot above). The reason is that their proposal measures the deepness of the $Y$-sample not in the $X$-sample but in the joint sample $X_1, ..., X_m, Y_1, ..., Y_n$. For example, by analogous considerations to above, a sample of $U(0, 0.5)$ random variables **joint** with a $U(0.5, 1)$ sample (of the same size) will have similar behaviour as if both were $U(0, 1)$, thus a corresponding test suffers from power problems. Note that this joint samples depth approach follows the original Wilcoxon statistic but its power issue regarding mean differences arises from the fact that statistical depths only yield a centre-outward ordering with no distinction of direction from the centre, whereas ranks – that are limited to the univariate case – make use of that information.

### 2.3. *Joint asymptotic behaviour of the $\mathcal{LS}$-tuple under the null hypothesis*

In the previous section, we have illustrated that a combined consideration of the $\mathcal{LS}$-statistics with both possible orderings is symmetric by definition and contains much more information than each of the two projections in the axes separately. Indeed, due to the better separation between null hypothesis and alternatives, the joint consideration offers a large potential for a power gain.

Usually, non-rejection regions for a test based on a bi-variate statistic such as the $\mathcal{LS}$-tuple are constructed as confidence regions of the asymptotic distribution under the null, i.e. regions that contain $1 - \alpha$ of its probability mass. In the following, we will shed light onto the joint asymptotic behaviour of the $\mathcal{LS}$-tuple. We will obtain the following key-observations, which will eventually motivate our proposed non-rejection regions:

- Theorem 2.3 shows that the joint asymptotic null distribution is degenerate, i.e. perfectly negatively correlated asymptotically, where the proof shows that this is due to the same leading term with different signs.
- This leading term is best captured by the properly scaled difference of the two statistics (henceforth called $\mathcal{LS}$-difference statistic). As a consequence, the latter has a much better empirical size behaviour than each of the original $\mathcal{LS}$-tests, see Figure 2.
- Indeed a suitable rotation transforms the $\mathcal{LS}$-tuple into a tuple of the $\mathcal{LS}$-difference statistic in addition to an $\mathcal{LS}$-sum statistic, which is defined as the properly scaled sum of the original $\mathcal{LS}$-statistics. Corollary 2.5 shows that the latter contracts faster than the $\mathcal{LS}$-difference statistic.
- For the special case of the Tukey depth in the one-dimensional Euclidean space, we prove the asymptotic independence of the $\mathcal{LS}$-difference and $\mathcal{LS}$-sum statistic in Theorem 2.6 also providing the joint asymptotic null distribution of the $\mathcal{LS}$-tuple.
- Some results for the simplicial depth in the one-dimensional Euclidean space suggest that the limit of the $\mathcal{LS}$-sum statistic and possibly also the contraction rate may depend on the underlying depth function, see Theorem 2.7 and Remark 2.9. In consequence, obtaining confidence regions from the joint limit distribution may not always be feasible.
- Finally, Section 2.4 discusses how to choose non-rejection regions in light of the above findings.

For the asymptotic behaviour of the $\widetilde{\mathcal{LS}}$-statistic, i.e. the projections of the $\widetilde{\mathcal{LS}}$-tuple, Liu and Singh [52] already presented some first results for several depths. Later, Rousson [71] presented ideas for a more general proof, before Zuo and He [93] gave a general proof based on regularity conditions on the depths and distributions. We will follow their approach in understanding the joint behaviour of the $\mathcal{LS}$-tuple, as such an approach has the advantage that additional depth functions on a variety of probability spaces can be used. In particular, we focus on applying the $\mathcal{LS}$-test for functional data in Section 4. In addition, we fill a gap in the proof in Zuo and He [93] (see Remark 2.4), as a small adaptation of the original regularity conditions is required.

### *2.3.1. Assumptions on the depths*

We state here the necessary regularity conditions on the depths and distributions, giving the corresponding asymptotic results later in the section. We will verify these regularity conditions for several important depths in Section 3.

Let $D : \Omega \times \mathcal{P}_\Omega \to \mathbb{R}_{\geqslant 0}$ be a depth function. In Assumption $\text{Asm}_1$ we require the CDF of $D(X, P)$ to be $\beta$-Hölder-continuous. Essentially, this requires the random variable $X$ to have a continuous distribution, which is a standard assumption even for univariate rank-based procedures, in addition to the depth being sufficiently regular in its first argument.

**Assumption $\mathbf{Asm_1}$.** $\mathbb{P}(D(X,P) \in [y_1, y_2]) \leqslant C|y_2 - y_1|^\beta$ *for positive constants* $C, \beta$ *and any* $y_1, y_2 \in \mathbb{R}_{\geqslant 0}$ *for some* $1/2 < \beta \leqslant 1$.

In Assumption $\text{Asm}_2$ we require sufficiently fast uniform contraction rates in the $L^{2\beta}$-norm for the empirical depth processes.

**Assumption $\mathbf{Asm_2}$.** $\mathbb{E}(\sup_{x\in\Omega} |D(x, \widehat{P}_m) - D(x,P)|^{2\beta}) = O(m^{-\beta})$ *with* $1/2 < \beta \leqslant 1$ *as in* $Asm_1$.

Assumption $\text{Asm}_3$ gives a deterministic upper bound between the usual empirical depth function and the empirical depth function where one element was taken out of the sample.

**Assumption $\mathbf{Asm_3}$.** *There exist a* ***deterministic*** *constant* $\mathrm{C}_{\text{det}}$ *and an index* $m_0$ *such that for every* $m \geqslant m_0$*:*

$$\sup_{x\in\Omega} |D(x, \widehat{P}_m) - D(x, \widehat{P}_m^{-\{1\}})| \leqslant \frac{\mathrm{C}_{\text{det}}}{m} \qquad \textit{almost surely,}$$

*where $\widehat{P}_m^{-I}$ for some $I \subset \{1,...,n\}$ is the empirical probability measure with respect to $\{X_i : i \in \{1,...,n\} \setminus I\}$.*

If we take $M_n$ elements out of the set (instead of just 1), then we immediately get an upper bound of $M_n \cdot \mathrm{C_{det}}/m$ in the above assumption. In particular, we get the same bound (with a different deterministic constant) if we take only a bounded number of elements $M_n \leqslant M$ out. In the proofs we make use of

$$\sup_{x\in\Omega} |D(x, \widehat{P}_m) - D(x, \widehat{P}_m^{-\{1,2\}})| \leqslant \frac{2\,\mathrm{C_{det}}}{m} \qquad \text{almost surely.} \tag{6}$$

In a sense, $\mathrm{Asm}_4$ is the empirical analogue of $\mathrm{Asm}_1$, but making use of the $L^2$-norm.

**Assumption $\mathrm{Asm}_4$.** *Let $X$ be an independent copy of $X_1, ..., X_m$. Then, for any constant $C > 0$, it holds*

$$\mathbb{E}\left(\left[\mathbb{P}\left(|D(X_1, \widehat{P}_m^{-\{1\}}) - D(X, \widehat{P}_m^{-\{1\}})| \leqslant \frac{C}{m}\,\middle|\, X_2, ..., X_m\right)\right]^2\right) = O(m^{-2\beta})$$

*with $1/2 < \beta \leqslant 1$ as in $\mathrm{Asm}_1$.*

In particular, $\mathrm{Asm}_4$ implies, for any $C > 0$,

$$\mathbb{E}\left(\left[\mathbb{P}\left(|D(X_1, \widehat{P}_m^{-\{1,2\}}) - D(X, \widehat{P}_m^{-\{1,2\}})| \leqslant \frac{C}{m}\,\middle|\, X_3, ..., X_m\right)\right]^2\right) = O(m^{-2\beta}). \tag{7}$$

The proofs of Theorem 2.3 and Corollary 2.10 go through for $0 < \beta \leqslant 1/2$ as long as we get the bound $o(m^{-1})$ in $\mathrm{Asm}_4$. However, this seems unrealistic in view of the $\beta$-Hölder-continuity in $\mathrm{Asm}_1$.

The above assumptions related to the Hölder-continuity provide a uniform framework – but while they are sufficient, they are not necessary: Theorem 2.7 shows that Theorem 2.3 and consequently Corollary 2.10 also holds for the univariate simplicial depth which is only Hölder with $\beta = 1/2$ and thus is not covered by the above theory. Indeed, our proof in Section A.1 shows that the dominating term of the statistic as given by $A_1$ in (20) is asymptotically normal as long as $D(X_1, P)$ is continuous, i.e. in particular under $\mathrm{Asm}_1$ for any $\beta > 0$ (see Lemma A.1). The main difficulty, where the assumptions really come into play, is the derivation of the rate of convergence of the remaining term $A_2$ (see (20)) in Lemma A.3.

**Remark 2.2.** *If we use the original generalised ranks as in (1) instead of (2), then we need additionally that (see Remark A.5)*

$$\mathbb{E}\left(\left[\sum_{i=1}^{B_m}(\mathbb{P}(D(X, \widehat{P}_m) = c_i | X_1, ..., X_m))^2\right]^2\right) = o(m^{-1}), \tag{8}$$

*where $X$ is an independent copy of $X_1, ..., X_m$ and $c_1, ..., c_{B_m}$ are the only values for which $\mathbb{P}(D(X, \widehat{P}_m) = c_i | X_1, ..., X_m) > 0$ holds true.*

### 2.3.2. *Joint asymptotic behaviour of the $\mathcal{LS}$-tuple*

The following theorem gives the joint limit of the $\mathcal{LS}$-tuple $(\mathcal{LS}(\widehat{P}_m, \widehat{Q}_n), \mathcal{LS}(\widehat{Q}_n, \widehat{P}_m))$. The localisation constant of $1/2$ stems from the fact that under the null hypothesis and if $D(X_1, P)$ is continuous (as under $\text{Asm}_1$), it holds $R(Y_1, P) \sim U(0,1)$, therefore $\mathcal{LS}(P,P) = \frac{1}{2}$.

**Theorem 2.3.** *Let $Asm_1$ - $Asm_4$ hold true and, for $m, n \to \infty$, let*

$$\frac{m}{m+n} \to \tau \qquad \text{for some } 0 < \tau < 1. \tag{9}$$

*Then, under the null hypothesis of $P = Q$ and as $m \to \infty$, we have that*

$$\sqrt{\frac{12 \cdot m \cdot n}{m+n}} \left[ \begin{pmatrix} \mathcal{LS}(\widehat{P}_m, \widehat{Q}_n) \\ \mathcal{LS}(\widehat{Q}_n, \widehat{P}_m) \end{pmatrix} - \frac{1}{2} \right] \xrightarrow{\mathcal{D}} N\left( \begin{pmatrix} 0 \\ 0 \end{pmatrix}, \begin{pmatrix} 1 & -1 \\ -1 & 1 \end{pmatrix} \right). \tag{10}$$

**Remark 2.4.** *The main result of Zuo and He [93] – this is [93], Theorem 1 – is actually incorrect under the alternative hypothesis. The limit distribution of*

$$\sqrt{\frac{12 \cdot m \cdot n}{m+n}} \left( \widetilde{\mathcal{LS}}_T(\widehat{P}_m, \widehat{Q}_n) - \frac{1}{2} \right)$$

*is claimed there to be standard normal (see also (32), (33)). However, by Theorem 2.1 (a) the true asymptotic variance is given by $5 - 4\tau$, which is strictly greater than 1, and thus, we have a counterexample for their theorem. This theoretic result is illustrated in Figure 1. Zuo and He [93] omit the proofs of points (i) and (ii) in [93], Lemma 1. They only prove point (iii). While (i) can be proven along the lines of (iii) (see the proof of Lemma A.1), this is not true for (ii). The reason is that the independence argument used in the proof of (iii) does not hold for (ii). This concerns the term $A_{2,2}(\widehat{P}_m)$ (respectively $\widetilde{A}_{2,2}(\widehat{P}_m)$) in the proof of Lemma A.3 (see Remark A.4 for a detailed explanation of the problem). In order to fill this gap, we need some slightly different regularity conditions (cf. $Asm_2$–$Asm_4$) compared to those in the paper by Zuo and He [93]. As [93], Lemma 1 is used to prove [93], Theorem 1, we have that, under the null hypothesis, our Theorem 2.3 and the results that follow from it also require the above mentioned stronger assumptions.*

The limit in Theorem 2.3 is degenerate in the sense that the two statistics are perfectly negatively correlated asymptotically under the null hypothesis. Consequently, under the null hypothesis both statistics are asymptotically symmetric in both arguments, while this is not the case under alternatives (see Theorem 2.1 and Figure 1). This is also the reason why the power behaviour of the $\mathcal{LS}$-statistic depends crucially on the choice of labelling, as illustrated in Section 2.2. In Section 2.4.1 we propose some symmetric first-order variations based on the $\mathcal{LS}$-tuple, whose limit distribution can be obtained based only on Theorem 2.3 above. However, these statistics suffer from size or power problems because they do not take the full geometry of the $\mathcal{LS}$-tuple into account.

Indeed, the proof of Theorem 2.3 shows that $\mathcal{LS}(\widehat{P}_m,\widehat{Q}_n)$ and $\mathcal{LS}(\widehat{Q}_n,\widehat{P}_m)$ have the same symmetric leading term ($A_1$ in (20)) with different (asymptotically negligible) remainder terms ($A_2$ in (20)). This leading term is best captured by the $\mathcal{LS}$-difference statistic $(\mathcal{LS}(\widehat{P}_m,\widehat{Q}_n)-\mathcal{LS}(\widehat{Q}_n,\widehat{P}_m))/\sqrt{2}$ (see also Corollary 2.10 (a) below, where the different variance stems from the fact that the $\mathcal{LS}$-difference statistic contains the leading term $A_1$ in (20) twice).

This can also be thought of as rotation of the coordinate system, with the dark blue dash-dotted secondary diagonal and red dashed diagonal in Figure 2 below as new axes and the factor $1/\sqrt{2}$ corresponding to the rotation matrix. The projection of the $\mathcal{LS}$-tuple on the red dashed main diagonal, essentially given by the sum of the two $\mathcal{LS}$-statistics, corresponds to the remainder terms ($A_2$ in (20)) and contracts faster than the main term (i.e. it is asymptotically negligible) as shown by the following Corollary (see also Lemma A.3). In Figure 2, this negligibility is the reason that the variability of the $\mathcal{LS}$-tuple under the null hypothesis is greater along the dark blue dash-dotted secondary diagonal than along the red dashed diagonal.

Using not only the information about the limit with respect to the main term as reflected by the $\mathcal{LS}$-difference statistic, but also the asymptotic information about the complementing $\mathcal{LS}$-sum statistic given by $\mathcal{LS}(\widehat{P}_m,\widehat{Q}_n)+\mathcal{LS}(\widehat{Q}_n,\widehat{P}_m)$, can further improve upon the size-power behaviour of the tests. Such symmetric second-order tests are discussed in Section 2.4.2. Their construction is motivated by and depends crucially on the following mathematical results.

**Corollary 2.5.** *Under Asm1 - Asm4 and* (9)*, it holds*

$$\mathcal{LS}(\widehat{P}_m,\widehat{Q}_n)+\mathcal{LS}(\widehat{Q}_n,\widehat{P}_m)-1=O_P(\delta_{m,n}),$$

*where*

$$\delta_{m,n}=\begin{cases}\left(\frac{m+n}{mn}\right)^{\beta} & \text{for } \frac{1}{2}<\beta\leqslant\frac{2}{3},\\ \left(\frac{m+n}{mn}\right)^{\frac{2+\beta}{4}} & \text{for } \frac{2}{3}\leqslant\beta\leqslant 1.\end{cases}$$

*In particular, under Lipschitz-continuity the rate $\delta_{m,n}=\left(\frac{m+n}{mn}\right)^{\frac{3}{4}}$ is achieved.*

The following theorem shows that the true convergence rate associated to a specific depth can be much faster than suggested by $\delta_{m,n}$ in Corollary 2.5.

**Theorem 2.6.** *Let $\mathcal{LS}_T(\cdot,\cdot)$ be the $\mathcal{LS}$-statistic with respect to the Tukey depth in the one-dimensional Euclidean space.*

*(a) Then, under the null hypothesis for any continuous probability measure $P$, it holds that*

$$-\frac{mn}{m+n}\cdot\Big(\mathcal{LS}_T(\widehat{P}_m,\widehat{Q}_n)+\mathcal{LS}_T(\widehat{Q}_n,\widehat{P}_m)-1\Big)\xrightarrow{\mathcal{D}}1+\frac{1}{2}\chi_1^2,$$

*where $\chi_1^2$ is a $\chi$-square-distribution with 1 degree of freedom. Furthermore, $\mathcal{LS}_T(\widehat{P}_m,\widehat{P}_n)+\mathcal{LS}_T(\widehat{P}_n,\widehat{P}_m)-1\leqslant 0$ almost surely.*

*(b)* $\mathcal{LS}_T(\widehat{P}_m, \widehat{Q}_n) + \mathcal{LS}_T(\widehat{Q}_n, \widehat{P}_m)$ *and* $\mathcal{LS}_T(\widehat{P}_m, \widehat{Q}_n) - \mathcal{LS}_T(\widehat{Q}_n, \widehat{P}_m)$ *are asymptotically independent statistics, in the sense that the joint limit distribution of the properly scaled and shifted versions is the product of their marginal limit distributions.*

*For the* $\widetilde{\mathcal{LS}}$*-statistics, one obtains the limit result of* $\frac{1}{2}\chi_1^2$ *in (a) instead and (b) holds analogously.*

For the two-sample location-scale problem of univariate data, there exist several further tests, e.g. the Lepage-test which is a combination of the Wilcoxon statistic with the Ansari-Bradley test [47] or the Cucconi-test which combines information of squared ranks and anti-ranks [18]. Further details can be found in [56]. An alternative test statistic can be obtained from Theorem 2.6 by making use of the joint distribution to obtain a non-rejection region that is as tight as possible. Such $\mathcal{LS}$-tests are neither based on the Lepage test (respectively the Wilcoxon test) nor on the Cucconi test. However, even in this simple univariate situation, the proof of the above theorem is rather involved so that one cannot expect such a result in more complicated situations.

Finally, we show that for the univariate simplicial depth, despite being only 1/2-Hölder-continuous, the remainder term also contracts faster than the leading term. This happens with a surprisingly fast rate. The conjecture in Remark 2.9 below suggests that, even with known joint distribution, it might be difficult to decide on an appropriate shape for the non-rejection regions in general.

**Theorem 2.7.** *Let* $\mathcal{LS}_S(\cdot,\cdot)$ *be the* $\mathcal{LS}$*-statistic with respect to the univariate simplicial depth in the U-statistic-representation with respect to closed simplices,*

$$D_S(x, \widehat{P}_m) = \binom{m}{2}^{-1} \sum_{1 \leqslant i < j \leqslant m} \mathbb{1}_{\{x \in [\min(X_i, X_j), \max(X_i, X_j)]\}},$$

*see* (16) *for details. Then, under the null hypothesis* $P = Q$ *and* (9) *for any continuous probability measure* $P$*, it holds*

$$\mathcal{LS}_S(\widehat{P}_m, \widehat{Q}_n) + \mathcal{LS}_S(\widehat{Q}_n, \widehat{P}_m) - 1 = O_P\left(\frac{\log(m+n)}{m+n}\right).$$

*Moreover, the limit in* (10) *holds for the* $\mathcal{LS}$*-tuple based on* $\mathcal{LS}_S$*, as well as on* $\widetilde{\mathcal{LS}}_S$*.*

**Remark 2.8.** *The proof of the above theorem is based on the asymptotic behaviour of the difference of the* $\mathcal{LS}$*-statistics with the respect to the univariate Tukey and simplicial depth: indeed,* (62) *shows that the difference between* $\widetilde{\mathcal{LS}}_S$ *and* $\widetilde{\mathcal{LS}}_T$ *is asymptotically negligible, while* (64) *shows the asymptotic negligibility between* $\widetilde{\mathcal{LS}}_S$ *and* $\mathcal{LS}_S$*.*

**Remark 2.9.** *Due to some preliminary considerations we conjecture that the following asymptotics holds in the situation of Theorem 2.7:*

$$- \frac{mn}{(m+n)(\log(\sqrt{2m+3}) + \log(\sqrt{2n+3}))} \cdot \left( \widetilde{\mathcal{LS}}_S(\widehat{P}_m, \widehat{P}_n) + \widetilde{\mathcal{LS}}_S(\widehat{P}_n, \widehat{P}_m) - 1 \right) \xrightarrow{P} \frac{1}{2}.$$

*The corresponding limit is degenerate and therefore highly non-standard and cannot be used to derive asymptotic confidence regions even if the joint limit with the $\mathcal{LS}$-difference statistic were known. Instead, one would need to subtract* $1/2$ *from the left-hand side (including the normalizing sequence) and then consider the limit distribution of that one.*

### 2.4. Symmetrisation based on the joint $\mathcal{LS}$-tuple

In this section, we discuss possibilities to construct symmetric tests based on the joint $\mathcal{LS}$-tuple. We distinguish between *first-order tests*, discussed in Section 2.4.1, and *second-order tests*, in Section 2.4.2. The limit of the first-order tests is a direct consequence of Theorem 2.3. They do not make any additional use of the information about the limit of the $\mathcal{LS}$-sum statistic. In consequence, these statistics suffer from unnecessary size or power issues, see Figure 2, where the power issues are indicated by unnecessary blind spots. In contrast, second-order tests additionally take the asymptotic information from the $\mathcal{LS}$-sum statistic into account, thus improving on the power while keeping a good size behaviour.

#### *2.4.1. First-order tests*

The null asymptotics of the original statistics given by the two projection statistics of the $\mathcal{LS}$-tuple, i.e. $\mathcal{LS}(\widehat{P}_m, \widehat{Q}_n)$ and $\mathcal{LS}(\widehat{Q}_n, \widehat{P}_m)$, is a direct consequence of Theorem 2.3.

The following corollary introduces two more statistics and gives the corresponding limit under the null hypothesis.

**Corollary 2.10.** *Let* $Asm_1$ *-* $Asm_4$ *and* (9) *hold. Then, under the null hypothesis of* $P = Q$*, it holds for* $m \to \infty$

$$(a) \quad \sqrt{\frac{6 \cdot m \cdot n}{m+n}} \cdot \left[ \frac{\mathcal{LS}(\widehat{P}_m, \widehat{Q}_n) - \mathcal{LS}(\widehat{Q}_n, \widehat{P}_m)}{\sqrt{2}} \right] \xrightarrow{\mathcal{D}} N(0,1),$$

$$(b) \quad \frac{12 \cdot m \cdot n}{m+n} \left[ \max \left( \left| \mathcal{LS}(\widehat{P}_m, \widehat{Q}_n) - \frac{1}{2} \right|, \left| \mathcal{LS}(\widehat{Q}_n, \widehat{P}_m) - \frac{1}{2} \right| \right) \right]^2 \xrightarrow{\mathcal{D}} \chi^2_1.$$

While both asymptotic statements in the corollary follow immediately from Theorem 2.3, the approximation of the $\mathcal{LS}$-difference statistics, in Corollary 2.10 (a), is quite good, while the approximation of the $\mathcal{LS}$-maximum statistic, in Corollary 2.10 (b), is not so good – resulting in a size problem for small sample sizes. The problem is illustrated in the empirical size shown in Figure 2 (d), corresponding to the $\mathcal{LS}$-maximum statistic, which is of the 11%, being far higher than the desired 5% that holds asymptotically.

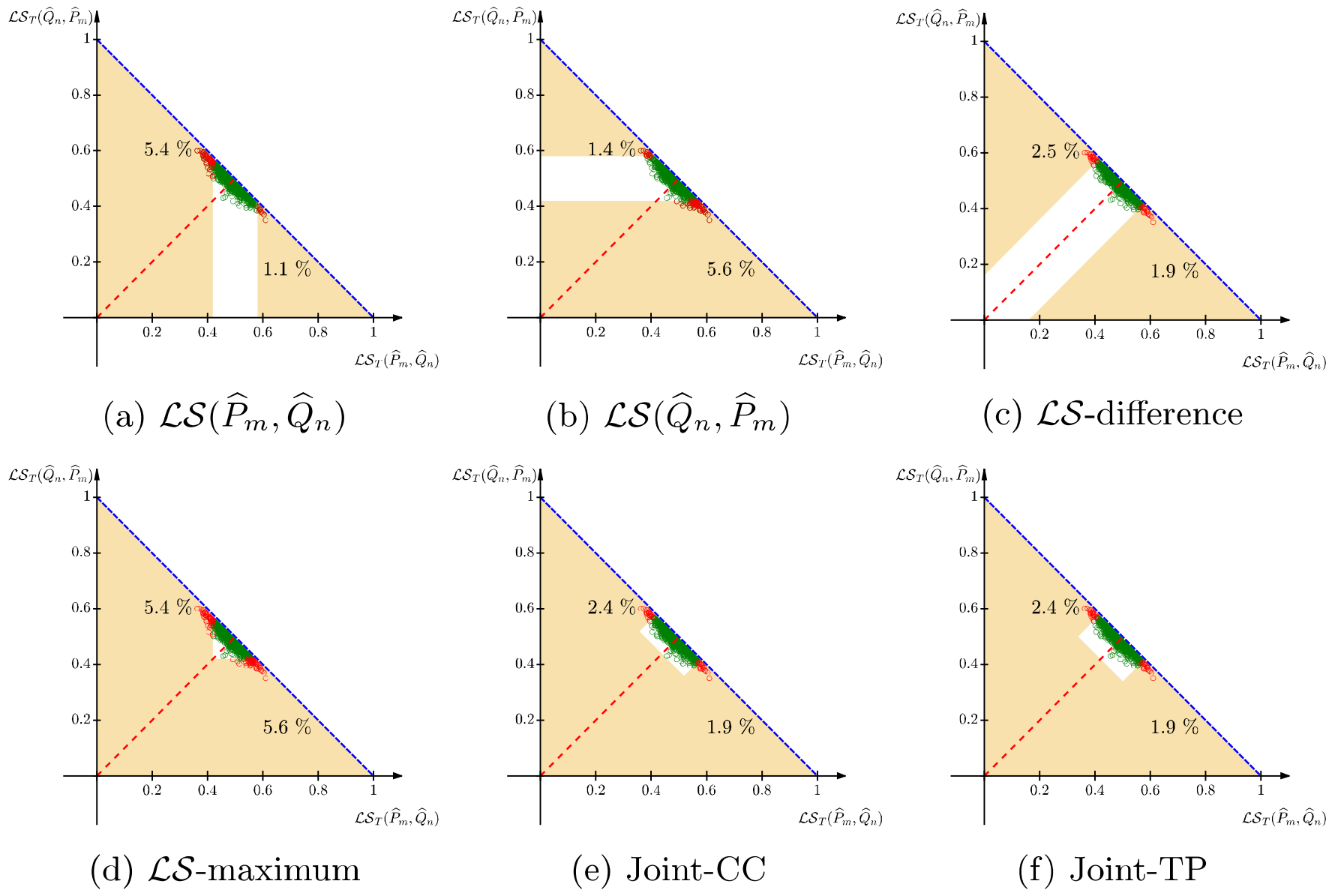

Fig 2. *Illustration of the blind spots: in each panel, the circles represent* 1000 *realizations of the $\mathcal{LS}$-tuple ($\mathcal{LS}(\widehat{P}_m, \widehat{Q}_n)$ in the horizontal axis and $\mathcal{LS}(\widehat{Q}_n, \widehat{P}_m)$ in the vertical axis) obtained with the univariate Tukey depth under the null hypothesis of $P$ and $Q$ following a $U(0,1)$ distribution ($m = n = 100$). The red circles (unlike the green ones) result in rejections made by the $\mathcal{LS}$-statistic corresponding to each panel. From left to right and top to bottom, the used $\mathcal{LS}$-statistics are: (a) $\mathcal{LS}(\widehat{P}_m, \widehat{Q}_n)$, (b) $\mathcal{LS}(\widehat{Q}_n, \widehat{P}_m)$, (c) $\mathcal{LS}$-difference, (d) $\mathcal{LS}$-maximum, (e) Joint-CC and (f) Joint-TP. In each panel, the shaded orange region is the rejection region of the corresponding variation of the $\mathcal{LS}$-statistic. Each rejection region is separated in two by a red dashed main diagonal. The numbers at each side of the red dashed main diagonal correspond to the rejection rate of that side.*

It is worth commenting that recently, Shi, Zhang and Fu [79] have run extensive simulation studies concerned with the empirical power, but not the empirical size, for the $\mathcal{LS}$-maximum statistic in addition to a $\mathcal{LS}$-ellipsoidal test with ellipsoid non-rejection regions based on a convex combination of the two projection tests; see Section C.3 in the for more details.

Figure 2 also contains the corresponding results for the $\mathcal{LS}$-projection statistics in Subfigures (a) and (b). Additionally, Figure 2 (c) corresponds to the $\mathcal{LS}$-difference statistic in Corollary 2.10 (a). Clearly, the empirical size is much better in Figure 2 (c) for the $\mathcal{LS}$-difference statistic than for the $\mathcal{LS}$-maximum statistic or any of the two $\mathcal{LS}$-statistics in Subfigures (a) or (b). Note that the empirical size shown in Figure 2 (c) is 4.4%, and thus much closer to the asymptotically given 5% than those in Figure 2 (a) and (b), which are 6.5% and 7%, respectively. Furthermore, the size problem of the two $\mathcal{LS}$-projection statistics clearly stems from the fact that the geometry of the $\mathcal{LS}$-tuple is not taken into account resulting also in an asymmetric rejection behaviour in the two disjoint rejection regions of the $\mathcal{LS}$-statistics. The $\mathcal{LS}$-maximum statistic in Figure 2 (d) has even greater size problems because it inherits the size behaviour from the 'bad' region of both $\mathcal{LS}$-statistics. On the other hand, it has clearly eliminated the *blind spots*, i.e. areas that do not belong to the rejection region despite the fact that $\mathcal{LS}$-tuples will only fall into these areas under alternatives. Figure 1 gives examples where the $\mathcal{LS}$ tuple does indeed fall into these *blind spots* under alternatives. Thus, the choice of rejection region as indicated in shaded orange in Figure 2 is truly relevant for the power of the test against different alternatives.

#### *2.4.2. Second-order tests: Joint $\mathcal{LS}$-tests*

The problems of the above first-order tests stem from the fact that they do not consider all available information about the asymptotic behaviour of the $\mathcal{LS}$-sum statistic under the null hypothesis. Taking this information into account, in this section, we propose second-order tests that can improve upon the power while keeping the good size behaviour of the $\mathcal{LS}$-difference statistic. Essentially, we aim at obtaining tests that

- are symmetric in both samples,
- have no systematic blind spots (making use of Corollary 2.5), and thus are powerful against more alternatives,
- make use of the asymptotic approximation in Corollary 2.10 (a), with better small-sample behaviour than e.g. the statistic in Corollary 2.10 (b),
- and have asymptotic size $\alpha$.

To understand their construction, the following key observations are important:

(a) The normal approximation of the $\mathcal{LS}$-difference statistic is very good even for small samples.
(b) We cannot control the behaviour of the $\mathcal{LS}$-sum statistic in general, but we know that it contracts faster than the $\mathcal{LS}$-difference statistic.
(c) We know that the $\mathcal{LS}$-difference and the $\mathcal{LS}$-sum statistic are asymptotically independent when making use of the Tukey depth in the one-dimensional Euclidean space, see Theorem 2.6. We conjecture that this remains true in more general situations.

Observations (a) and (c) motivate the use of rectangular non-rejection regions that are oriented as the diagonal lines in Figure 2.

More precisely, we use the $\mathcal{LS}$-tuple with a rectangular non-rejection region described by the following inequalities

$$\begin{aligned} (I) \quad & \frac{\left|\mathcal{LS}(\widehat{P}_m, \widehat{Q}_n) - \mathcal{LS}(\widehat{Q}_n, \widehat{P}_m)\right|}{\sqrt{2}} \leqslant \sqrt{\frac{m+n}{12 \cdot mn}}\, z_{1-\frac{\alpha}{2}}, \\ (II) \quad & \frac{\mathcal{LS}(\widehat{P}_m, \widehat{Q}_n) + \mathcal{LS}(\widehat{Q}_n, \widehat{P}_m)}{\sqrt{2}} \geqslant \sqrt{2}\left(\frac{1}{2} - \gamma_{m,n}\right), \end{aligned} \tag{11}$$

for some $\gamma_{m,n} \to 0$ discussed below. That is, the null hypothesis is rejected if at least one of those two inequalities is violated. For an illustration containing such types of rectangles, see Subfigures (e) and (f) of Figure 2. For the sake of completeness, note that $(I)$ and $(II)$ describe a region that is not restricted to the lower triangle below the dark blue dash-dotted secondary diagonal in Figure 2 and that with some particular depths and finite samples the $\mathcal{LS}$-tuple can take values above that secondary diagonal, see Section C.4.

In Figure 2, the term (without the absolute value) on the left-hand side of $(I)$ corresponds to the projection of the $\mathcal{LS}$-tuple on the dark blue dash-dotted secondary diagonal, essentially corresponding to the $\mathcal{LS}$-difference statistic. Note that here the rejection rule is obtained from the corresponding normal approximation, see observation (a) above. On the other hand, the term on the left-hand side of $(II)$ corresponds to the projection of the $\mathcal{LS}$-tuple on the red dashed main diagonal, essentially corresponding to the $\mathcal{LS}$-sum statistic. Here, the rejection rule is obtained from observation (b). As a consequence of the construction of $(I)$ and $(II)$, the test holds the level $\alpha$ asymptotically by Corollaries 2.10 (a) as long as $\gamma_{m,n}$ contracts slower to 0 than the contraction rate of $\mathcal{LS}(\widehat{P}_m, \widehat{Q}_n) + \mathcal{LS}(\widehat{Q}_n, \widehat{P}_m) - 1$ (corresponding to the remainder term) as in Corollary 2.5. If faster convergence rates are available, as in Theorems 2.6 or 2.7, then it is sufficient for $\gamma_{m,n}$ to contract slower than these.

A conservative choice for the contraction in $(II)$ that works for all depth functions, where the remainder term is indeed negligible compared to the leading term, is given by

$$\gamma_{1,m,n} = \sqrt{\chi^2_{1,0.95} \frac{m+n}{12mn}}, \tag{12}$$

which contracts at the same rate as the leading term and thus at the same rate as the $\mathcal{LS}$-difference statistic. The constant reflects the chosen quantiles for the normal approximation leading to a non-rejection rectangle that extends to the tip of the non-rejection triangle of the $\mathcal{LS}$-maximum statistic as in Figure 2 (d), which is a conservative choice. We call the corresponding test Joint-TP test (for **T**riangle-**P**oint). This is illustrated in Figure 2 (f).

Alternatively, if faster contraction rates $\delta_{m,n}$ for the $\mathcal{LS}$-sum statistic $\mathcal{LS}(\widehat{P}_m, \widehat{Q}_n) + \mathcal{LS}(\widehat{Q}_n, \widehat{P}_m) - 1$ (corresponding to the remainder term) are known, e.g. by Corollary 2.5, then they can be used to further decrease the rectangular non-rejection region. Any such rectangle should not extend beyond the one corresponding to $\gamma_{1,m,n}$ for any value of $m$ and $n$ and at the same time should only start to significantly differ for suitably large sample sizes $m$ and $n$. Both can be achieved at the same time by the following convex combination of the two rates

$$\gamma_{2,m,n}(\xi_{m,n}) = (1 - \xi_{m,n})\, \delta^+_{m,n} + \xi_{m,n}\gamma_{1,m,n},$$

where $\delta^+_{m,n} \leqslant \gamma_{1,m,n}$ with $\delta_{m,n}/\delta^+_{m,n} \to 0$ (e.g. by choosing $\delta^+_{m,n} = \min(\delta_{m,n} \cdot \log(mn/(m+n)), \gamma_{1,m,n})$) and $\xi_{m,n} \geqslant 0$ is a sequence that is bounded away from 1. In particular, it is recommended to use a sequence $\xi_{m,n} \to 0$ sufficiently fast to have the optimal contraction rates. The enlargement of $\delta_{m,n}$ by $\delta^+_{m,n}$ is necessary to guarantee that $(II)$ holds with asymptotic probability 1 under the null hypothesis. In the simulation study, we use the choices $\xi_{m,n} = e^{-100\frac{m+n}{mn}}$ and $\delta^+_{m,n} = \left(\frac{m+n}{mn}\right)^{\frac{3}{4}}$, i.e. the contraction rate for Lipschitz-continuous depths from Corollary 2.10 without enlargement which is justified by Theorem 2.6 (a). We refer to the corresponding test as Joint-CC test (for **C**onvex-**C**ombination). It is illustrated in Figure 2 (e).

Figure 2 suggests indeed that these statistics have much improved small sample size behaviour compared to the two $\mathcal{LS}$-statistics or the $\mathcal{LS}$-maximum-statistic as in Corollary 2.10 (b). Additional simulation results in this spirit can be found in Section C.1, confirming for functional and multivariate data the observations made in this section.

## 3. Validation of assumptions for selected depth functions

This section is dedicated to study whether there are statistical depth functions that satisfy $\text{Asm}_1$ - $\text{Asm}_4$, required in the results of Section 2. Our main objective is to study functional instances and, in turn, we consider the $h$-depth [22] and the multivariate functional generalisation [17] of the integrated depth [35]. For simplicity, we refer to the generalisation as integrated depth in what follows.

Let the sample space $\mathbb{H}$ be a separable Hilbert space, $P$ a probability measure on $\mathbb{H}$ (shortly $P \in \mathcal{P}_{\mathbb{H}}$), $X \sim P$, and consider a univariate kernel function $K(\cdot)$ as well as a bandwidth $h > 0$. Then, the $h$-depth of an element $x \in \mathbb{H}$ with respect to $P$ is defined as

$$D_h(x, P) := \mathbb{E}\left(\frac{1}{h} \cdot K\left(\frac{\|x - X\|_{\mathbb{H}}}{h}\right)\right). \tag{13}$$

Its sample version with respect to $\widehat{P}_m$ (based on iid $X_1, ..., X_m \sim P$) reads

$$D_h(x, \widehat{P}_m) = \frac{1}{mh}\sum_{i=1}^{m} K\left(\frac{\|x - X_i\|_{\mathbb{H}}}{h}\right).$$

The opportunity of choosing the norm (respectively the underlying sample space), the kernel and the bandwidth offers a high amount of flexibility for several classes of data.

Let $C^d([0,1])$ be the space of $d$-variate continuous functions over $[0,1]$. For a multivariate depth function $D_d$, $d \in \mathbb{N}$, the corresponding integrated depth of an element $x \in C^d([0,1])$, with respect to a distribution $P \in \mathcal{P}_{C^d([0,1])}$ reads

$$ID(x, P) := \int_0^1 D_d(x(t), P(t))\, \mathrm{d}t. \tag{14}$$

For the sample analog, unless stated otherwise, we make use of the integrated depth with respect to $\widehat{P}_m$. The flexibility here comes from the multivariate depth selected, as different depth functions highlight different characteristics of the underlying distribution.

$Asm_1$ is related to the continuity assumption on the underlying distribution for univariate rank-based methods. This assumption, and even the stronger assumption of Lipschitz-continuity, has also been made in related previous works. In particular, Zuo and He [93] as well as Ramsay and Chenouri [66, 68, 67] do not check whether the different depth functions satisfy it. We proceed here similarly. Note that this is not an eccentric assumption as, for instance, [26], Lemma 6.1 shows that the Tukey depth is continuous in the first argument when computed with respect to absolutely continuous distributions. Proposition 3.3, below, proves that the univariate Tukey depth computed with respect to absolutely continuous distributions satisfies this assumption with $\beta = 1$.

For the $h$-depth the below result holds.

**Theorem 3.1.** *Let $\mathbb{H}$ be a separable Hilbert space and $P \in \mathcal{P}_{\mathbb{H}}$. $D_h$ fulfills $Asm_2$ and $Asm_3$ for any bounded kernel $K(\cdot)$ and any fixed bandwidth $h > 0$.*

The following result reveals that the integrated depth fulfills assumptions $Asm_2$ and $Asm_3$ when based on a multivariate depth taking values in $[0, 1]$ that fulfills them (see Section 3.1). While Serfling and Zuo [74] restrict the notion of multivariate depths to have an image in $[0, 1]$, the latter has not been fulfilled by all multivariate depths proposed more recently, such as e.g. the $h$-depth with a kernel not bounded from above by 1. Nevertheless, all depths considered here meet this requirement. The result is established in $C^d[0, 1]$, however, it can be replaced by another appropriate space.

**Theorem 3.2.** *Let $P \in \mathcal{P}_{C^d[0,1]}$ with $P(t)$ (the probability distribution of a functional observation evaluated at $t$) continuous for all $t \in [0, 1]$. $ID$ satisfies $Asm_2$, and respectively $Asm_3$, when based on a multivariate depth function taking values in $[0, 1]$ that fulfills, uniformly in $P$, $Asm_2$, and respectively $Asm_3$.*

As explained in Section 2.3, $Asm_4$ is different in Zuo and He [93]. We assume $Asm_4$, in the same manner to $Asm_1$. This conjecture seems plausible for the integrated depth family and the $h$-depth because they can generally take a continuous range of values even when computed with respect to the empirical distribution. It is also worth pointing out that $Asm_4$ is satisfied if the following conditional Hölder-continuity property holds:

$$\mathfrak{P}_m \leqslant \left(\frac{L}{m}\right)^{\beta} \text{ with } \mathbb{E}(L^{2\beta}) < C < \infty$$

for the Hölder-exponent $\beta$ from $Asm_1$ and $C$ not depending on $m$, where

$$\begin{aligned}\mathfrak{P}_m &= \mathbb{P}\left( |D(X_1, \widehat{P}_m^{-\{1\}}(\cdot)) - D(X, \widehat{P}_m^{-\{1\}}(\cdot))| \leqslant \frac{2C_{det}}{m} \middle| X_1, ..., X_m \right) \\ &= \mathbb{P}\Bigg( D(X, \widehat{P}_m^{-\{1\}}(\cdot)) \in \\ &\qquad \left[ D(X_1, \widehat{P}_m^{-\{1\}}(\cdot)) - \frac{2C_{det}}{m}, D(X_1, \widehat{P}_m^{-\{1\}}(\cdot)) + \frac{2C_{det}}{m} \right] \Bigg| X_1, ..., X_m \Bigg),\end{aligned}$$

which relates to $Asm_4$.

### 3.1. Multivariate depths

Theorem 3.2 requires a multivariate depth function satisfying Asm$_2$ and Asm$_3$. Multivariate depths were the ones analysed in studying the $\mathcal{LS}$-statistic in Zuo and He [93]; however, the assumptions here differ, as they have been adapted in Section 2.

As the $h$-depth, cf. (13), can also be computed for multivariate data, we have already a multivariate depth satisfying the required assumptions. In Subsection 3.1.1 we provide the definition of other multivariate depths and in Subsection 3.1.2 we study whether they satisfy Asm$_2$ and Asm$_3$.

#### 3.1.1. Multivariate depth definitions

The Tukey depth [85] is the most well known multivariate depth function, because of the properties it satisfies. Its only issue is its computational time [28], which is easily solved by approximating it by the random Tukey depth [20]. The Tukey depth of $x \in \mathbb{R}^d$ with respect to a distribution $P$ on $\mathbb{R}^d$ is

$$D_T(x, P) = \inf_{H \in \mathcal{H}(x)} P(H), \tag{15}$$

where $\mathcal{H}(x)$ denotes the family of closed half spaces in $\mathbb{R}^d$ that contain $x$. For the sample analogue, we make use of the Tukey depth evaluated at $\widehat{P}_m$.

The simplicial depth [48] of $x \in \mathbb{R}^d$ with respect to a probability measure $P$ on $\mathbb{R}^d$ is

$$D_S(x, P) = \mathbb{P}(x \in \operatorname{conv}\{X_1, ..., X_{d+1}\}),$$

where $\operatorname{conv}\{\cdot\}$ denotes the closed convex hull and $X_1, ..., X_{d+1}$ iid $d$-variate random variables with distribution $P$. The empirical simplicial depth of $x \in \mathbb{R}^d$ is

$$D_S(x, \widehat{P}_m) = \binom{m}{d+1}^{-1} \sum_{1 \leqslant i_1 < \ldots < i_{d+1} \leqslant m} \mathbb{1}_{\{x \in \operatorname{conv}\{X_{i_1}, \ldots, X_{i_{d+1}}\}\}}. \tag{16}$$

Hence, it can be regarded as a U-statistic. In the proof of Theorem 3.4 we also consider an empirical depth version that is a V-statistic in an intermediate step. For the simplicial depth the latter was proposed in [27], page 2 and reads

$$D_{S,\mathrm{mod}}(x, \widehat{P}_m) = \int \mathbb{1}_{\{x \in \operatorname{conv}\{x_1, \ldots, x_{d+1}\}\}} \, \mathrm{d}\widehat{P}_m^{d+1}(x_1, ..., x_{d+1}).$$

A less technical theory for the simplicial depth than the one used here is developed in [29]; it is based on the interior of the convex hull.

The spherical depth, proposed by [30], is very similar to the simplicial depth. In contrast to closed simplices, closed balls are taken into account. In $\mathbb{R}^d$ only two points are required to define such a closed ball. The spherical depth of $x \in \mathbb{R}^d$ with respect to a probability measure $P$ on $\mathbb{R}^d$ is

$$D_{Sph}(x, P) = \mathbb{P}\left((X_1 - x)^T (X_2 - x) \leqslant 0\right),$$

with $X_1, X_2$ iid $d$-variate random variables with distribution $P$.

The lens depth by Liu and Modarres [50] takes hyperlenses into account instead and is defined as

$$D_{Le}(x,P) = \mathbb{P}\left(\|X_1 - X_2\| \geqslant \max(\|x - X_1\|, \|x - X_2\|)\right),$$

with $X_1, X_2$ iid $d$-variate random variables with distribution $P$.

The multivariate band depth (based on $k \in \mathbb{N} \setminus \{1\}$ bands) of a point $x \in \mathbb{R}^d$ with respect to the probability measure $P$ (cf. [54], Section 3) is given by

$$D_{Ba,k}(x,P) = \mathbb{P}\left(\bigcap_{j=1}^{d} \left\{ \min_{\ell \in \{1,\dots,k\}} X_\ell^{(j)} \leqslant x^{(j)} \leqslant \max_{\ell \in \{1,\dots,k\}} X_\ell^{(j)} \right\}\right), \tag{17}$$

with $X_1, ..., X_\ell$ iid $d$-variate random variables with distribution $P$.

### *3.1.2. Multivariate results*

The case of $d = 1$ is a special one as there exists an intrinsic order in $\mathbb{R}$.

**Proposition 3.3.** *Tukey depth with respect to a continuous distribution on $\mathbb{R}$ satisfies Asm$_1$–Asm$_4$, with $\beta = 1$ in Asm$_1$.*

The following theorem is for $\mathbb{R}^d$ and needs for $d \geqslant 2$ the following measurability assumption, (M), to prove that Asm$_2$ holds. It is required in particular for the application of the Dvoretzky-Kiefer-Wolfowitz inequality [2], Theorem 3.1. Such inequality is also used in the proof of Asm$_2$ in Proposition 3.3. $\mathcal{H}$ denotes the family of closed halfspaces in $\mathbb{R}^d$ and $\widehat{P}_{2m}$ the empirical probability measure corresponding to $X_1, ..., X_{2m}$ iid random variables with distribution $P$ on $\mathbb{R}^d$.

(M) For every $m \in \mathbb{N}$, $\sqrt{m} \cdot \sup_{H \in \mathcal{H}} |\widehat{P}_m(H) - P(H)|$ is measurable and

$$\sqrt{m} \cdot \sup_{H \in \mathcal{H}} |\widehat{P}_m(H) - P(H)|,$$

$$\sqrt{m} \cdot \sup_{H \in \mathcal{H}} \left| \widehat{P}_m(H) - \frac{1}{m} \sum_{i=m+1}^{2m} \mathbb{1}_{\{X_i \in H\}} \right| \text{ and}$$

$$\sup_{H \in \mathcal{H}} |\widehat{P}_{2m}(H) - P(H)|$$

are measurable with respect to the completed product measure of $X_1, ..., X_{2m}$.

In the theorem we also make use of the function $\mathfrak{L} : \bigtimes_{j=1}^{k} \mathbb{R}^d \to 2^{\mathbb{R}^d}$, where $2^{\mathbb{R}^d}$ denotes the power set, and of the notation

$$\mathrm{A}(x, y_1, ..., y_{k-1}) = \{y \in \mathbb{R}^d : x \in \mathfrak{L}(y_1, ..., y_{k-1}, y)\}$$

for the class of sets generated by $\mathfrak{L}(\cdot)$. It is also worth saying that the theorem makes use of U-statistics, as we consider depth functions whose empirical versions can be regarded as U-statistics of order $k > 1$, i.e.

$$D(x, \widehat{P}_m) = \frac{1}{\binom{m}{k}} \sum_{1 \leqslant i_1 < \dots < i_k \leqslant m} \mathbb{1}_{\{x \in \mathfrak{L}(X_{i_1}, \dots, X_{i_k})\}}. \tag{18}$$

**Theorem 3.4.** *Let $D : \mathbb{R}^d \times \mathcal{P}_{\mathbb{R}^d} \to \mathbb{R}$ be a function with a sample version, $D(x, \widehat{P}_m)$, that is a U-statistic of fixed order $k > 1$ with a bounded kernel $K_x(\cdot)$. Then, $D(x, \widehat{P}_m)$ satisfies $Asm_3$. If additionally the kernel has the form $K_x(\cdot) = \mathbb{1}_{\{x \in \mathfrak{L}(\cdot)\}}$, the class of sets $\mathrm{A}(x, y_1, ..., y_{n-1})$ has finite VC-dimension and (M) is fulfilled, we have that the pair $D(x, P), D(x, \widehat{P}_m)$ satisfies $Asm_2$.*

Proposition 3.5 below makes use of Theorem 3.4 to prove the assumptions for different depth functions.

**Proposition 3.5.** *Tukey depth, simplicial depth, spherical depth, lens depth and multivariate band depth satisfy $Asm_3$. Furthermore, they satisfy $Asm_2$ when computed with respect to a probability measure that fulfills (M).*

**Remark 3.6.** *There also exists a version of the multivariate band depth (cf. [54], Section 3) defined as*

$$D_{Ba}(x, P) = \sum_{k=2}^{\mathcal{K}} D_{Ba,k}(x, P) \text{ for } \mathcal{K} \in \mathbb{N} \setminus \{1\}.$$

*Its sample version is a sum of the U-statistic-estimators and, consequently, satisfies $Asm_3$ and $Asm_2$ when computed with respect to a probability measure that fulfills (M).*

## 4. Simulations and real data analysis

In the following, we give short summary of the simulation results that are presented in more detail in the appendix. In Subsection 4.2, we apply the Joint-TP test to a dataset of temperature curves from ocean drifters in order to demonstrate its ability to discriminate between El Niño and La Niña years.

### *4.1. Summary of the simulation study*

The simulation study is divided into four parts. In the first we concentrate on the performance of the different tests based on the $\mathcal{LS}$-tuple considered in Section 2.4, regarding the impact of its blind spots. Based on the obtained results, the Joint-TP test is chosen for further analysis in the second and third part. In the second part we compare its performance when basing it on different depth functions. Based on these results we select three that performed well for location differences, for scale difference respectively overall to use them in a comparison with other tests from the literature that are available as R-packages. The fourth part is dedicated to evaluate the performance of the Joint-TP test in the presence of outliers.

*Variations of $\mathcal{LS}$-statistics and their blind spots*

We construct two classes of alternatives to compare the performance of the tests based on the $\mathcal{LS}$-tuple in Section 2.4. One of the classes consists of a simple difference in the location of the data and the other of a simultaneous difference in location and scaling. This is done for multivariate data, making use of the Tukey depth, and for functional data, with the integrated Tukey depth.

The main observation is that the $\mathcal{LS}$-statistics have blind spots like those that occur with univariate data (Figure 2). This results in the $\mathcal{LS}$-difference statistic not being suitable to detect differences in location. Additionally, either $\mathcal{LS}(\widehat{P}_m, \widehat{Q}_n)$ or $\mathcal{LS}(\widehat{Q}_n, \widehat{P}_m)$ has very little power against the class with differences in location and scaling, due to their asymmetry. Moreover, we observe that the empirical size of the $\mathcal{LS}$-maximum statistic is above 20%.

Consequently, the only two suitable candidates are the Joint-CC statistic and the Joint-TP statistic. By construction, the last mentioned one is more conservative, and we restrict our attention to it for the rest of the study. We refer to Section C.1 for further details.

*Joint-TP test for different functional depths*

We compare the small sample behaviour of the Joint-TP test when basing it on different depth functions for functional data. We make use of the following depth functions given in Section 3: integrated depths based on Tukey, simplicial as well as a modification of the simplicial depth and two versions of the $h$-depth – one with $h = 1$ fixed and the other with a data-adaptively chosen $h$. Additionally we use the spatial depth [11], the lens metric depth [39], the functional random Tukey depth [20] and the random projection depth [22]. We consider three different types of functional data – non-smooth, smooth, smooth with high fluctuation – for generating the sets of alternatives for the following scenarios: simple location difference, sine curve location difference, simple scale difference, simultaneous difference in location and scaling. Under the null hypothesis, even for samples consisting of only $m = n = 50$ observations, the empirical size with most depths is below or close to the desired 5% level. Only with the integrated simplicial depth in the setting of smooth fluctuating functions, the empirical size is more than 80%, i.e. too large. That effect vanishes asymptotically, with $m = n = 100$ observations the size is 6.5%. This is because of the disturbed inside-out ordering entailed by the univariate simplicial depth, see Example C.1, which carries over to the integrated simplicial depth. A slight modification, described in Section C.2, eliminates this behaviour, but it leads to a small loss of power.

Detecting complex deviations in the second order structure (i.e. differences between Model 1-3) works well, where the integrated Tukey depth has a somewhat lower power that is still above 55% in the worst case. For detecting differences in the location of the data, the integrated depths give the highest power. With smooth functional data that rarely intersect, the adaptive $h$-depth and the spatial depth offer good performance in the setting of a difference in the location. When it comes to scale differences, the depths that take the distance between the curves into account, instead of pointwise comparisons on the time grid, have a higher power.

In the setting of a simultaneous location-scale difference, using the integrated simplicial depth leads to a high sensitivity for location differences if the data do not fluctuate much. If the adaptive $h$-depth is used instead, there is an increase in the sensitivity regarding differences in the second order structure. When it comes to alternatives with equal location and second order structure, but different higher order structures (i.e. different shape), the tests based on the integrated depths have less power compared to the tests that use depths based on a metric. The reason for this phenomenon is that integrated depths are computed by averaging pointwise computed values over the time axis, but the shape difference does not occur pointwise.

Furthermore, functional random Tukey depth as well as the random projection depth have the advantage of low computational cost, but lead – when applied with the Joint-TP test – to less power than the other functional depths. The lens metric depth as well as the $h$-depth with $h = 1$ show average performance among the other depths that we used with the Joint-TP test.

We recommend applying the Joint-TP test either with the adaptive $h$-depth if the focus is more on differences in the second order structure (respectively higher order structures), or with the integrated simplicial depth if the goal is to give more weight on detecting differences in the location. The integrated Tukey depth is also a good candidate as an allround solution because it leads to a conservative behaviour in comparison to the integrated simplicial depth. The detailed results as well as additional information are provided in Section C.2.

*Comparison of the Joint-TP test with the state-of-the-art functional two-sample tests*

This part is dedicated to compare the performance of the Joint-TP tests based on the integrated Tukey depth, the integrated simplicial depth and the adaptive $h$-depth with several state-of-the-art two-sample tests for functional data. We are using exactly the same sets of functional alternatives as in the previous part, where we compared the Joint-TP test for different functional depths. The considered tests are: ($\mathfrak{a}$) two ANOVA tests – a bias-reduced F-type test [78, 89] and a globalised F-test [90], ($\mathfrak{b}$) a test based on PCA [36], ($\mathfrak{c}$) a test based on the energy distance [70, 83] and ($\mathfrak{d}$) a ball divergence test [65]. Gretton et al. [41] and Zhang and Smaga [91] also proposed tests based on the maximum mean discrepancy similar to the one considered in ($\mathfrak{c}$). We have selected ($\mathfrak{c}$) using 10000 bootstrap replicates because of its availability on CRAN. We have relegated Gretton et al. [41] to the appendix because, while it is also on CRAN, it performs significantly worse than the test in ($\mathfrak{c}$) except for the scenarios in the third line of Figure 12, where both perform badly.

As expected, the ANOVA-based methods ($\mathfrak{a}$) have high power in presence of a difference in the location, but fail to detect any differences in the second order structure. In the case of a simultaneous location-scale difference, they are a good choice as long as the location difference dominates. The PCA-based method ($\mathfrak{b}$) is designed to detect differences in the second order structure and, hence, has no power against pure location-difference alternatives. Nevertheless, the Joint-TP test with the data-adaptive $h$-depth shows an equal, respectively even better, performance than all the other tests, of which some are particularly designed to detect this type of difference, which is surprising. The energy-distance test ($\mathfrak{c}$) has high power against differences in location, but difficulties to discriminate between alternatives in which only the scaling is different. The ball divergence statistic ($\mathfrak{d}$) has less power than the three versions of the Joint-TP test when it comes to pure second order differences in non-smooth or smooth data with less fluctuation. In presence of location differences, it seems to perform better than the Joint-TP tests with smooth data, but in contrast, with non-smooth data the Joint-TP test with the integrated simplicial depth is more sensitive.

All in all, there exist alternatives in which some of these test classes dominate all the other tests, but in general, the choice of the test statistic is closely related to the type of difference in the underlying distributions. The Joint-TP is a good choice when there is few or no information on the type of difference to detect. The detailed simulation results are given in Section C.5.

*Performance in the presence of outliers*

Here, we compare the $\mathcal{LS}$-tests with the above ANOVA-based tests in the presence of magnitude outliers. In particular, we analyse the case of having an additive outlier, where a constant function with a value of +50 was added to the first curve in either one or both samples. The $\mathcal{LS}$-test is robust with respect to this type of outlier both under the null hypothesis, with only a slightly more conservative size, and under the alternative, with comparable power behaviour. This is not true for the ANOVA-based tests where the size and power break down completely to empirical rejection rates of 0% both under the null and alternative hypotheses, with the exception where the additive outlier is added to the sample with higher mean. The detailed simulation results also including some results with more than one outlier per sample are given in Section C.6.

### 4.2. Analysis of ocean drifter temperature data

Temperature curves of ocean drifters contain information on the occurrence of El Niño respectively its counter event La Niña. For instance, Sun and Genton [82] use sea surface data from the 1990s related to El Niño as a data example for their depth based functional boxplot.

For each year between 2002 and 2021 the temperature data for each drifter deploying in the area with coordinates less than 10°N and in between 69°W and 180°W (this area contains the region of the Pacific Ocean where El Niño or La Niña usually appear) were downloaded from [55] accessed via the ERDAP-platform [80]. These data are given on the same equidistant time grid (time lag: 6 hours between two consecutive observations). Curves that lack at least a whole calendar day of observations were removed. For each temperature curve of each drifter we smoothed the data by computing the median temperature

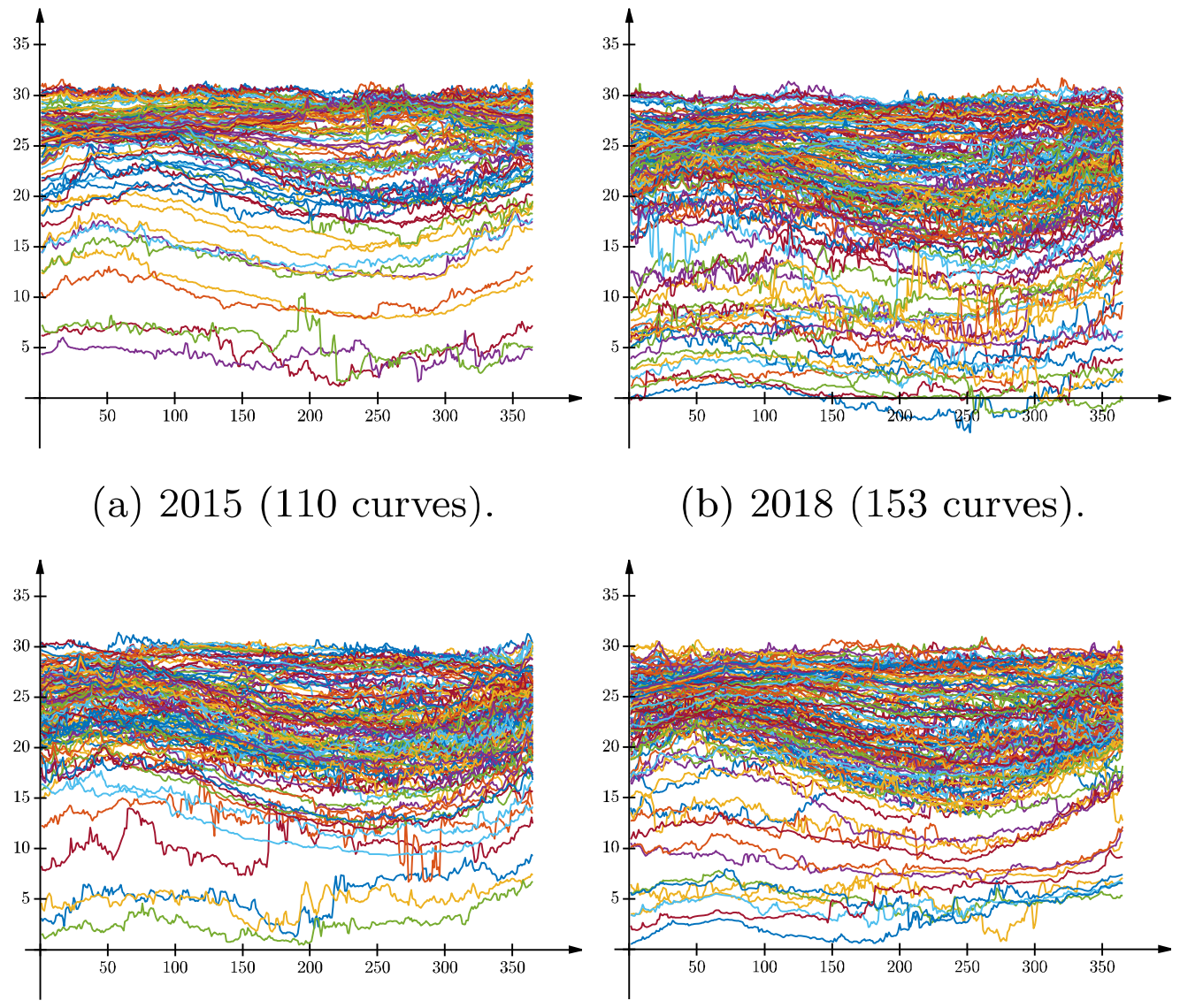

(a) 2015 (110 curves). (b) 2018 (153 curves).

(c) 2020 (130 curves). (d) 2021 (142 curves).

Fig 3. *Temperature curves of all considered drifters for the years 2015, 2018, 2020 and 2021. The sample sizes are balanced. The horizontal axis denotes the day of the year, the vertical axis denotes the temperature in °C. 2015 is a strong El Niño year, 2018 a moderate one, 2020 and 2021 are moderate La Niña years.*

for each day, with NA-values being ignored. In the leap years, we removed the temperature value for February 29 and replaced the values for February 28, respectively March 1, by the mean temperature of February 28 and February 29 respectively by the mean of February 29 and March 1.

Our objective is to investigate whether the Joint-TP test can help to distinguish between El Niño- and La Niña-periods, respectively periods of its different intensities. The periods are selected with respect to the classification given by [64].

First, we compare samples of balanced sizes, namely the strong El Niño year 2015, the weak El Niño year 2018 as well as the moderate La Niña phase 2020 and 2021 (see Figure 3 for the curves). For all pairs of these 4 years, we compute the values of the Joint-TP statistic with the integrated Tukey depth (IT), the integrated simplicial depth (IS) and the $h$-depth with adaptive bandwidth $h$. The corresponding $p$-values and sample sizes are given in Table 1.

With balanced sample sizes the $p$-values given in Table 1 help to clearly distinguish between the strong El Niño year 2015 and the weak El Niño year 2018 as well as the moderate La Niña phase 2020–2021. In Figure 3 the temperature curves of these four periods are plotted. Obviously, there are fewer drifters with low temperature in 2015 than in the other depicted years. Moreover, there are less drifters with temperatures above 30°C in 2021 than in 2018. Between the pair 2020 and 2021, the analysis does not indicate a difference.

TABLE 1

*p-values of the Joint-TP statistic for the data from Figure 3 with the integrated Tukey depth (IT), the integrated simplicial depth (IS) and the adaptive h-depth (h). 2015 is considered as very strong El Niño year, 2018 is considered a weak El Niño period and 2020–2021 is considered a moderate La Niña period.*

| **Sample 1** | **Sample 2** 2018 (153 curves) | 2020 (130 curves) | 2021 (142 curves) |
|---|---|---|---|
| 2015 (110 curves) | $0.341 \cdot 10^{-3}$ (IT)<br>$0.122 \cdot 10^{-3}$ (IS)<br>$0.703 \cdot 10^{-3}$ ($h$) | $0.791 \cdot 10^{-3}$ (IT)<br>$0.488 \cdot 10^{-3}$ (IS)<br>$0.007 \cdot 10^{-3}$ ($h$) | $0.721 \cdot 10^{-3}$ (IT)<br>$0.392 \cdot 10^{-3}$ (IS)<br>$0.009 \cdot 10^{-3}$ ($h$) |
| 2018 (153 curves) | | 0.068 (IT)<br>0.067 (IS)<br>0.017 ($h$) | 0.365 (IT)<br>0.351 (IS)<br>0.051 ($h$) |
| 2020 (130 curves) | | | 0.266 (IT)<br>0.292 (IS)<br>0.362 ($h$) |

Next, we compare samples of unbalanced sizes, namely the temperature curves from 2002 (a moderate El Niño year) and the pooled temperature curves of the years 2011–2012 (a moderate La Niña period) as well as of the years 2020–2021 (another moderate La Niña period). The unbalancedness arises because of an increase of interest in the Global drifter program [55], which results in more drifters being added every year. The plots of the curves for these three periods are given in Figure 3 (c)–(d) and Figure 4. The corresponding $p$-values for the comparison between the pooled 2020–2021 sample (with 272 curves) versus the years 2002 (with 52) as well as (pooled) 2011-12 (with 54 curves) are given in Table 2 again for the Joint-TP statistic with the integrated Tukey depth (IT), the integrated simplicial depth (IS) and the $h$-depth with adaptive bandwidth $h$.

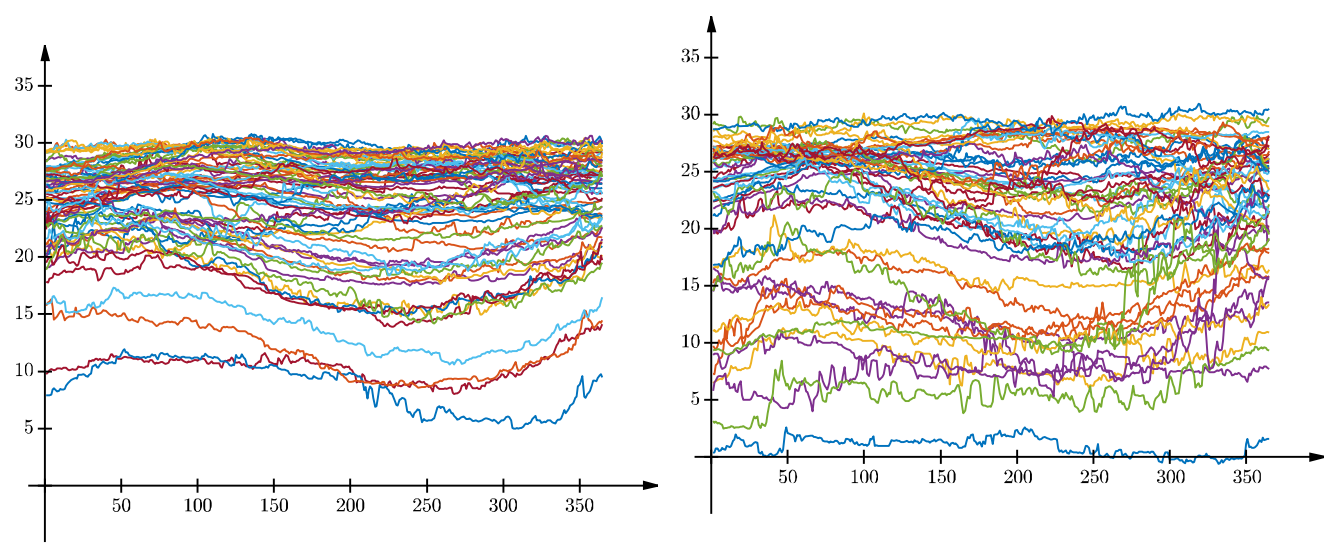

(a) 2002 (52 curves).

(b) 2011-2012 (pooled, 54 curves).

FIG 4. *Temperature data of the year 2002 and the period 2011–2012 with unbalanced sample sizes. The horizontal axis denotes the day of the year, the vertical axis denotes the temperature in °C. The period in (a) is considered as moderate El Niño year, the years in (b) are a moderate La Niña period.*

Table 2

*p-values of the Joint-TP statistic for the data from Figure 4 for unbalanced sample sizes with the integrated Tukey depth (IT), the integrated simplicial depth (IS) and the adaptive h-depth (h): 2011–2012 as well as 2020–2021 are considered as moderate La Niña periods. In 2002 a moderate El Niño occured.*

| **Sample 1** | **Sample 2** | |
|---|---|---|
| | 2002 (52 curves) | Pooled 2011–2012 (54 curves) |
| Pooled 2020–2021 (272 curves) | 0.039 (IT)<br>0.021 (IS)<br>0.001 ($h$) | 0.169 (IT)<br>0.347 (IS)<br>0.203 ($h$) |

With such unbalanced sample sizes the $p$-values for the Joint-TP statistic still indicate a difference between the underlying distributions of the temperature curves from 2002 (a moderate El Niño year) and the pooled temperature curves of the years 2020–2021 (a moderate La Niña period) but are much larger than what has been observed in the balanced situation (with both samples having between 100 and 150 curves). Indeed, in contrast to both of the pooled periods (Figure 3 (c)–(d), Figure 4 (b)), there are less drifters with lower temperature in 2002 (Figure 4 (a)). Similarly to before, the analysis indicates no difference between the distributions of the pooled temperature curves of the years 2011–2012 and the pooled temperature curves of the years 2020–2021 (both moderate La Niña periods).

Furthermore, in Section C.7, we make use of all the curves that were in the relevant area on January 1st of the respective year, including those that are only observed over part of the grid domain. This increases significantly the number of sample curves, depending on the year by at least a factor of two, often much more. There, we only use integrated data depth functions as they can be easily adapted to the situation of partially observed data by integrating only over the observed part of the domain, following Elías et al. [31]. Probably due to the larger sample size, all significant results remain significant with smaller (partly by magnitudes) $p$-values. Interestingly, that analysis suggests that the weak El Niño year of 2018 may be closer to the moderate La Niña period of 2020 than that of 2021. This is consistent with the shapes of the curves of the Oceanic Niño index given in [64].

## 5. Discussion and outlook

In this paper, we revisit a depth-based test for the two-sample problem originally proposed by Liu and Singh [52], later picked up by several other authors for further investigation. While the original test is able to detect various location, scale and location-scale differences, it was originally proposed as a one-sided test for a scale increase, possibly combined with a location shift. In this paper, we shed some additional light onto this situation by considering the $\mathcal{LS}$-tuple jointly: This indicates in particular, that scale increases, possibly combined with location shifts do typically not fall into the *blind spots* as made visible in Figure 2. For such a one-sided test, the inherent asymmetry of the original test is not problematic. Nevertheless, Liu and Singh [52] already propose a symmetric version of their test, where the depth values are calculated with respect to the pooled sample, noting that such a test will only have power against scale differences but not against differences in location. Motivated by the asymptotic properties of the joint $\mathcal{LS}$-tuple of the original statistics with exchanged labels, we propose a different class of symmetric $\mathcal{LS}$-statistics improving upon both the size and the power behaviour. The corresponding observations also give insight into how one could construct unbiased one-sided tests with better size behaviour, see Section C.4 for more information. Detailed analysis is left for future work.

We extend the asymptotic theory for the $\mathcal{LS}$-test as given by Zuo and He [93] to our proposed generalisation, the Joint-TP test. Indeed, the main term and the remainder terms in their proof do have a mathematical correspondence to the location of the $\mathcal{LS}$-tuple which directly corresponds to our construction of the rejection regions of the Joint-TP test. Based on this asymptotic theory, the critical values can be directly computed in a time efficient way. Incidentally, we fill a gap in the original proof and give an illustrative and mathematically interesting counterexample to their main result under the alternative.

Furthermore, we extend the test from multivariate to functional data by combining it with several functional depths, which did not exist when the test was first proposed by Liu and Singh [52]. In particular, we apply the test with the integrated Tukey and simplicial depth or the $h$-depth. In the context of functional data, it is worthwhile mentioning that our Joint-TP test makes use of the full information within the data and does not require dimension reduction such that a loss of information by considering a finite-dimensional subspace is excluded.

A large simulation study and an application to ocean drifter data illustrate that the corresponding functional two-sample tests are competitive with the potential of being more robust in comparison to state-of-the-art methodology for functional data.

In practise, functional data are often not observed over dense and equally spaced grid domains. Different scenarios are possible: For instance, a curve can be observed over a grid domain with sparse and dense regions or it could be observed only on a few grid points, which may or may not be concentrated on certain regions of the domain. As such, it is an interesting research question in its own right which functional depth functions can be well adapted to deal with this situation, see [31, 75] for some first results in that direction. When dealing with a dataset that contains curves not observed over the whole grid domain, our proposed test statistic is in principle applicable. It simply requires that an empirical version of the depth can be calculated and evaluated. For integrated depths this can be achieved by integrating over the observed domains, see [31]. Because our assumptions to obtain the limit distributions depend on the convergence rates of the empirical to the population version (see e.g. Assumption $\mathrm{Asm}_2$), obtaining such rates for empirical functional depths in the light of a data containing curves not observed over the whole grid domain is another interesting question for future research. The results on the ocean drifter dataset including those curves not observed over the whole grid domain indicate that our proposed tests also work well in practice under this type of scenario.

Furthermore, future work may include constructing tests for change point problems based on the notion of depth-based ranks as well as investigating whether this test is applicable to other types of data for which statistical depth functions exist, e.g. fuzzy data [37, 38], random graphs and networks [34] or text data [8].

## Appendix A: Proofs of Section 2

We will move the proof of Theorem 2.1 to Section A.2 behind the proofs of Section 2.3 as the notation and part of the arguments are being used.

### A.1. Proof of Theorem 2.3 and Corollary 2.5

In this section, we are under the null hypothesis with $Q = P$.

The proof is an extension of the proof of Theorem 1 in [93], filling a gap there (see also Remark 2.4). As such we adopt part of their notation and only sketch the parts of the proof that are already done there.

Throughout the proof we make use of the following notation

$$I(x, y, P) = \mathbb{1}_{\{D(x,P)<D(y,P)\}} + \frac{1}{2}\mathbb{1}_{\{D(x,P)=D(y,P)\}}. \tag{19}$$

If instead the ranks (1) (as in [93]) are used, we denote $\widetilde{I}(x, y, P) = \mathbb{1}_{\{D(x,P)\leqslant D(y,P)\}}$ and similarly, for all decompositions, the tilde-term is used in combination with (1). As $\mathcal{LS}(P, P) = \frac{1}{2}$, it holds

$$\begin{aligned}
&\mathcal{LS}(\widehat{P}_m, \widehat{Q}_n) - \frac{1}{2} \\
&= \int\int I(x, y, \widehat{P}_m)\, \mathrm{d}\widehat{P}_m(x)\, \mathrm{d}\widehat{Q}_n(y) - \int\int I(x, y, P)\, \mathrm{d}P(x)\, \mathrm{d}Q(y) \\
&= A_1(\widehat{P}_m, \widehat{Q}_n) + A_2(\widehat{P}_m, \widehat{Q}_n),
\end{aligned} \tag{20}$$

where

$$\begin{aligned}
A_1(\widehat{P}_m, \widehat{Q}_n) &= \int\int I(x, y, P)\, \mathrm{d}\widehat{P}_m(x)\, \mathrm{d}\widehat{Q}_n(y) - \int\int I(x, y, P)\, \mathrm{d}P(x)\, \mathrm{d}Q(y), \\
A_2(\widehat{P}_m, \widehat{Q}_n) &= \int\int (I(x, y, \widehat{P}_m) - I(x, y, P))\, \mathrm{d}\widehat{P}_m(x)\, \mathrm{d}\widehat{Q}_n(y).
\end{aligned}$$

Under the null hypothesis it holds by definition $I(X, Y, P) = 1 - I(Y, X, P)$ hence $A_1(\widehat{P}_m, \widehat{Q}_n) = -A_1(\widehat{Q}_n, \widehat{P}_m)$. In fact, $A_1$ is the leading term responsible for the asymptotic behaviour and also the reason why the two statistics are perfectly negatively correlated asymptotically while $A_2$ converges to 0 faster than $A_1$ and is asymptotically negligible (see Figure 2 for an empirical illustration). Note, however, that this is only true under the null hypothesis because under the alternative for the permuted samples the integrand becomes $I(x,y,Q) \neq I(x,y,P)$ (see also Theorem 2.1 and Remark 2.4). The following lemma gives the asymptotic behaviour of the leading term $A_1$.

**Lemma A.1.** *If $D(X_1,P)$ is continuously distributed (as e.g. guaranteed by Asm$_1$) and under (9), it holds $A_1(\widehat{P}_m,\widehat{Q}_n) = -A_1(\widehat{Q}_n,\widehat{P}_m)$ and*

$$\sqrt{\frac{12\cdot m\cdot n}{m+n}}\,A_1(\widehat{P}_m,\widehat{Q}_n) \xrightarrow{\mathcal{D}} N(0,1).$$

*Proof.* We have already observed the first statement above the lemma. Decompose further into

$$\begin{aligned} A_1(\widehat{P}_m,\widehat{Q}_n) &= \int\int I(x,y,P)\,\mathrm{d}P(x)\,\mathrm{d}(\widehat{Q}_n(y)-Q(y)) \\ &\quad + \int\int I(x,y,P)\,\mathrm{d}(\widehat{P}_m(x)-P(x))\,\mathrm{d}Q(y) \\ &\quad + \int\int I(x,y,P)\,\mathrm{d}(\widehat{P}_m(x)-P(x))\,\mathrm{d}(\widehat{Q}_n(y)-Q(y)) \\ &= A_{1,1}(\widehat{P}_m,\widehat{Q}_n)+A_{1,2}(\widehat{P}_m,\widehat{Q}_n)+A_{1,3}(\widehat{P}_m,\widehat{Q}_n). \end{aligned} \tag{21}$$

Concerning the third term, it holds

$$A_{1,3}(\widehat{P}_m,\widehat{Q}_n) = O_P\left(\frac{1}{\sqrt{mn}}\right). \tag{22}$$

This fact was already stated in Lemma 1 (i) in [93], but the proof was omitted. Thus, we give a sketch here: by the independence of the two samples and the independence within the samples and $I(x,y,P)\leqslant 1$ it holds

$$\begin{aligned} &\mathbb{E}\left(A_{1,3}(\widehat{P}_m,\widehat{Q}_n)\right)^2 \\ &= \mathbb{E}\left(\mathbb{V}\mathrm{ar}\left(\frac{1}{n}\sum_{j=1}^{n}\int I(x,Y_j,P)\,\mathrm{d}(\widehat{P}_m(x)-P(x))\,\middle|\,X_1,\dots,X_m\right)\right) \\ &\leqslant \frac{1}{n^2}\sum_{j=1}^{n}\mathbb{E}\left(\int I(x,Y_j,P)\,\mathrm{d}(P_n(x)-P(x))\right)^2 \\ &= \frac{1}{n}\,\mathbb{E}\left(\mathbb{V}\mathrm{ar}\left(\frac{1}{m}\sum_{i=1}^{m}I(X_i,Y_1,P)\,\middle|\,Y_1\right)\right) = \frac{1}{mn}. \end{aligned}$$

So, the statement follows from Markov's inequality.

By (2), it holds $\int I(x,y,P)\,\mathrm{d}P(x) = R(y,P)$ and due to the continuity of $D(X_1,P)$ under Asm$_1$ $R(Y_1,P)\sim U(0,1)$ under the null hypothesis as well as $\int I(x,y,P)\,\mathrm{d}Q(y) = 1-R(x,P)$ with $1-R(X_1,P)\sim U(0,1)$.

Therefore, $A_{1,1}$ and $A_{2,2}$ are independent centred sums of length $n$ and $m$ respectively of iid $U(0,1)$ random variables, i.e. have variance $1/12$. So, the results follows by a standard application of the central limit theorem. □

We make use of the following inequality in the proof of Lemma A.3.

**Lemma A.2.** *It holds for all probability measures $P_j$, $j = 1, 2$*

$$|I(x,y,P_1) - I(x,y,P_2)| \leqslant \mathbb{1}_{\{|D(x,P_1)-D(y,P_1)|\leqslant 2\sup_{z\in\Omega}|D(z,P_1)-D(z,P_2)|\}}.$$

*Proof.* Since the left hand side is 0 when $D(x,P_1) < D(y,P_1)$ and $D(x,P_2) < D(y,P_2)$ or when in both expressions the inequality is given by $>$, it is sufficient to show that the right hand side is 1 if $D(x,P_1) \leqslant D(y,P_1)$ and $D(x,P_2) \geqslant D(y,P_2)$ or vice versa. Indeed, in the first case,

$$\begin{aligned}
|D(x,P_1) - D(y,P_1)| &= D(y,P_1) - D(x,P_1) \\
&\leqslant D(y,P_1) - D(x,P_1) + D(x,P_2) - D(y,P_2) \\
&\leqslant |D(y,P_1) - D(y,P_2)| + |D(x,P_2) - D(x,P_1)| \\
&\leqslant 2\sup_{z\in\Omega}|D(z,P_1) - D(z,P_2)|.
\end{aligned}$$

The second case is analogous. □

**Lemma A.3.** *Under Asm₁ - Asm₄ and (9), it holds*

$$\sqrt{\frac{12\cdot m\cdot n}{m+n}}\,A_2(\widehat{P}_m,\widehat{Q}_n) \xrightarrow{P} 0, \qquad \sqrt{\frac{12\cdot m\cdot n}{m+n}}\,A_2(\widehat{Q}_n,\widehat{P}_m) \xrightarrow{P} 0.$$

*More precisely,*

$$A_2(\widehat{Q}_n,\widehat{P}_m) = O_P(\delta_{m,n}), \qquad A_2(\widehat{P}_m,\widehat{Q}_n) = O_P(\delta_{m,n}),$$

*where*

$$\delta_{m,n} = \begin{cases} \left(\frac{m+n}{mn}\right)^{\beta} & \text{for } \frac{1}{2} < \beta \leqslant \frac{2}{3}, \\ \left(\frac{m+n}{mn}\right)^{\frac{2+\beta}{4}} & \text{for } \frac{2}{3} \leqslant \beta \leqslant 1. \end{cases}$$

*Proof.* We only prove the first equation as the second follows analogously by equality of distributions under the null hypothesis. Firstly, we decompose $A_2$ into

$$\begin{aligned}
A_2(\widehat{P}_m,\widehat{Q}_n) &= \int\int (I(x,y,\widehat{P}_m) - I(x,y,P))\,\mathrm{d}\widehat{P}_m(x)\,\mathrm{d}(\widehat{Q}_n(y) - Q(y)) \\
&\quad + \int\int (I(x,y,\widehat{P}_m) - I(x,y,P))\,\mathrm{d}\widehat{P}_m(x)\,\mathrm{d}Q(y) \\
&= A_{2,1}(\widehat{P}_m,\widehat{Q}_n) + A_{2,2}(\widehat{P}_m)
\end{aligned} \tag{23}$$

We first prove that $A_{2,1}(\widehat{P}_m, \widehat{Q}_n) = O_P(1/(n^{1/2}\, m^{\beta/4}))$. While the proof is analogous to the proof of Lemma 1 (iii) in [93] we sketch it for the sake of completeness because of our many (minor) modifications: By Jensen's inequality

$$\begin{aligned}
&\mathbb{E}(A_{2,1}(\widehat{P}_m, \widehat{Q}_n)^2) \\
&\leqslant \mathbb{E}\left(\int \left(\int (I(x,y,\widehat{P}_m) - I(x,y,P))\, \mathrm{d}(\widehat{Q}_n(y) - Q(y))\right)^2 \mathrm{d}\widehat{P}_m(x)\right) \\
&= \mathbb{E}\left[\left(\int (I(X_1,y,\widehat{P}_m) - I(X_1,y,P))\, \mathrm{d}(\widehat{Q}_n(y) - Q(y))\right)^2\right] \\
&= \mathbb{E}\left(\mathbb{V}\mathrm{ar}\left(\frac{1}{n}\sum_{j=1}^{n}(I(X_1,Y_j,\widehat{P}_m) - I(X_1,Y_j,P)) \,\middle|\, X_1,...,X_m\right)\right) \\
&\leqslant \frac{1}{n}\mathbb{E}\left(\mathbb{E}\left(\left(I(X_1,Y_1,\widehat{P}_m) - I(X_1,Y_1,P)\right)^2 \,\middle|\, X_1,...,X_m\right)\right).
\end{aligned}$$

By Lemma A.2, $\mathrm{Asm}_1$ and $\mathrm{Asm}_2$ it follows

$$\begin{aligned}
&\mathbb{E}(A_{2,1}(\widehat{P}_m, \widehat{Q}_n)^2) \\
&\leqslant \frac{1}{n}\,\mathbb{E}\left(\mathbb{P}\left(|D(Y_1,P) - D(X_1,P)| \leqslant 2\sup_{z\in\Omega}|D(z,\widehat{P}_m) - D(z,P)|\right)\right) \\
&= O\left(\frac{1}{n}\right)\,\mathbb{E}\left(\sup_{z\in\Omega}|D(z,\widehat{P}_m) - D(z,P)|^\beta\right) = O\left(\frac{1}{n\, m^{\beta/2}}\right) = O\left(\delta^2_{m,n}\right). \qquad (24)
\end{aligned}$$

So, the assertion follows from applying the Markov inequality. Finally, we need to investigate the term $A_{2,2}(\widehat{P}_m)$, for which the rate was only stated without proof in [93] (see also Remark A.4 below). Unlike above, in $A_{2,2}$, we sum over $X_i$ while at the same time all summands depend on the whole $X$-sample via $\widehat{P}_m$, therefore conditional independence as above cannot be used. Nevertheless, we use an $L^2$-argument, explicitly considering the mixed term (effectively corresponding to the covariances) that arises. We will deal with this mixed term in a similar fashion as in (24), by looking at the difference when $\widehat{P}_m$ and $\widehat{P}_m^{-\{1,2\}}$ (notation as introduced in $\mathrm{Asm}_3$) instead of $\widehat{P}_m$ and $P$ are being used. For the latter, the same type of arguments as above can be used. However, additional regularity conditions are required to treat the difference. With the notation $J(x,\widehat{P}_m) = \int (I(x,y,\widehat{P}_m) - I(x,y,P))\, \mathrm{d}Q(y)$ and $\widehat{P}_m^{-\{1,2\}}$ as defined in $\mathrm{Asm}_3$ it holds

$$\begin{aligned}
\mathbb{E}\left(A_{2,2}(\widehat{P}_m)^2\right) &= \mathbb{E}\left(\frac{1}{m}\sum_{i=1}^{m} J(X_i,\widehat{P}_m)\right)^2 \\
&= \mathbb{E}\left[\left(\frac{1}{m}\sum_{i=1}^{m}\left(J(X_i,\widehat{P}_m) - J(X_i,\widehat{P}_m^{-\{i\}})\right) + \frac{1}{m}\sum_{i=1}^{m} J(X_i,\widehat{P}_m^{-\{i\}})\right)^2\right]
\end{aligned}$$

$$\leqslant 2 \cdot \mathbb{E}\left[\left(\frac{1}{m}\sum_{i=1}^{m}\left(J(X_i, \widehat{P}_m) - J(X_i, \widehat{P}_m^{-\{i\}})\right)\right)^2\right]$$
$$+ 2 \cdot \mathbb{E}\left[\left(\frac{1}{m}\sum_{i=1}^{m} J(X_i, \widehat{P}_m^{-\{i\}})\right)^2\right] = 2 \cdot B_1 + 2 \cdot B_2.$$

For $B_1$, we obtain by the analytic version of Jensen's inequality (applied to the squared sum inside the expectation)

$$B_1 \leqslant \mathbb{E}\left[\left(J(X_1, \widehat{P}_m) - J(X_1, \widehat{P}_m^{-\{1\}})\right)^2\right].$$

By Lemma A.2 and $\mathrm{Asm}_3$ we get

$$\left|J(X_1, \widehat{P}_m) - J(X_1, \widehat{P}_m^{-\{1\}})\right| \leqslant \int |I(X_1, y, \widehat{P}_m) - I(X_1, y, \widehat{P}_m^{-\{1\}})| \,\mathrm{d}Q(y)$$
$$\leqslant P\left(D(X_1, \widehat{P}_m^{-\{1\}}) - \frac{2C_{\mathrm{det}}}{m} \leqslant D(Y, \widehat{P}_m^{-\{1\}}) \leqslant D(X_1, \widehat{P}_m^{-\{1\}}) + \frac{2C_{\mathrm{det}}}{m} \,\middle|\, X_1, X_3, ..., X_m\right). \tag{25}$$

Because of the above inequality and $\mathrm{Asm}_4$, it holds

$$B_1 = O(m^{-2\beta}) = O\left(\delta_{m,n}^2\right).$$

For $B_2$ we have

$$\mathbb{E}\left(B_2^2\right) = \mathbb{E}\left[\left(\frac{1}{m}\sum_{i=1}^{m} J(X_i, \widehat{P}_m^{-\{i\}})\right)^2\right]$$
$$= \frac{1}{m}\mathbb{E}J^2(X_1, \widehat{P}_m^{-\{1\}}) + \left(1 - \frac{1}{m}\right)\mathbb{E}\left(J(X_1, \widehat{P}_m^{-\{1\}}) \cdot J(X_2, \widehat{P}_m^{-\{1\}})\right).$$

For any $x$ and $\tilde{P}$ it holds by Lemma A.2 and $\mathrm{Asm}_1$

$$\left|J(x, \tilde{P})\right| \leqslant \int \mathbb{1}_{\{|D(y,P) - D(x,P)| \leqslant 2\sup_{z\in\Omega}|D(z,P) - D(z,\tilde{P})|\}} \,\mathrm{d}Q(y)$$
$$\leqslant 2C \sup_{z\in\Omega} |D(z, P) - D(z, \tilde{P})|^{\beta}.$$

By this and $\mathrm{Asm}_2$ it holds

$$\frac{1}{m}\mathbb{E}J^2(X_1, \widehat{P}_m^{-\{1\}}) = O\left(\frac{1}{m}\mathbb{E}\sup_{z\in\Omega}|D(z, P) - D(z, \widehat{P}_m^{-\{1\}})|^{2\beta}\right)$$
$$= O\left(\frac{1}{m^{1+\beta}}\right) = O\left(\delta_{m,n}^2\right). \tag{26}$$

Finally, for any probability measure $P_1$ (independent of $X$ and $Y$) it holds under the null hypothesis of $X \overset{D}{=} Y$, i.e. $P = Q$ and $I(X, Y, P_1) \overset{\mathcal{D}}{=} I(Y, X, P_1)$, that

$$\int\int I(x, y, P_1)\ \mathrm{d}P(x)\ \mathrm{d}Q(y) \\ = \frac{1}{2}\int\int (I(x, y, P_1) + I(y, x, P_1))\ \mathrm{d}P(x)\ \mathrm{d}P(y) = \frac{1}{2}, \quad (27)$$

because $I(x, y, P_1) + I(y, x, P_1) = 1$. So,

$$\mathbb{E}\left(J(X_1, \widehat{P}_m^{-\{1,2\}}) \,\middle|\, X_3, ..., X_m\right) = 0 \quad \text{a.s.}, \\ \mathbb{E}\left(J(X_1, \widehat{P}_m^{-\{1\}}) \,\middle|\, X_2, ..., X_m\right) = 0 \quad \text{a.s.} \quad (28)$$

Moreover

$$\begin{aligned} &\mathbb{E}\left(J(X_1, \widehat{P}_m^{-\{1\}}) \cdot J(X_2, \widehat{P}_m^{-\{2\}})\right) \\ &= \mathbb{E}\left((J(X_1, \widehat{P}_m^{-\{1\}}) - J(X_1, \widehat{P}_m^{-\{1,2\}}))(J(X_2, \widehat{P}_m^{-\{2\}}) - J(X_2, \widehat{P}_m^{-\{1,2\}}))\right) \\ &\quad + 2 \cdot \mathbb{E}\left((J(X_1, \widehat{P}_m^{-\{1\}}) - J(X_1, \widehat{P}_m^{-\{1,2\}}))J(X_2, \widehat{P}_m^{-\{1,2\}})\right) \\ &\quad + \mathbb{E}\left(J(X_1, \widehat{P}_m^{-\{1,2\}})J(X_2, \widehat{P}_m^{-\{1,2\}})\right), \end{aligned}$$

where the first term is $O(m^{-2\beta}) = O(\delta_{m,n}^2)$ by first upper bounding by an adapted version of (25) (replace $\widehat{P}_m$ by $\widehat{P}_m^{-\{1\}}$ and $\widehat{P}_m^{-\{1\}}$ by $\widehat{P}_m^{-\{1,2\}}$), then using conditional independence given $X_3, ..., X_m$ and then $\mathrm{Asm}_4$. For the second term, observe that

$$\begin{aligned} &\mathbb{E}\left((J(X_1, \widehat{P}_m^{-\{1\}}) - J(X_1, \widehat{P}_m^{-\{1,2\}}))J(X_2, \widehat{P}_m^{-\{1,2\}})\right) \\ &= \mathbb{E}\left(\mathbb{E}\left((J(X_1, \widehat{P}_m^{-\{1\}}) - J(X_1, \widehat{P}_m^{-\{1,2\}}))J(X_2, \widehat{P}_m^{-\{1,2\}}) \,\middle|\, X_2, ..., X_m\right)\right) \\ &= \mathbb{E}\left(\mathbb{E}\left(J(X_1, \widehat{P}_m^{-\{1\}}) - J(X_1, \widehat{P}_m^{-\{1,2\}}) \,\middle|\, X_2, ..., X_m\right) \cdot J(X_2, \widehat{P}_m^{-\{1,2\}})\right) = 0. \end{aligned}$$

The last term is 0 by (28) and a conditional independence argument. This concludes the proof. □

**Remark A.4.** *In [93] it is shown that*

$$\widetilde{A}_{2,2}(\widehat{P}_m) = \int\int (\widetilde{I}(x, y, \widehat{P}_m) - \widetilde{I}(x, y, P))\ \mathrm{d}(\widehat{P}_m(x) - P(x))\ \mathrm{d}Q(y) + o_P(m^{-1/2})$$

*(see the last lines in the proof of Lemma 3 in [93]), where $\widetilde{A}$ is the version of $A$ as in the above proof but with the ranks (1). Then, in their Lemma 1 (ii), [93] claim that the term $\int\int (I(x, y, \widehat{P}_m) - I(x, y, P))\ \mathrm{d}(\widehat{P}_m(x) - P(x))\ \mathrm{d}Q(y)$ is also $o_P(m^{-1/2})$ under their regularity conditions. However, the proof of this*

*part of their lemma is omitted. It seems that the authors had a similar proof in mind as for* $\widetilde{A}_{2,1}(\widehat{P}_m, \widehat{Q}_n)$ *which they did prove correctly in the same lemma. Personal communication with the authors could not solve this issue. However, there is a fundamental difference between* $\widetilde{A}_{2,1}$ *and* $\widetilde{A}_{2,2}$ *because of the following: In both cases the integrand* $\widetilde{I}(x, y, \widehat{P}_m)$ *depends on* $\widehat{P}_m$*, but in the former case we integrate over* $\widehat{Q}_n$ *(which yields independent summands conditionally on* $X_1, ..., X_m$*) while in the latter case we integrate with respect to* $\widehat{P}_m$ *such that the summands are* ***not*** *independent (and no corresponding conditioning argument can be made). Therefore, the term* $\widetilde{A}_{2,2}$ *is more difficult to deal with than* $\widetilde{A}_{2,1}$ *despite the fact that it looks simpler. This is the reason why we need Asm₃ and Asm₄.*

**Remark A.5.** *If the ranks* $\widetilde{R}$ *as in (1) are used, then instead of (27) it holds*

$$\begin{aligned}\int\int \widetilde{I}(x,y,P_1)\ \mathrm{d}P(x)\ \mathrm{d}Q(y) &= \frac{1}{2}\int\int \Big(\widetilde{I}(x,y,P_1) + \widetilde{I}(y,x,P_1)\Big)\ \mathrm{d}P(x)\ \mathrm{d}Q(y)\\ &= \frac{1}{2} + \frac{1}{2}\int\int \mathbb{1}_{\{D(x,P_1)=D(y,P_1)\}}\ \mathrm{d}P(x)\ \mathrm{d}Q(y),\end{aligned}$$

*where the second term for* $P_1 = \widehat{P}_m^{-\{1,2\}}$ *is* $o_P(1/\sqrt{m})$ *if (8) holds and* $0$ *for* $P_1 = P$.

*Proof of Theorem 2.3.* The decomposition in (20) in combination with Lemma A.3 shows that the joint asymptotic distribution of $\mathcal{LS}(\widehat{P}_m, \widehat{Q}_n) - 1/2$ and $\mathcal{LS}(\widehat{Q}_n, \widehat{P}_m) - 1/2$ coincide with that of $A_1(\widehat{P}_m, \widehat{Q}_n)$. By Lemma A.1, we have that $A_1(\widehat{Q}_n, \widehat{P}_m) = -A_1(\widehat{P}_m, \widehat{Q}_n)$ a.s. and the result follows. □

*Proof of Corollary 2.5.* By (20) and Lemma A.1 it holds $\mathcal{LS}(\widehat{P}_m, \widehat{Q}_n) + \mathcal{LS}(\widehat{Q}_n, \widehat{P}_m) - 1 = A_2(\widehat{P}_m, \widehat{Q}_n) + A_2(\widehat{Q}_n, \widehat{P}_m)$ a.s., and consequently the assertion follows from Lemma A.3. □

**Remark A.6.** *Instead of considering* $\mathcal{LS}(\widehat{P}_m, \widehat{Q}_n)$*, it is possible to consider*

$$\frac{1}{mn}\sum_{i=1}^{m}\sum_{j=1}^{n}\mathbb{1}_{\{D(X_i,\widehat{P}_m^{-\{i\}})<D(Y_j,\widehat{P}_m^{-\{i\}})\}} + \frac{1}{2}\frac{1}{mn}\sum_{i=1}^{m}\sum_{j=1}^{n}\mathbb{1}_{\{D(X_i,\widehat{P}_m^{-\{i\}})=D(Y_j,\widehat{P}_m^{-\{i\}})\}}$$

*inspired by leave-one-out cross validation. This was implicitely used in the proof of Lemma A.3 and the rest of the proofs go through analogously. Thus, for this modification, Theorem 2.3 and Corollary 2.5 also hold. In practice, this statistic is too computationally expensive.*

*The limit distribution as in Theorem 2.6 (a) for this version of the statistic is now centered, namely is given by* $\frac{1}{2}\chi_1^2 - \frac{1}{2}$.

## A.2. Proof of Theorem 2.1 and Remark 2.4

The depth regions of Tukey depth as in (5) with respect to $\widehat{P}_m$ are given by $[X_{(i)}, X_{(i+1)}) \cup (X_{(m-i)}, X_{(m-i+1)}]$, $i = 0, ..., \lfloor m/2 \rfloor$, (with $X_{(0)} = -\infty$ and

$X_{(m+1)} = \infty$). Thus,

$$\mathbb{E}\left(\sum_{i=1}^{m} \mathbb{1}_{\{D_T(X_{(i)},\widehat{P}_m)=D_T(Y_\ell,\widehat{P}_m)\}} \middle| X_1, ..., X_m\right) \leqslant 4\sum_{i=1}^{m}(X_{(i+1)} - X_{(i)}) \leqslant 4.$$

By an application of the Markov inequality, this yields $\mathcal{LS}(\widehat{P}_m, \widehat{Q}_n) = \widetilde{\mathcal{LS}}(\widehat{P}_m, \widehat{Q}_n) + O_P(1/m)$. Hence, it is sufficient to derive the asymptotic distribution of $\widetilde{\mathcal{LS}}(\widehat{P}_m, \widehat{Q}_n)$. Denote the distribution function corresponding to $P$ resp. $Q$ by $F_P$ resp. $F_Q$. Then, by continuity of $F_P$ and $F_Q$,

$$\begin{aligned} D_T(X, P) &= \min(F_P(X), 1 - F_P(X)) = \min(X, 1 - X) \sim U(0, 1/2), \qquad (29)\\ D_T(Y, Q) &= \min(F_Q(Y), 1 - F_Q(Y)) = \min(2Y, 1 - 2Y) \sim U(0, 1/2),\\ D_T(Y, P) &= \min(F_P(Y), 1 - F_P(Y)) = \min(Y, 1 - Y) = Y \sim U(0, 1/2),\\ D_T(X, Q) &= \min(F_Q(X), 1 - F_P(X)) = \begin{cases} 0 & \text{for } X > 1/2,\\ 1 - 2X & \text{for } 1/4 < X \leqslant 1/2,\\ 2X & \text{for } X \leqslant 1/4. \end{cases} \end{aligned}$$

Thus, with $\widetilde{I}$ as in (19) but with $\widetilde{R}$ rather than $R$, simple calculations lead to

$$\int\int \widetilde{I}(x, y, P)\ \mathrm{d}P(x)\ \mathrm{d}Q(y) = \mathbb{P}(D_T(X, P) \leqslant D_T(Y, P)) = \frac{1}{2}, \qquad (30)$$

$$\int\int \widetilde{I}(y, x, Q)\ \mathrm{d}Q(y)\ \mathrm{d}P(x) = \mathbb{P}(D_T(Y, Q) \leqslant D_T(X, Q)) = \frac{1}{4}. \qquad (31)$$

For the statement of Remark 2.4 we need to calculcate $\sigma^2_{PQ}$ and $\sigma^2_{QP}$ as in Theorem 1 in [93]. Indeed, it follows from the above considerations that

$$\sigma^2_{PQ} = \mathbb{E}\left([2\, D_T(Y, P)]^2\right) - \left(\frac{1}{2}\right)^2 = \frac{1}{12}, \qquad (32)$$

$$\sigma^2_{QP} = \mathbb{E}\left([1 - 2\, D_T(X, P)]^2\right) - \left(\frac{1}{2}\right)^2 = \frac{1}{12}. \qquad (33)$$

A careful check entails that the choice $P = U(0, 1)$, $Q = U(0, 0.5)$ satisfies all the prerequisites of Zuo and He [93, Theorem 1]. In particular, Assumptions (A1)–(A4) of [93] are fulfilled, see Zuo and He [93, Example 2] respectively Proposition 3.3. Note, that for verifying Assumption (A4) of [93], some calculations using the joint density of two order statictics of the $U(0, 1)$-distribution show that there exists a constant $C > 0$ such that $\mathbb{E}((F_Q(X_{(i+1)}) - F_Q(X_{(i)}))^2) \leqslant C/m^2$ for all $i \in \{1, ..., m-1\}$ and then applying Hölder's inequality entails the rate.

To prove Theorem 2.1 (a) we decompose

$$\begin{aligned} &\widetilde{\mathcal{LS}}_T(\widehat{P}_m, \widehat{Q}_n)\\ &= \int\int \widetilde{I}(x, y, P)\ \mathrm{d}P(x)\ \mathrm{d}\widehat{Q}_n(y) \end{aligned}$$

$$
\begin{aligned}
&+ \left( \int\int \widetilde{I}(x,y,\widehat{P}_m) \,\mathrm{d}\widehat{P}_m(x)\,\mathrm{d}Q(y) - \int\int \widetilde{I}(x,y,P)\,\mathrm{d}P(x)\,\mathrm{d}Q(y) \right) \\
&+ \int\int \widetilde{I}(x,y,P)\,\mathrm{d}(\widehat{P}_m(x) - P(x))\,\mathrm{d}(\widehat{Q}_n(y) - Q(y)) \\
&+ \int\int (\widetilde{I}(x,y,\widehat{P}_m) - \widetilde{I}(x,y,P))\,\mathrm{d}\widehat{P}_m(x)\,\mathrm{d}(\widehat{Q}_n(y) - Q(y)) \\
&= \mathrm{K}_1(\widehat{Q}_n;P) + \mathrm{K}_2(\widehat{P}_m;P,Q) + \mathrm{K}_3(\widehat{P}_m,\widehat{Q}_n;P,Q) + \mathrm{K}_4(\widehat{P}_m,\widehat{Q}_n;P,Q), \qquad (34)
\end{aligned}
$$

where $\mathrm{K}_3(\widehat{P}_m,\widehat{Q}_n;P,Q)$ corresponds to $\widetilde{A}_{1,3}$ as in (21) and $\mathrm{K}_4(\widehat{P}_m,\widehat{Q}_n;P,Q)$ corresponds to $\widetilde{A}_{2,1}$ as in (23) but based on the $\widetilde{R}$ depth-ranks and under the given alternatives. Nevertheless, analogous proofs (see also Lemma 1 (i) and (iii) in [93]) show $\mathrm{K}_3(\widehat{P}_m,\widehat{Q}_n;P,Q) = o_P(1/\sqrt{n})$ as well as $\mathrm{K}_4(\widehat{P}_m,\widehat{Q}_n;P,Q) = o_P(1/\sqrt{n})$. Concerning $\mathrm{K}_1(\widehat{Q}_n;P)$ it holds

$$
\begin{aligned}
&\int\int \widetilde{I}(x,y,P)\,\mathrm{d}P(x)\,\mathrm{d}\widehat{Q}_n(y) \\
&\qquad = \frac{1}{n}\sum_{j=1}^{n} \mathbb{P}\left( D_T(X,P) \leqslant D_T(Y_j,P) \mid Y_j \right) = \frac{1}{n}\sum_{j=1}^{n} 2\,Y_j,
\end{aligned}
$$

and by the central limit theorem

$$
\sqrt{\frac{12\cdot m\cdot n}{m+n}} \cdot \left( \mathrm{K}_1(\widehat{Q}_n;P) - \frac{1}{2} \right) \xrightarrow{\mathcal{D}} N(0,\tau).
$$

This term is independent of $\mathrm{K}_2(\widehat{P}_m;P,Q)$ which by (30) fufills

$$
\mathrm{K}_2(\widehat{P}_m;P,Q) = \frac{1}{m}\sum_{i=1}^{m} \mathbb{P}\left( D_T(X_{(i)},\widehat{P}_m) \leqslant D_T(Y,\widehat{P}_m) \middle| X_1,\dots,X_m \right) - \frac{1}{2}.
$$

It holds (irrespective of the distributional assumptions on $P$ and $Q$)

$$
D_T(X_{(i)},\widehat{P}_m) \leqslant D_T(Y,\widehat{P}_m) \iff X_{(\min(i,m-i+1))} \leqslant Y \leqslant X_{(\max(i,m-i+1))}. \qquad (35)
$$

Thus, we obtain for $\mathrm{K}_2(\widehat{P}_m;P,Q)$

$$
\begin{aligned}
&\frac{1}{m}\sum_{i=1}^{m} \mathbb{P}\left( D_T(X_{(i)},\widehat{P}_m) \leqslant D_T(Y,\widehat{P}_m) \middle| X_1,\dots,X_m \right) - \frac{1}{2} \\
&\overset{\text{a.s.}}{=} \frac{2}{m}\sum_{i\leqslant\frac{m}{2}} \mathbb{P}\left( Y \in [X_{(i)},X_{(m-i+1)}] \middle| X_1,\dots,X_m \right) - \frac{1}{2} \\
&= \frac{2}{m}\sum_{i\leqslant\frac{m}{2}} \left( F_Q(X_{(m-i+1)}) - F_Q(X_{(i)}) \right) - \frac{1}{2}
\end{aligned}
$$

$$= \frac{4}{m} \sum_{i \leqslant \frac{m}{2}} \left( \min \left( X_{(m-i+1)}, \frac{1}{2} \right) - \min \left( X_{(i)}, \frac{1}{2} \right) \right) - \frac{1}{2} \tag{36}$$

$$= \frac{8}{m} \sum_{i \leqslant \frac{m}{2}} \min \left( X_{(m-i+1)} - \frac{1}{2}, 0 \right) + \frac{4}{m} \sum_{i=1}^{m} \left( \frac{3}{8} - \min \left( X_i, \frac{1}{2} \right) \right) + O\left( \frac{1}{m} \right). \tag{37}$$

Concerning the second sum in (37) first note that with $Y \sim U(0, 1/2)$

$$\mathbb{E}\left( \min \left( X_i, \frac{1}{2} \right) \right) = \frac{1}{2} \mathbb{P}\left( X_i \geqslant \frac{1}{2} \right) + \mathbb{E}(Y) \, \mathbb{P}\left( X_i \leqslant \frac{1}{2} \right) = \frac{3}{8},$$

and similary $\mathbb{E}(\min (X_i, 1/2)^2) = 1/6$ such that $\mathbb{V}\mathrm{ar}(\min (X_i, 1/2)) = 5/(4^2 \cdot 12)$. Thus, an application of the central limit theorem yields

$$\sqrt{\frac{12 \cdot m \cdot n}{m+n}} \cdot \frac{4}{m} \sum_{i=1}^{m} \left( \frac{3}{8} - \min \left( X_i, \frac{1}{2} \right) \right) \xrightarrow{\mathcal{D}} N(0, 5 \cdot (1-\tau)).$$

We will now complete the proof by showing that the first sum in (36) is $o_P(1/\sqrt{m})$. It holds

$$\left| \sum_{i \leqslant \frac{m}{2}} \min \left( X_{(m-i+1)} - \frac{1}{2}, 0 \right) \right| = \left| \sum_{i \leqslant \frac{m}{2}} \left( X_{(m-i+1)} - \frac{1}{2} \right) \cdot \mathbb{1}_{\{X_{(m-i+1)} \leqslant \frac{1}{2}\}} \right|$$

$$\leqslant \sum_{i \leqslant \frac{m}{2} - m^{3/4}} \mathbb{1}_{\{X_{(m-i+1)} \leqslant \frac{1}{2}\}} + m^{3/4} \cdot \left| X_{(\lceil \frac{m}{2} \rceil)} - \frac{1}{2} \right|,$$

where the last summand is $O_P(m^{1/4}) = o_P(\sqrt{m})$ due to the asymptotic normality of the sample median. For the sum we obtain the result by Markov's inequality on noting

$$\mathbb{E}\left( \left| \sum_{i \leqslant \frac{m}{2} - m^{3/4}} \mathbb{1}_{\{X_{(m-i+1)} \leqslant \frac{1}{2}\}} \right| \right) = \sum_{i \leqslant \frac{m}{2} - m^{3/4}} \mathbb{P}\left( X_{(m-i+1)} \leqslant \frac{1}{2} \right)$$

$$= \sum_{i \leqslant \frac{m}{2} - m^{3/4}} \mathbb{P}\left( \sum_{j=1}^{m} \left( \mathbb{1}_{\{X_j \leqslant \frac{1}{2}\}} - \frac{1}{2} \right) \geqslant m - i + 1 - \left\lfloor \frac{m}{2} \right\rfloor \right)$$

$$\leqslant C \cdot \sum_{i \leqslant \frac{m}{2} - m^{3/4}} \frac{m}{m - i + 1 - \lfloor \frac{m}{2} \rfloor} = O(m) \int_{m^{3/4}}^{m/2} \frac{1}{x^2} \, \mathrm{d}x$$

$$= O(m) \, m^{-3/4} = o(\sqrt{m}).$$

So, the proof of (a) is completed.

For (b), we use exactly the same decomposition as in (34) but with $P$ and $Q$ (respectively $\widehat{P}_m$ and $\widehat{Q}_n$) switched. In that case, we obtain

$$\sqrt{\frac{12 \cdot m \cdot n}{m+n}} \left( \mathrm{K}_1(\widehat{P}_m; Q) - \frac{1}{4} \right) \xrightarrow{\mathcal{D}} N\left( 0, \frac{5}{4} \cdot (1-\tau) \right), \tag{38}$$

due to

$$\int\int \widetilde{I}(y,x,Q)\,\mathrm{d}\widehat{P}_m(x)\,\mathrm{d}Q(y)$$
$$= \frac{1}{m}\sum_{i=1}^{m}\mathbb{E}\left(\mathbb{1}_{\{D_T(Y,Q)\leqslant D_T(X_i,Q)\}}\,\middle|\, X_1,...,X_m\right) = \frac{2}{m}\sum_{i=1}^{m} D_T(X_i,Q).$$

By using the same arguments as in (a) both $\mathrm{K}_3(\widehat{Q}_n,\widehat{P}_m;Q,P)$ and $\mathrm{K}_4(\widehat{Q}_n,\widehat{P}_m;Q,P)$ are $o_P(1/\sqrt{n})$.

Analogously to (36) it holds

$$\mathrm{K}_2(\widehat{Q}_n;Q,P) \overset{\text{a.s.}}{=} \frac{2}{n}\sum_{j\leqslant\frac{n}{2}}\left(Y_{(n-j+1)}-Y_{(j)}\right)-\frac{1}{4} = \frac{2}{n}\sum_{j=1}^{n}\left|Y_j - Y_{(\frac{n}{2})}\right| - \frac{1}{4}.$$

[7, Theorem 2.5] entails that

$$\frac{1}{\sqrt{n}}\left(\sum_{j=1}^{n}\left|Y_j - Y_{(\frac{n}{2})}\right| - \frac{1}{8}\right) \xrightarrow{\mathcal{D}} N\left(0,\frac{1}{192}\right).$$

Thus,

$$\sqrt{\frac{12mn}{m+n}}\cdot\left(\frac{2}{n}\sum_{j=1}^{n}\left|Y_j - Y_{(\frac{n}{2})}\right| - \frac{1}{4}\right) \xrightarrow{\mathcal{D}} N\left(0,\frac{1}{4}\cdot\tau\right)$$

and by combining this with (38) due to the independence of these two terms, Theorem 2.1 (b) follows.

### *A.3. Proof of Theorems 2.6 and 2.7*

In this section, we make use of the following square-bracket notation (later also for random variables $U,V$)

$$X_{[i]} = X_{(\min(i,m-i+1))},\ X^{[i]} = X_{(\max(i,m-i+1))}.$$

Furthermore, we introduce

$$U_i := F(X_i) \sim U(0,1), \qquad V_j := F(Y_j) \sim U(0,1) \tag{39}$$

and denote $V_{(a)} = V_{(\lceil a\rceil)}$ for $a\in(0,n]$ (and analogously for $U$).

The proof of Theorem 2.6 is based on the following Lemmas.

**Lemma A.7.** *Let $U_1,...,U_m$ be a sample of iid $U(0,1)$-distributed random variables. Then it holds:*

*(a) $U_{(k)}\sim Beta(k,m-k+1)$ for each $k\in\{1,...,m\}$ and for any sample of independent and identically distributed random variables $X_1,...,X_m$ with continuous distribution function $F$, we have that $F(X_{(k)}) \overset{\mathcal{D}}{=} U_{(k)}$ for each $k\in\{1,...,m\}$.*

*(b)* $U_{(b)} - U_{(a)} \sim Beta(b-a, m-(b-a)+1)$ *for any* $a, b \in \{1, ..., m\}$, $a < b$.
*(c)* $\mathbb{E}(U_{(k)}) = k/(m+1)$ *for each* $k \in \{1, ..., m\}$.
*(d)* $U_{(m)} \xrightarrow{P} 1$ *as well as* $U_{(1)} \xrightarrow{P} 0$ *for* $m \to \infty$.
*(e) For* $m \to \infty$,

$$\sqrt{m}\left(U_{(\frac{m}{2})} - \frac{1}{2}\right) \xrightarrow{\mathcal{D}} N\left(0, \frac{1}{4}\right).$$

*(f) For the joint distribution*

$$\left(U_{(1)}, ..., U_{(m)}\right) \stackrel{\mathcal{D}}{=} \left(\frac{Z_1}{Z_1 + ... + Z_{m+1}}, ..., \frac{Z_1 + ... + Z_m}{Z_1 + ... + Z_{m+1}}\right),$$

*where* $Z_1, ..., Z_{m+1} \sim Exp(1)$ *independent.*

*(g)*

$$\left(U_{(1)}, ..., U_{(m)}\right)|_{U_{(m/2)}=u} \stackrel{\mathcal{D}}{=} \left(U'_{(1)}, ..., U'_{(\lceil \frac{m}{2} - 1 \rceil)}, u, U''_{(1)}, ..., U''_{(\lfloor \frac{m}{2} \rfloor)}\right),$$

*where* $\left(U_{(1)}, ..., U_{(m)}\right)|_{U_{(m/2)}=u}$ *is the conditional distribution of* $\left(U_{(1)}, ..., U_{(m)}\right)$ *given* $U_{(\frac{m}{2})} = u$ *for some* $u \in (0,1)$, *and independent random variables* $U'_i(u)$, $i < \frac{m}{2}$, $U''_j(u)$, $j \leqslant \frac{m}{2}$, *(also independent of each other) with* $U'_i(u) \sim U(0, u)$ *and* $U''_j(u) \sim U(u, 1)$.

All of the above assertions are well known in the nonparametric statistics literature on order statistics. For instance, the first part of (a) can be found in Ahsanullah, Nevzorov and Shakil [1, Remark 2.1], the second one in Scheffé and Tukey [72, Section 4], (b) is given in David and Nagaraja [23, Example 2.3], (c) is presented in Gupta and Nadarajah [42, Chapter 2, Section IV], (d) is a direct consequence of Ahsanullah, Nevzorov and Shakil [1, Exercise 11.4], (e) equals [73, Corollary A in Section 2.3.3], (f) and (g) correspond to Ahsanullah, Nevzorov and Shakil [1, Example 4.6 respectively Example 5.1].

**Lemma A.8.** *Under the assumptions of Theorem 2.6 it holds*

$$-\frac{mn}{m+n}\left(\widetilde{\mathcal{LS}}_T(\widehat{P}_m, \widehat{Q}_n) + \widetilde{\mathcal{LS}}_T(\widehat{Q}_n, \widehat{P}_m) - 1\right) = \frac{2mn}{m+n}\left(U_{(\frac{m}{2})} - V_{(\frac{n}{2})}\right)^2 + o_P(1) \xrightarrow{\mathcal{D}} \frac{1}{2}\chi^2(1).$$

*Proof.* By (35) it holds under $H_0$ with $F$ denoting the distribution function of $X$ and $Y$

$$\begin{aligned} D_T(X_{(i)}, \widehat{P}_m) \leqslant D_T(Y_j, \widehat{P}_m) \quad &\Longleftrightarrow \quad X_{[i]} \leqslant Y_j \leqslant X^{[i]} \\ &\overset{\text{a.s.}}{\Longleftrightarrow} \quad F(X_{[i]}) \leqslant F(Y_j) \leqslant F(X^{[i]}). \end{aligned} \tag{40}$$

The last equivalence holds except on the set of mass zero where $F(X_{(i)}) = F(Y_j)$ but at the same time $X_{[i]} > Y_j$ or $F(Y_j) = F(X^{[i]})$ at the same time as $Y_j > X^{[i]}$. Equivalent transformations also hold for $D_T(Y_{(j)}, \widehat{Q}_n) \leqslant D_T(X_i, \widehat{Q}_n)$.

With $U_i, V_j$ as in (39) this implies

$$
\begin{aligned}
&\widetilde{\mathcal{LS}}_T(\widehat{P}_m, \widehat{Q}_n) + \widetilde{\mathcal{LS}}_T(\widehat{Q}_n, \widehat{P}_m) \\
&= \frac{1}{mn}\sum_{i=1}^{m}\sum_{j=1}^{n}\Big(\mathbb{1}_{\{D_T(X_i,\widehat{P}_m)\leqslant D_T(Y_j,\widehat{P}_m)\}} + \mathbb{1}_{\{D_T(Y_j,\widehat{Q}_n)\leqslant D_T(X_i,\widehat{Q}_n)\}}\Big) \\
&\overset{\text{a.s.}}{=} \frac{1}{mn}\sum_{i=1}^{m}\sum_{j=1}^{n}\mathbb{1}_{\{U_{[i]}\leqslant V_{(j)}\leqslant U^{[i]}\}} + \frac{1}{mn}\sum_{i=1}^{m}\sum_{j=1}^{n}\mathbb{1}_{\{V_{[j]}\leqslant U_{(i)}\leqslant V^{[j]}\}}. \qquad (41)
\end{aligned}
$$

By symmetry

$$
\begin{aligned}
&\frac{1}{mn}\sum_{i=1}^{m}\sum_{j=1}^{n}\mathbb{1}_{\{U_{[i]}\leqslant V_{(j)}\leqslant U^{[i]}\}} \\
&= \frac{2}{mn}\sum_{i\leqslant\frac{m}{2}}\sum_{j=1}^{n}\mathbb{1}_{\{U_{[i]}\leqslant V_{(j)}\leqslant U^{[i]}\}} + \mathbb{1}_{\{m \text{ is odd}\}}\frac{1}{mn}\sum_{j=1}^{n}\mathbb{1}_{\left\{V_{(j)}=U_{\left(\frac{m}{2}\right)}\right\}} \\
&\overset{\text{a.s.}}{=} \frac{2}{mn}\sum_{i\leqslant\frac{m}{2}}\sum_{j=1}^{n}\mathbb{1}_{\{U_{[i]}\leqslant V_{(j)}\leqslant U^{[i]}\}}. \qquad (42)
\end{aligned}
$$

Furthermore,

$$
\begin{aligned}
&\frac{2}{mn}\sum_{i\leqslant\frac{m}{2}}\sum_{j=1}^{n}\mathbb{1}_{\{U_{[i]}\leqslant V_{(j)}\leqslant U^{[i]}\}} \\
&\qquad= \frac{2}{mn}\sum_{i\leqslant\frac{m}{2}}\sum_{j\leqslant\frac{n}{2}}\Big(\mathbb{1}_{\{U_{[i]}\leqslant V_{[j]}\leqslant U^{[i]}\}} + \mathbb{1}_{\{U_{[i]}\leqslant V^{[j]}\leqslant U^{[i]}\}}\Big) \\
&\qquad\quad + \mathbb{1}_{\{n \text{ is odd}\}}\frac{2}{mn}\sum_{i\leqslant\frac{m}{2}}\mathbb{1}_{\left\{U_{[i]}\leqslant V_{\left(\frac{n}{2}\right)}\leqslant U^{[i]}\right\}}. \qquad (43)
\end{aligned}
$$

Analogous assertions also hold for the second summand in (41), i.e.

$$
\begin{aligned}
&\frac{2}{mn}\sum_{i=1}^{m}\sum_{j\leqslant\frac{n}{2}}\mathbb{1}_{\{V_{[j]}\leqslant U_{(i)}\leqslant V^{[j]}\}} \\
&\qquad= \frac{2}{mn}\sum_{i\leqslant\frac{m}{2}}\sum_{j\leqslant\frac{n}{2}}\Big(\mathbb{1}_{\{V_{[j]}\leqslant U_{[i]}\leqslant V^{[j]}\}} + \mathbb{1}_{\{V_{[j]}\leqslant U^{[i]}\leqslant V^{[j]}\}}\Big) \\
&\qquad\quad + \mathbb{1}_{\{m \text{ is odd}\}}\frac{2}{mn}\sum_{j\leqslant\frac{n}{2}}\mathbb{1}_{\left\{V_{[j]}\leqslant U_{\left(\frac{m}{2}\right)}\leqslant V^{[j]}\right\}}. \qquad (44)
\end{aligned}
$$

We will now show that $\frac{2}{n}\sum_{j\leqslant n/2}\big(\mathbb{1}_{\{V_{[j]}\leqslant U_{(m/2)}\leqslant V^{[j]}\}} - 1\big)$ respectively $\frac{2}{m}\sum_{i\leqslant m/2}\big(\mathbb{1}_{\{U_{[i]}\leqslant V_{(n/2)}\leqslant U^{[i]}\}}-1\big)$ are asymptotically negligible in a P-stochastic sense for any $m$ respectively $n$ (where we will use the case for even integers later). Observe that

$$\frac{2}{n}\sum_{j\leqslant \frac{n}{2}}\left(\mathbb{1}_{\left\{V_{[j]}\leqslant U_{\left(\frac{m}{2}\right)}\leqslant V^{[j]}\right\}}-1\right)$$
$$=-\frac{2}{n}\sum_{j\leqslant \frac{n}{2}}\mathbb{1}_{\left\{U_{\left(\frac{m}{2}\right)}<V_{(j)}\right\}}-\frac{2}{n}\sum_{j\leqslant \frac{n}{2}}\mathbb{1}_{\left\{U_{\left(\frac{m}{2}\right)}>V_{(n-j+1)}\right\}}, \quad (45)$$

such that

$$\left|\frac{2}{n}\sum_{j\leqslant \frac{n}{2}}\left(\mathbb{1}_{\left\{V_{[j]}\leqslant U_{\left(\frac{m}{2}\right)}\leqslant V^{[j]}\right\}}-1\right)\right|$$
$$=\frac{2}{n}\sum_{j\leqslant \frac{n}{2}}\mathbb{1}_{\left\{U_{\left(\frac{m}{2}\right)}<V_{(j)}\right\}}+\frac{2}{n}\sum_{j\leqslant \frac{n}{2}}\mathbb{1}_{\left\{U_{\left(\frac{m}{2}\right)}>V_{(n-j+1)}\right\}}. \quad (46)$$

We will now prove the assertion for the first summand on the right hand side of (46), the assertion for the second sum can be dealt with analogously. By splitting the sum in summands with $j \leqslant \lceil n/2 - n^{7/8}\rceil$ and the complement we get

$$\begin{aligned}
&\mathbb{E}\left(\frac{2}{n}\sum_{j\leqslant \frac{n}{2}}\mathbb{1}_{\left\{U_{\left(\frac{m}{2}\right)}<V_{[j]}\right\}}\right)\leqslant \mathbb{P}\left(U_{\left(\frac{m}{2}\right)}<V_{\left(\frac{n}{2}-n^{7/8}\right)}\right)+O(n^{-1/8})\\
&\leqslant \mathbb{P}\left(U_{\left(\frac{m}{2}\right)}<V_{\left(\frac{n}{2}-n^{7/8}\right)}, U_{\left(\frac{m}{2}\right)}\geqslant \frac{1}{2}-n^{-1/4}\right)\\
&\quad+\mathbb{P}\left(\sqrt{m}\left(U_{\left(\frac{m}{2}\right)}-\frac{1}{2}\right)\leqslant -\frac{m^{1/2}}{n^{1/4}}\right)+o(1)\\
&\leqslant \mathbb{P}\left(V_{\left(\frac{n}{2}-n^{7/8}\right)}\geqslant \frac{1}{2}-n^{-1/4}\right)+o(1),
\end{aligned}$$

where we used Lemma A.7 (e) in the last line. We can now conclude with an application of the Chebychev inequality

$$\begin{aligned}
&\mathbb{P}\left(V_{\left(\frac{n}{2}-n^{7/8}\right)}\geqslant \frac{1}{2}-n^{-1/4}\right)=\mathbb{P}\left(\sum_{j=1}^{n}\mathbb{1}_{\left\{V_j<\frac{1}{2}-n^{-1/4}\right\}}<\frac{n}{2}-n^{7/8}\right)\\
&\leqslant \mathbb{P}\left(\left|\sum_{j=1}^{n}\left(\mathbb{1}_{\left\{V_j<\frac{1}{2}-n^{-1/4}\right\}}-\frac{1}{2}+n^{-1/4}\right)\right|\geqslant n^{7/8}\left(1-n^{-1/8}\right)\right)\leqslant O(n^{-3/4})\\
&=o(1).
\end{aligned}$$

Thus, we have shown that

$$\frac{2}{mn}\sum_{j\leqslant \frac{n}{2}}\mathbb{1}_{\left\{U_{\left(\frac{m}{2}\right)}<V_{(j)}\right\}}=o_P\left(\frac{1}{m}\right), \quad \frac{2}{mn}\sum_{j\leqslant \frac{n}{2}}\mathbb{1}_{\left\{U_{\left(\frac{m}{2}\right)}>V_{(n-j+1)}\right\}}=o_P\left(\frac{1}{m}\right). \quad (47)$$

Analogously, we get

$$\frac{2}{mn}\sum_{i\leqslant\frac{m}{2}}\mathbb{1}_{\left\{V_{\left(\frac{n}{2}\right)}<U_{(i)}\right\}}=o_P\left(\frac{1}{n}\right),\quad \frac{2}{mn}\sum_{i\leqslant\frac{m}{2}}\mathbb{1}_{\left\{V_{\left(\frac{n}{2}\right)}>U_{(m-i+1)}\right\}}=o_P\left(\frac{1}{n}\right). \tag{48}$$

Firstly, note that

$$1=\frac{2}{mn}\cdot\left(2\left\lfloor\frac{m}{2}\right\rfloor\left\lfloor\frac{n}{2}\right\rfloor+\mathbb{1}_{\{m\text{ is odd}\}}\left\lfloor\frac{n}{2}\right\rfloor+\mathbb{1}_{\{n\text{ is odd}\}}\left\lfloor\frac{m}{2}\right\rfloor+\frac{\mathbb{1}_{\{m,n\text{ are both odd}\}}}{2}\right)$$

as well as

$$\begin{aligned}\mathbb{1}_{\{U_{[i]}\leqslant V_{[j]}\leqslant U^{[i]}\}}+\mathbb{1}_{\{V_{[j]}<U_{[i]}\leqslant V^{[j]}\}}-1&=\mathbb{1}_{\{V_{[j]}\leqslant U^{[i]}\}}+\mathbb{1}_{\{U_{[i]}\leqslant V^{[j]}\}}-2\\&=-\mathbb{1}_{\{V_{[j]}>U^{[i]}\}}-\mathbb{1}_{\{U_{[i]}>V^{[j]}\}}\end{aligned} \tag{49}$$

and analogously

$$\mathbb{1}_{\{U_{[i]}\leqslant V^{[j]}\leqslant U^{[i]}\}}+\mathbb{1}_{\{V_{[j]}\leqslant U^{[i]}<V^{[j]}\}}-1=-\mathbb{1}_{\{V^{[j]}<U_{[i]}\}}-\mathbb{1}_{\{U^{[i]}<V_{[j]}\}}.$$

Consequently, we get by (41)–(46) as well as (47) and (48)

$$\begin{aligned}&\widetilde{\mathcal{LS}}_T(\widehat{P}_m,\widehat{Q}_n)+\widetilde{\mathcal{LS}}_T(\widehat{Q}_n,\widehat{P}_m)-1\\&\overset{\text{a.s.}}{=}-\frac{4}{mn}\sum_{i\leqslant\frac{m}{2}}\sum_{j<\frac{n}{2}}\mathbb{1}_{\{V_{(j)}>U_{(m-i+1)}\}}-\frac{4}{mn}\sum_{i<\frac{m}{2}}\sum_{j\leqslant\frac{n}{2}}\mathbb{1}_{\{U_{(i)}>V_{(n-j+1)}\}}\\&\quad+o_P\left(\frac{1}{n}+\frac{1}{m}\right),\end{aligned} \tag{50}$$

where another application of (47) and (48) permits the modification in the indices of the sums for $m$ respectively $n$ even.

Analogous considerations noting that the term in (45) as well as (49) are negative yield that

$$\widetilde{\mathcal{LS}}_T(\widehat{P}_m,\widehat{Q}_n)+\widetilde{\mathcal{LS}}_T(\widehat{Q}_n,\widehat{P}_m)-1\leqslant 0\quad\text{a.s.} \tag{51}$$

We will now show that the main terms in (50) are dominated by functionals of the sample medians

$$U_{\left(\frac{m}{2}\right)}=U_{\left(\lceil\frac{m}{2}\rceil\right)},\qquad V_{\left(\frac{n}{2}\right)}=V_{\left(\lceil\frac{n}{2}\rceil\right)}.$$

We can conclude for the first two terms in (50)

$$\mathbb{P}\left(\frac{mn}{m+n}\cdot\frac{4}{mn}\left(\sum_{i<\frac{m}{2}}\sum_{j\leqslant\frac{n}{2}}\mathbb{1}_{\{U_{(i)}>V_{(n-j+1)}\}}+\sum_{i\leqslant\frac{m}{2}}\sum_{j<\frac{n}{2}}\mathbb{1}_{\{V_{(j)}>U_{(m-i+1)}\}}\right)\leqslant z\right)$$

$$= \mathbb{E}\left(\mathbb{P}\left(\frac{mn}{m+n}\cdot\frac{4}{mn}\left(\sum_{i<\frac{m}{2}}\sum_{j\leqslant\frac{n}{2}}\mathbb{1}_{\{U_{(i)}>V_{(n-j+1)}\}}\right.\right.\right.$$
$$\left.\left.\left.+\sum_{i\leqslant\frac{m}{2}}\sum_{j<\frac{n}{2}}\mathbb{1}_{\{V_{(j)}>U_{(m-i+1)}\}}\right)\leqslant z\,\middle|\,U_{(\frac{m}{2})},V_{(\frac{n}{2})}\right)\right).$$

Let $(u_m)_{m\in\mathbb{N}},(v_n)_{n\in\mathbb{N}}$ be two real-valued sequences with

$$\lim_{m\to\infty}u_m=\lim_{n\to\infty}v_n=\frac{1}{2}\text{ and }|u_m-v_n|=o\left((m+n)^{-1/3}\right). \tag{52}$$

From Lemma A.7 (g) (applied jointly to the independent sequences $\{U_j\}$ as well as $\{V_j\}$) we obtain with $U_i'\sim U(0,u_m)$ independent of $V_j''\sim U(v_n,1)$ that for $u_m>v_n$ (otherwise the term is a.s. equal to zero)

$$\mathbb{V}\mathrm{ar}\left(\frac{mn}{m+n}\cdot\frac{4}{mn}\sum_{i<\frac{m}{2}}\sum_{j\leqslant\frac{n}{2}}\mathbb{1}_{\{U_{(i)}>V_{(n-j+1)}\}}\,\middle|\,U_{(\frac{m}{2})}=u_m,V_{(\frac{n}{2})}=v_n\right)$$
$$=\mathbb{V}\mathrm{ar}\left(\frac{mn}{m+n}\cdot\frac{4}{mn}\sum_{i<\frac{m}{2}}\sum_{j\leqslant\frac{n}{2}}\mathbb{1}_{\{U_i'>V_j''\}}\right). \tag{53}$$

Straightforward calculation yields

$$\mathbb{E}\left(\mathbb{1}_{\{U_i'>V_j''\}}\right)=\frac{(u_m-v_n)_+^2}{2u_m(1-v_n)}=o\left((m+n)^{-2/3}\right) \tag{54}$$

due to (52), where $a_+=\max(a,0)$ and $a_+^2=(\max(a,0))^2$. This implies

$$\mathbb{V}\mathrm{ar}\left(\mathbb{1}_{\{U_i'>V_j''\}}\right)\leqslant\mathbb{E}\left(\mathbb{1}_{\{U_i'>V_j''\}}\right)=o\left((m+n)^{-2/3}\right). \tag{55}$$

Moreover we obtain for $i\neq k$ and $j\neq l$ by independence

$$\mathbb{C}\mathrm{ov}\left(\mathbb{1}_{\{U_i'>V_j''\}},\mathbb{1}_{\{U_k'>V_l''\}}\right)=0 \tag{56}$$

and for $i\neq k$ and $j=l$ by straightforward calculation

$$\begin{aligned}&\mathbb{C}\mathrm{ov}\left(\mathbb{1}_{\{U_i'>V_j''\}},\mathbb{1}_{\{U_k'>V_j''\}}\right)\\&=\mathbb{E}\left(\mathbb{1}_{\{U_i'>V_j''\}}\mathbb{1}_{\{U_k'>V_j''\}}\right)-\mathbb{E}\left(\mathbb{1}_{\{U_i'>V_j''\}}\right)\mathbb{E}\left(\mathbb{1}_{\{U_k'>V_j''\}}\right)\\&=\frac{(u_m-v_n)_+^3}{3u_m^2(1-v_n)}+o\left((m+n)^{-4/3}\right)=o\left((m+n)^{-1}\right)\end{aligned} \tag{57}$$

with (52). Analogously, for $i=k$ and $j\neq l$

$$\mathbb{C}\mathrm{ov}\left(\mathbb{1}_{\{U_i'>V_j''\}},\mathbb{1}_{\{U_i'>Y_l''\}}\right)=o\left((m+n)^{-1}\right). \tag{58}$$

Transforming the expression in (53) into a sum of variance and covariance expressions and then plugging in the derived rates respectively identities (55)–(58) delivers

$$\mathbb{V}\mathrm{ar}\left(\frac{mn}{m+n}\cdot\frac{4}{mn}\sum_{i<\frac{m}{2}}\sum_{j\leqslant\frac{n}{2}}\mathbb{1}_{\{U_{(i)}>V_{(n-j+1)}\}}\middle| U_{(\frac{m}{2})}=u_m, V_{(\frac{n}{2})}=v_n\right)=o(1).$$

Then, by Lemma A.7 (g), (54) and the Chebychev inequality it holds for any $\epsilon>0$

$$\mathbb{P}\left(\left|\frac{mn}{m+n}\cdot\frac{4}{mn}\sum_{i<\frac{m}{2}}\sum_{j\leqslant\frac{n}{2}}\left[\mathbb{1}_{\{U_{(i)}>V_{(n-j+1)}\}}-\frac{(u_m-v_n)_+^2}{2u_m\,(1-v_n)}\right]\right|>\epsilon\middle| U_{(\frac{m}{2})}=u_m, V_{(\frac{n}{2})}=v_n\right)=o(1).$$

By Lemma A.7 (e) the sequences $u_m=U_{(\frac{m}{2})}$ and $v_n=V_{(\frac{n}{2})}$ fulfill (52) in a P-stochastic sense. So, an application of the subsequence principle yields

$$\mathbb{P}\left(\left|\frac{mn}{m+n}\cdot\frac{4}{mn}\sum_{i<\frac{m}{2}}\sum_{j\leqslant\frac{n}{2}}\left[\mathbb{1}_{\{U_{(i)}>V_{(n-j+1)}\}}-\frac{\left(U_{(\frac{m}{2})}-V_{(\frac{n}{2})}\right)_+^2}{2U_{(\frac{m}{2})}\left(1-V_{(\frac{n}{2})}\right)}\right]\right|>\epsilon\middle| U_{(\frac{m}{2})}, V_{(\frac{n}{2})}\right)=o_P(1).$$

Since the expression on the left hand side is bounded from above by 1, we get by an application of the dominated convergence theorem (after taking the expectation of that expression)

$$\frac{mn}{m+n}\cdot\frac{4}{mn}\sum_{i<\frac{m}{2}}\sum_{j\leqslant\frac{n}{2}}\left[\mathbb{1}_{\{U_{(i)}>V_{(n-j+1)}\}}-\frac{\left(U_{(\frac{m}{2})}-V_{(\frac{n}{2})}\right)_+^2}{2\cdot U_{(\frac{m}{2})}\cdot\left(1-V_{(\frac{n}{2})}\right)}\right]=o_P(1).$$

By Lemma A.7 (e) it holds $\frac{mn}{m+n}\left(U_{(\frac{m}{2})}-V_{(\frac{n}{2})}\right)_+^2=O_P(1)$, $U_{(m/2)}\xrightarrow{P}1/2$ as well as $V_{(n/2)}\xrightarrow{P}1/2$. So, it follows

$$\frac{mn}{m+n}\cdot\frac{4}{mn}\sum_{i<\frac{m}{2}}\sum_{j\leqslant\frac{n}{2}}\mathbb{1}_{\{U_{(i)}>V_{(n-j+1)}\}}=\frac{2\,mn}{m+n}\left(U_{(\frac{m}{2})}-V_{(\frac{n}{2})}\right)_+^2+o_P(1).$$

Analogous arguments can be used to derive for the second sum of indicators in (50)

$$\frac{mn}{m+n}\cdot\frac{4}{mn}\sum_{i\leqslant\frac{m}{2}}\sum_{j<\frac{n}{2}}\mathbb{1}_{\{V_{(j)}>U_{(m-i+1)}\}}=\frac{2\,mn}{m+n}\left(V_{(\frac{n}{2})}-U_{(\frac{m}{2})}\right)_+^2+o_P(1).$$

The assertion now follows by

$$\left(U_{(\frac{m}{2})} - V_{(\frac{n}{2})}\right)_+^2 + \left(V_{(\frac{n}{2})} - U_{(\frac{m}{2})}\right)_+^2 = \left(U_{(\frac{m}{2})} - V_{(\frac{n}{2})}\right)^2$$

and Lemma A.7 (assertion (e)). □

**Lemma A.9.** *Under the assumptions of Theorem 2.6 it holds*

$$\frac{mn}{m+n} \cdot \frac{1}{2} \cdot \frac{1}{mn} \cdot \sum_{i=1}^{m} \sum_{j=1}^{n} \left( \mathbb{1}_{\left\{D_T(X_i,\widehat{P}_m)=D_T(Y_j,\widehat{P}_m)\right\}} + \mathbb{1}_{\left\{D_T(Y_j,\widehat{Q}_n)=D_T(X_i,\widehat{Q}_n)\right\}} \right) \xrightarrow{P} 1.$$

*Proof.* Firstly, note that with the same notation as above for $2i \neq m+1$ (if $m$ odd)

$$\begin{aligned} &D_T(X_{(i)}, \widehat{P}_m) = D_T(Y_j, \widehat{P}_m) \\ &\iff X_{(\min(i,m-i+1))} \leqslant Y_j < X_{(\min(i,m-i+1)+1)} \\ &\qquad\qquad \text{or} \quad X_{(\max(i,m-i+1)-1)} < Y_j \leqslant X_{(\max(i,m-i+1))} \\ &\overset{\text{a.s.}}{\iff} U_{(\min(i,m-i+1))} \leqslant V_j < U_{(\min(i,m-i+1)+1)} \\ &\qquad\qquad \text{or} \quad U_{(\max(i,m-i+1)-1)} \leqslant V_j < U_{(\max(i,m-i+1))}. \end{aligned}$$

For $2i \neq m+1$ (if $m$ odd) $D_T(X_{(m+1)/2}, \widehat{P}_m) \neq D_T(Y_j, \widehat{P}_m)$ almost surely. Consequently,

$$\sum_{i=1}^{m} \mathbb{1}_{\left\{D_T(X_i,\widehat{P}_m)=D_T(Y_j,\widehat{P}_m)\right\}} \overset{\text{a.s.}}{=} 2 \sum_{i=1}^{m-1} \mathbb{1}_{\{U_{(i)} \leqslant V_j < U_{(i+1)}\}} = 2\,\mathbb{1}_{\{U_{(1)} \leqslant V_j < U_{(m)}\}}.$$

Furthermore,

$$\begin{aligned} &\mathbb{E}\left[ \left( \frac{1}{n} \sum_{j=1}^{n} \left( \mathbb{1}_{\{U_{(1)} \leqslant V_j < U_{(m)}\}} - (U_{(m)} - U_{(1)}) \right) \right)^2 \right] \\ &= \mathbb{E}\left( \mathbb{V}\mathrm{ar}\left( \frac{1}{n} \sum_{j=1}^{n} \mathbb{1}_{\{U_{(1)} \leqslant V_j < U_{(m)}\}} \,\middle|\, U_{(1)}, U_{(m)} \right) \right) \leqslant \frac{1}{n} = o(1). \end{aligned}$$

Since $U_{(m)} \xrightarrow{P} 1$ and $U_{(1)} \xrightarrow{P} 0$, cf. Lemma A.7 (d), it yields

$$\begin{aligned} &\frac{mn}{m+n} \cdot \frac{1}{2} \cdot \frac{1}{mn} \sum_{i=1}^{m} \sum_{j=1}^{n} \mathbb{1}_{\left\{D_T(X_i,\widehat{P}_m)=D_T(Y_j,\widehat{P}_m)\right\}} \\ &\qquad \overset{\text{a.s.}}{=} \frac{m}{m+n} \frac{1}{n} \sum_{j=1}^{n} \mathbb{1}_{\{U_{(1)} \leqslant V_j < U_{(m)}\}} \xrightarrow{P} 1 - \tau. \end{aligned} \tag{59}$$

Analogously we get

$$\frac{mn}{m+n}\cdot\frac{1}{2}\cdot\frac{1}{mn}\sum_{i=1}^{m}\sum_{j=1}^{n}\mathbb{1}_{\left\{D_T(Y_j,\widehat{Q}_n)=D_T(X_i,\widehat{Q}_n)\right\}}\xrightarrow{P}\tau, \tag{60}$$

such that the assertion follows. □

**Lemma A.10.** *It holds*

$$\sqrt{n+1}\cdot\begin{pmatrix}\frac{1}{n+1}\sum_{j=1}^{n}\left(\min(V_j,1-V_j)-\frac{1}{4}\right)\\ V_{(\frac{n}{2})}-\frac{1}{2}\end{pmatrix}\xrightarrow{\mathcal{D}}N(\mathbf{0},\Sigma),$$

*where* $\Sigma$ *denotes a non-singular diagonal-matrix.*

*Proof.* We make use of the methods outlined in [33].

$$\sum_{j\leqslant\frac{n+1}{2}}\frac{j}{n+1}+\sum_{j>\frac{n+1}{2}}\frac{n+1-j}{n+1}=\frac{n}{4}+O(1). \tag{61}$$

Decompose

$$\min\left(V_{(j)},1-V_{(j)}\right)=\begin{cases}V_{(j)}-R_j & \text{for } j\leqslant\dfrac{n+1}{2},\\ 1-V_{(j)}-R_j & \text{for } j>\dfrac{n+1}{2},\end{cases}$$

where $R_j\geqslant 0$ by construction. Because $\mathbb{E}(V_{(j)})=j/(n+1)$, cf. Lemma A.7 (a) and (c), as well as $\mathbb{E}(\min(V_j,1-V_j))=1/4$ (as $\min(V_j,1-V_j)\sim U(0,1/2)$), we obtain by (61) that $\mathbb{E}|\sum_{j=1}^{n}R_j|=O(1)$. Thus, $\frac{1}{\sqrt{n}}\sum_{j=1}^{n}R_j=o_P(1)$ and consequently

$$\begin{aligned}&\frac{1}{\sqrt{n+1}}\sum_{j=1}^{n}\left(\min\left(V_j,1-V_j\right)-\frac{1}{4}\right)\\ &\qquad=\frac{1}{\sqrt{n+1}}\cdot\left(\sum_{j\leqslant\frac{n+1}{2}}V_{(j)}+\sum_{j>\frac{n+1}{2}}\left(1-V_{(j)}\right)-\frac{n+1}{4}\right)+o_P(1).\end{aligned}$$

Let $Z_1,\dots,Z_{n+1}\sim\mathrm{Exp}(1)$ independent. Denoting

$$S_t=\frac{1}{\sqrt{n+1}}\sum_{j=1}^{\lfloor(n+1)t\rfloor}(Z_j-1)$$

we get by Lemma A.7 (f), $\frac{1}{n+1}\sum_{j=1}^{n+1}Z_j\xrightarrow{P}1$ and another application of (61) (showing that the centering step cancels)

$$\frac{1}{\sqrt{n+1}}\sum_{j=1}^{n}\left(\min\left(V_j,1-V_j\right)-\frac{1}{4}\right)$$

$$
\begin{aligned}
&\overset{\mathcal{D}}{=} (1+o_P(1)) \cdot \frac{1}{n+1} \cdot \left( \sum_{j \leqslant \frac{n+1}{2}} S_{j/(n+1)} + \sum_{j > \frac{n+1}{2}} \left(S_1 - S_{j/(n+1)}\right) - \frac{n+1}{4} S_1 \right) \\
&\quad + o_P(1) \\
&= (1+o_P(1)) \cdot \left( \int_0^{\frac{1}{2}} S_t \, \mathrm{d}t + \int_{\frac{1}{2}}^1 (S_1 - S_t) \, \mathrm{d}t - \frac{1}{4} S_1 \right) + o_P(1) \\
&= (1+o_P(1)) \cdot \left( \int_0^{\frac{1}{2}} (S_t - tS_1) \, \mathrm{d}t - \int_{\frac{1}{2}}^1 (S_t - tS_1) \, \mathrm{d}t \right) + o_P(1),
\end{aligned}
$$

where we used $\frac{1}{4} = \int_0^{\frac{1}{2}} t \, \mathrm{d}t + \int_{\frac{1}{2}}^1 (1-t) \, \mathrm{d}t$ in the last line. For the sample median we obtain in a similar fashion

$$
\sqrt{n+1} \cdot \left( V_{\left(\frac{n}{2}\right)} - \frac{1}{2} \right) \overset{\mathcal{D}}{=} (1+o_P(1)) \cdot \left( S_{\frac{1}{2}} - \frac{1}{2} S_1 + o_P(1) \right).
$$

Thus, by an application of the functional central limit theorem in combination with the continuous mapping theorem the joint limit distribution is given by

$$
\left( \int_0^{\frac{1}{2}} B_t \, \mathrm{d}t - \int_{\frac{1}{2}}^1 B_t \, \mathrm{d}t \, , \, B(1/2) \right),
$$

where $\{B(\cdot)\}$ is a Brownian bridge with $\mathbb{E}(B(t)B(s)) = \mathbb{C}\mathrm{ov}(B(t), B(s)) = \min(s,t) - st$. The assertion follows by an application of Fubini's theorem as

$$
\begin{aligned}
&\mathbb{C}\mathrm{ov}\left( B\left(\frac{1}{2}\right), \int_0^{\frac{1}{2}} B(t) \, \mathrm{d}t - \int_{\frac{1}{2}}^1 B(t) \, \mathrm{d}t \right) \\
&\qquad = \int_0^{\frac{1}{2}} \mathbb{E}\left(B(1/2) \, B(t)\right) \, \mathrm{d}t - \int_{\frac{1}{2}}^1 \mathbb{E}\left(B(1/2) \, B(t)\right) \, \mathrm{d}t \\
&\qquad\qquad = \int_0^{\frac{1}{2}} \left( t - \frac{1}{2} \cdot t \right) \, \mathrm{d}t - \int_{\frac{1}{2}}^1 \left( \frac{1}{2} - \frac{1}{2} \cdot t \right) \, \mathrm{d}t = 0.
\end{aligned}
$$

□

*Proof of Theorem 2.6.* The proof of (a) for both $\mathcal{LS}$ and $\widetilde{\mathcal{LS}}$ follows immediately from Lemmas A.8 and A.9. The almost sure negativity follows from (51).

By (21) and (22) it holds

$$
\begin{aligned}
&\sqrt{\frac{3 \cdot mn}{m+n}} \cdot \left( \mathcal{LS}_T(\widehat{P}_m, \widehat{Q}_n) - \mathcal{LS}_T(\widehat{Q}_n, \widehat{P}_m) \right) \\
&\qquad\qquad = \sqrt{\frac{3 \cdot mn}{m+n}} \cdot 2 \cdot A_{1,2}(\widehat{P}_m, \widehat{Q}_n) + o_P(1).
\end{aligned}
$$

Furthermore, for $X \sim P$ independent of $Y_i$ (and $V_i$ as in (39)), by (29)

$$\begin{aligned}&\int I(x, Y_i, P)\, \mathrm{d}P(x)\\ &= \mathbb{P}(D_T(X,P) < D_T(Y_i,P)) + \frac{1}{2}\mathbb{P}(D_T(X,P) = D_T(Y_i,P))\\ &= 2 \cdot \min(V_i, 1-V_i).\end{aligned}$$

Analogously, with $Y \sim P = Q$ independent of $X_j$ (and $U_j$ as in (39))

$$\int I(X_j, y, P)\, \mathrm{d}Q(y) = 1 - 2 \cdot \min(U_j, 1-U_j).$$

This shows that the asymptotic distribution of the $\mathcal{LS}$-difference statistic is determined by sums over $\min(V_i, 1-V_i)$ resp. over $\min(U_j, 1-U_j)$. Lemmas A.8–A.9, on the other hand, show that the $\mathcal{LS}$-sum statistic is determined by the medians. Thus an application of Lemma A.10 in combination with the independence of the two samples yields the result. □

*Proof of Theorem 2.7.* By Lemma A.8 it holds

$$\widetilde{\mathcal{LS}}_T(\widehat{Q}_m, \widehat{P}_n) + \widetilde{\mathcal{LS}}_T(\widehat{P}_n, \widehat{Q}_m) - 1 = o_P(\log(m+n)/(m+n)).$$

So, it is sufficient to show

$$\widetilde{\mathcal{LS}}_T(\widehat{P}_m, \widehat{Q}_n) - \widetilde{\mathcal{LS}}_S(\widehat{P}_m, \widehat{Q}_n) = O_P\left(\frac{\log(m+n)}{m+n}\right), \tag{62}$$

from which we get, by relabelling the two samples, $\widetilde{\mathcal{LS}}_T(\widehat{Q}_n, \widehat{P}_m) - \widetilde{\mathcal{LS}}_S(\widehat{Q}_n, \widehat{P}_m) = O_P(\log(m+n)/(m+n))$ and thus the assertion. Firstly, for the here considered sample version of simplicial depth with respect to closed simplices it holds for $i = 1, \ldots, m$

$$D_S(X_{(i)}, \widehat{P}_m) = \frac{m - 1 + (m-i)\cdot(i-1)}{\binom{m}{2}}, \tag{63}$$

since there exist $(m-i)\cdot(i-1)$ intervals $[X_{(j)}, X_{(k)}]$ containing $X_{(i)}$ with $j < k$ and both $j, k \neq i$ and $m-1$ intervals that can be represented as $\operatorname{conv}\{X_{(i)}, X_{(j)}\}$, $j \neq i$. By using a similar combinatorial argument, we obtain for $x \in (X_{(i)}, X_{(i+1)})$, $i = 1, \ldots, m-1$,

$$D_S(x, \widehat{P}_m) = \frac{i \cdot (m-i)}{\binom{m}{2}}.$$

Denote $M = m/2$ and $f_m(i) = \left\lceil M - \sqrt{(M-i)^2 - (i-1)}\right\rceil$. Making use of the notation as in (39), for $Y_j$ with continuous distribution $P$ and $i \leqslant M - \sqrt{M - 3/4} + 1/2$, it holds that

$$D_S(X_{(i)}, \widehat{P}_m) \leqslant D_S(Y_j, \widehat{P}_m)$$

$$
\begin{aligned}
&\overset{\text{a.s.}}{\Longleftrightarrow} \quad X_{(k)} < Y < X_{(k+1)} \quad \text{for some } k \text{ with} \\
&\qquad\qquad\qquad\qquad\qquad\quad m-1+(m-i)(i-1) \leqslant k(m-k) \\
&\overset{\text{a.s.}}{\Longleftrightarrow} \quad X_{(f_m(i))} < Y_j < X_{(m-f_m(i)+1)} \\
&\overset{\text{a.s.}}{\Longleftrightarrow} \quad U_{(f_m(i))} < V_j < U_{(m-f_m(i)+1)}.
\end{aligned}
$$

The restriction $i \leqslant M - \sqrt{M-3/4} + 1/2$ is required to guarantee that the integrand is positive. An analogous assertion holds for $i \geqslant M+\sqrt{M-3/4}+1/2$ by the symmetry of the simplicial depth (looking forward or backward). In all other cases the depth of $Y_j$ cannot be greater or equal than that of $X_{(i)}$. Together with (40) this yields, as $f_m(i) - i \geqslant 0$,

$$\widetilde{\mathcal{LS}}_T(\widehat{P}_m, \widehat{Q}_n) - \widetilde{\mathcal{LS}}_S(\widehat{P}_m, \widehat{Q}_n) \geqslant 0.$$

Furthermore, $\mathbb{E}\left(\mathbb{1}_{\{V_j \in [U_{(a)}, U_{(b)}]\}} \middle| U_1, ..., U_m\right) = U_{(b)} - U_{(a)}$ entails in combination with A.7 (b) and (c)

$$\mathbb{E}\left(\mathbb{1}_{\{Y_j \in [X_{(a)}, X_{(b)}]\}}\right) = \frac{b-a}{m+1}.$$

For $i \geqslant 1$ it holds $(M-i) - \sqrt{(M-i)^2 - (i-1)} \leqslant \sqrt{i}$ by an application of the inverse triangular inequality for the square-root. Consequently,

$$
\begin{aligned}
&\mathbb{E}\left|\widetilde{\mathcal{LS}}_T(\widehat{P}_m, \widehat{Q}_n) - \widetilde{\mathcal{LS}}_S(\widehat{P}_m, \widehat{Q}_n)\right| \\
&= \frac{4}{m} \sum_{i \leqslant M - \sqrt{M-3/4}+1/2} \frac{f_m(i) - i}{m+1} + \frac{4}{m} \sum_{i > M - \sqrt{M-3/4}+1/2}^{M} \frac{M - i + \frac{1}{2}}{m+1} \\
&\leqslant \frac{1}{M^2} \sum_{i \leqslant M - 2\sqrt{M}} \left((M-i) - \sqrt{(M-i)^2 - (i-1)}\right) + o\left(\frac{\log(m+n)}{m+n}\right) \\
&\leqslant \frac{1}{M^2} \sum_{i \leqslant M - 2\sqrt{M}} \left((M-i) - \sqrt{(M-i)^2 - i}\right) + o\left(\frac{\log(m+n)}{m+n}\right).
\end{aligned}
$$

Now, by splitting the sum into the two summands and noting that for positive monotone functions $g$ it holds $\sum_{j=a+1}^{b} g(j) = \int_a^b g(t) \text{ d}t + O(\max(g(a), g(a+1), g(b), g(b+1)))$, we have

$$
\begin{aligned}
&\frac{1}{M^2} \sum_{i \leqslant M - 2\sqrt{M}} \left((M-i) - \sqrt{(M-i)^2 - i}\right) \\
&= \frac{1}{M^2} \int_0^{M-2\sqrt{M}} \left((M-t) - \sqrt{(M-t)^2 - t}\right) \text{ d}t + o\left(\frac{\log(m+n)}{m+n}\right) \\
&= \int_{2/\sqrt{M}}^{1} \left(s - \sqrt{s^2 - \frac{1-s}{M}}\right) \text{ d}s + o\left(\frac{\log(m+n)}{m+n}\right) = O\left(\frac{\log(m+n)}{m+n}\right).
\end{aligned}
$$

The last line follows by a Taylor expansion with some $s^2 - \frac{1-s}{M} \leqslant \xi_s \leqslant s^2$. Thus, fulfilling $\xi_s \geqslant 3s^2/4$, it holds

$$\begin{aligned}&\int_{2/\sqrt{M}}^{1}\left(\sqrt{s^2}-\sqrt{s^2-\frac{1-s}{M}}\right)\,\mathrm{d}s\\&=\frac{1}{2M}\int_{2/\sqrt{M}}^{1}\frac{1-s}{s}\,\mathrm{d}s+O\left(\frac{1}{M^2}\right)\int_{2/\sqrt{M}}^{1}\frac{1}{\xi_s^{3/2}}\,\mathrm{d}s\\&=O\left(\frac{\log(M)}{M}\right)+O\left(\frac{1}{M^2}\right)\int_{2/\sqrt{M}}^{1}\frac{1}{(3s^2/4)^{3/2}}\,\mathrm{d}s=O\left(\frac{\log(m+n)}{m+n}\right).\end{aligned}$$

This proves (62) and implies in combination with Lemma A.8 that

$$\widetilde{\mathcal{LS}}_S(\widehat{P}_m,\widehat{Q}_n)+\widetilde{\mathcal{LS}}_S(\widehat{Q}_n,\widehat{P}_m)-1=O_P\left(\frac{\log(m+n)}{m+n}\right).$$

Next, we show that this rate also holds true when $\widetilde{\mathcal{LS}}_S$ is replaced by $\mathcal{LS}_S$ in the above expression: First note, that for $Y_j \notin [X_{(1)}, X_{(m)}]$ it holds

$$\sum_{i=1}^{m}\mathbb{1}_{\left\{D_S(X_i,\widehat{P}_m)=D_S(Y_j,\widehat{P}_m)\right\}}=0=\sum_{i=1}^{m}\mathbb{1}_{\left\{D_T(X_i,\widehat{P}_m)=D_T(Y_j,\widehat{P}_m)\right\}}.$$

Now consider $Y_j \in [X_{(1)}, X_{(m)}]$ in addition to $Y_j \neq X_{(i)}$, $i = 1, \dots, m$, where the latter occurs with probability 1. Then, the above sum with respect to the Tukey depth equals 2. The same sum with respect to the simplicial depths equals at most 2 because the formula in (63) is a parabola in $i$ and, as such, can take each value at most twice. Consequently,

$$\sum_{i=1}^{m}\mathbb{1}_{\left\{D_S(X_i,\widehat{P}_m)=D_S(Y_j,\widehat{P}_m)\right\}}\leqslant\sum_{i=1}^{m}\mathbb{1}_{\left\{D_T(X_i,\widehat{P}_m)=D_T(Y_j,\widehat{P}_m)\right\}}\quad\text{a.s.}$$

Thus, (59) and (60) imply that

$$\widetilde{\mathcal{LS}}_S(\widehat{P}_m,\widehat{Q}_n)-\mathcal{LS}_S(\widehat{P}_m,\widehat{Q}_n)=O_P\left(\frac{m+n}{mn}\right)\tag{64}$$

and analogous assertions for $\widetilde{\mathcal{LS}}_S(\widehat{Q}_n,\widehat{P}_m)-\mathcal{LS}_S(\widehat{Q}_n,\widehat{P}_m)$ as well as the $\mathcal{LS}$-statistics with respect to the univariate Tukey depth. By Proposition 3.3 the Tukey depth fulfills the assumptions of Theorem 2.3. So, (10) also holds for $\mathcal{LS}_S$ and $\widetilde{\mathcal{LS}}_S$ by applications of (62) and (64). □

## Appendix B: Proofs of Section 3

*Proof of Theorem 3.1.* $\mathrm{Asm}_3$ is trivial because the $h$-depth can be regarded as a U-statistic with bounded kernel and $\mathrm{Asm}_2$ follows by an inequality in Wynne and Nagy [88, Proof of Theorem 5] in combination with the layer cake representation of moments and Hölder's inequality. □

*Proof of Theorem 3.2.* Following the reasoning in the proofs of Nagy et al. [60, Theorems 4.15 and 5.2 part II], we have that

$$\sup_{x \in C^d[0,1]} |ID(x, P_1) - ID(x, P_2)| \leqslant \int_0^1 \sup_{u \in \mathbb{R}^d} |D_d(u, P_1(t)) - D_d(u, P_2(t))| \,\mathrm{d}t \tag{65}$$

for any two probability measures $P_1, P_2 \in \mathcal{P}_{C^d[0,1]}$. $\mathrm{Asm}_3$ follows directly from this and $D_d$ satisfying this assumption.

For the proof of $\mathrm{Asm}_2$ we make use of (65), with distributions $\widehat{P}_m$ and $P$, and Jensen's inequality with $\beta > 1/2$, obtaining

$$\sup_{x \in C^d[0,1]} |ID(x, \widehat{P}_m) - ID(x, P)|^{2\beta} \leqslant \int_0^1 \sup_{u \in \mathbb{R}^d} |D_d(u, \widehat{P}_m(t)) - D_d(u, P(t))|^{2\beta} \,\mathrm{d}t, \tag{66}$$

which results in

$$\begin{aligned} E_m &:= \mathbb{E}\left( \sup_{x \in C^d[0,1]} |ID(x, \widehat{P}_m) - ID(x, P)|^{2\beta} \right) \\ &\leqslant \mathbb{E}\left( \int_0^1 \sup_{u \in \mathbb{R}^d} |D_d(u, \widehat{P}_m(t)) - D_d(u, P(t))|^{2\beta} \,\mathrm{d}t \right) < \infty \end{aligned}$$

by taking into account that the multivariate depths, considered here, take values in $[0, 1]$ [94]. Note, that Nagy et al. [60, Theorem 3.1] entails the measurability of the integrands. Then, by making use of the Fubini-Tonelli Theorem and (66), we obtain $E_m \leqslant \int_0^1 \mathbb{E}(\sup_{u \in \mathbb{R}^d} |D_d(u, \widehat{P}_m(t)) - D_d(u, P(t))|^{2\beta}) \,\mathrm{d}t$. Applying to this that $\mathrm{Asm}_2$ holds uniformly for $D_d$ with respect to the probability measure, we have that there exists a constant $C$ such that $E_m \leqslant \int_0^1 (C/m^\beta) \,\mathrm{d}t = C/m^\beta$, and, consequently $E_m = O(m^{-\beta})$. □

In what follows, for any $x, u \in \mathbb{R}^d$, $d \in \mathbb{N}$, $H_u(x) := \{w \in \mathbb{R}^d : \langle u, w \rangle \leqslant \langle u, x \rangle\}$ denotes a halfspace determined by $u$ and $x$.

**Lemma B.1.** *Let $P_1, P_2$ be probability measures on $\mathbb{R}^d$, $d \in \mathbb{N}$. Then, for any $x \in \mathbb{R}^d$,*

$$\begin{aligned} |D_T(x, P_1) - D_T(x, P_2)| &\leqslant \sup_{\|u\|_2 = 1} |P_1(H_u(x)) - P_2(H_u(x))| \\ &\leqslant \sup_{H \in \mathcal{H}} |P_1(H) - P_2(H)|. \end{aligned}$$

*Proof.* It follows directly from the definitions of the Tukey depth and of halfspaces. □

*Proof of Proposition 3.3.* Proof for $\mathrm{Asm}_1$. For any continuous $P$ on $\mathbb{R}$, we have that $D_T(x, P) = \min(P((-\infty, x]), 1 - P((-\infty, x]))$ and $P((-\infty, X]) \sim U(0, 1)$, for any $X \sim P$, by the probability integral transform. Thus, $D_T(X,P) \sim U(0,0.5)$, which satisfies $\text{Asm}_1$.

Proof for $\text{Asm}_2$. It follows from Lemma B.1 by combining the layer cake representation of moments and the univariate Dvoretzky-Kiefer-Wolfowitz inequality [57, Corollary 1].

Proof for $\text{Asm}_3$. We refer to the proof of Proposition 3.5, which contains the univariate case.

Proof for $\text{Asm}_4$. Note that

$$\begin{aligned}&\mathbb{P}\left(|D_T(X_1,\widehat{P}_m^{-\{1\}})-D_T(X,\widehat{P}_m^{-\{1\}})|\leqslant \frac{2C}{m}\middle| X_2,\dots,X_m\right)\\ &\quad\overset{\mathcal{D}}{=}\mathbb{P}\left(|D_T(X_m,\widehat{P}_m^{-\{m\}})-D_T(X,\widehat{P}_m^{-\{m\}})|\leqslant \frac{2C}{m}\middle| X_1,\dots,X_{m-1}\right)\end{aligned}\tag{67}$$

as $X_1,\dots,X_m$ are iid. Due to conditional independence and because of $\widehat{P}_m^{-\{m\}}=\widehat{P}_{m-1}$, we obtain

$$\begin{aligned}&\mathbb{P}\left(|D_T(X_m,\widehat{P}_m^{-\{m\}})-D_T(X,\widehat{P}_m^{-\{m\}})|\leqslant \frac{2C}{m}\middle| X_1,\dots,X_{m-1}\right)\\
&\overset{\text{a.s.}}{\leqslant}\mathbb{P}\left(|D_T(X_m,\widehat{P}_m^{-\{m\}})-D_T(X,\widehat{P}_m^{-\{m\}})|\leqslant \frac{2C}{m-1}\middle| X_1,\dots,X_{m-1}\right)\\
&\overset{\text{a.s.}}{=}\sum_{i=0}^{\lfloor\frac{m-1}{2}\rfloor}\mathbb{P}\Bigg(D_T(X_m,\widehat{P}_{m-1})=\frac{i}{m-1},\\
&\qquad\qquad D_T(X,\widehat{P}_{m-1})\in\left[\frac{i-2C}{m-1},\frac{i+2C}{m-1}\right]\Bigg| X_1,\dots,X_{m-1}\Bigg)\\
&=\sum_{i=0}^{\lfloor\frac{m-1}{2}\rfloor}\mathbb{P}\left(D_T(X_m,\widehat{P}_{m-1})=\frac{i}{m-1}\middle| X_1,\dots,X_{m-1}\right)\\
&\qquad\qquad\cdot\mathbb{P}\left(D_T(X,\widehat{P}_{m-1})\in\left[\frac{i-2C}{m-1},\frac{i+2C}{m-1}\right]\middle| X_1,\dots,X_{m-1}\right)\\
&=\sum_{i=0}^{\lfloor\frac{m-1}{2}\rfloor}\sum_{j=-2C}^{2C}\mathbb{P}\left(D_T(X_m,\widehat{P}_{m-1})=\frac{i}{m-1}\middle| X_1,\dots,X_{m-1}\right)\\
&\qquad\qquad\cdot\mathbb{P}\left(D_T(X,\widehat{P}_{m-1})=\frac{i+j}{m-1}\middle| X_1,\dots,X_{m-1}\right).\end{aligned}\tag{68}$$

Let $F$ be the cumulative distribution function associated to $P$, following Zuo and He [93, Example 2], we denote $\mathfrak{D}_i := F(X_{(i+1)})-F(X_{(i)}), i\in\{0,\dots,m\}$, with $F(X_{(0)})=0$ and $F(X_{(m+1)})=1$, by convention, which results in

$$\mathbb{P}\left(D_T(X,\widehat{P}_m)=\frac{i}{m}\middle| X_1,\dots,X_m\right)\overset{\text{a.s.}}{=}\mathfrak{D}_i+\mathfrak{D}_{m-i}\text{ for } i\in\{0,\dots,\lceil m/2\rceil-1\}\tag{69}$$

and

$$\mathbb{P}\left(D_T(X,\widehat{P}_m)=\frac{i}{m}\,\middle|\,X_1,...,X_m\right)\overset{\text{a.s.}}{=}\mathbb{1}_{\{\mathrm{mod}(m,2)=0\}}\cdot\mathfrak{D}_i \text{ for } i=\lceil m/2\rceil\,. \quad (70)$$

As $X_1,...,X_m$ are continuously distributed, Lemma A.7 (a)–(b) entails that $\mathfrak{D}_i$ follows a Beta$(1,m)$-distribution for $i\in\{0,...,m\}$. Then, by means of Hölder's inequality, we have that

$$\mathbb{E}\left(\prod_{j=1}^{4}\mathfrak{D}_{i_j}^{s_j}\right)\leqslant\mathbb{E}\left(\mathfrak{D}_{i_1}^4\right)=\frac{24}{(m+1)(m+2)(m+3)(m+4)},$$

for any combination of $s_1,...,s_4$ such that $\sum_{j=1}^4 s_j=4$ with $s_1\in\{1,...,4\}$, $s_2\in\{0,...,2\}$ and $s_3,s_4\in\{0,1\}$. Hence all of the mixed moments we have to consider are uniformly $O(m^{-4})$. In particular, replacing $m$ by $m-1$ in the above expressions has no impact on the rate so that we can use the representation (69)-(70) for both factors in (68). Consequently, we get with (67)

$$\mathbb{E}\left(\left[\mathbb{P}\left(|D_T(X_1,\widehat{P}_m^{-\{1\}})-D_T(X,\widehat{P}_m^{-\{1\}})|\leqslant\frac{2C}{m}\,\middle|\,X_2,...,X_m\right)\right]^2\right)$$
$$=O(m^2)\cdot O(m^{-4})=O(m^{-2}).$$

□

*Proof of Theorem 3.4.* Let us prove the statement for $\mathrm{Asm}_3$. Applying the definition of U-statistics with bounded kernel $K_x(\cdot)$ of order $k>1$,

$$\begin{aligned}&|D(x,\widehat{P}_m)-D(x,\widehat{P}_m^{-\{m\}})|\\&=\left|\left(\binom{m}{k}^{-1}-\binom{m-1}{k}^{-1}\right)\cdot\sum_{1\leqslant i_1<...<i_k\leqslant m-1}K_x(X_{i_1},...,X_{i_k})\right.\\&\quad\left.+\binom{m}{k}^{-1}\cdot O\left(\binom{m}{k}-\binom{m-1}{k}\right)\right|.\end{aligned}$$

Taking into account that $\sum_{1\leqslant i_1<...<i_k\leqslant m-1}K_x(X_{i_1},...,X_{i_k})=O\big(\binom{m-1}{k}\big)$ and that $\big(\binom{m}{k}^{-1}-\binom{m-1}{k}^{-1}\big)=\binom{m-1}{k}^{-1}\cdot O(m^{-1})$, the first summand is $O(m^{-1})$. After applying some algebra on the second summand, we also obtain it is $O(m^{-1})$. Thus, the triangle inequality and the exchangeability-argument give us

$$\sup_{x\in\mathbb{R}^d}|D(x,\widehat{P}_m)-D(x,\widehat{P}_m^{-\{1\}})|=O(m^{-1}).$$

Let us prove the statement for $\mathrm{Asm}_2$. With the kernel $K_x(\cdot)=\mathbb{1}_{\{x\in\mathfrak{L}(\cdot)\}}$ of order $k>1$ it holds $D(x,\widehat{P}_m)=\sum_{1\leqslant i_1<...<i_k\leqslant m}\mathbb{1}_{\{x\in\mathfrak{L}(X_{i_1},...,X_{i_k})\}}$. Consider the corresponding V-statistic $D_{\text{mod}}(x, \widehat{P}_m) = \int \mathbb{1}_{\{x \in \mathfrak{L}(y_1,...,y_k)\}} \, \mathrm{d}\widehat{P}_m^k(y_1, ..., y_k)$. Similar to Dümbgen [27, Section 1],

$$\begin{aligned} &D_{\text{mod}}(x, \widehat{P}_m) - D(x, P) \\ &\qquad = \int (\widehat{P}_m - P) \mathrm{A}(x, y_1, ..., y_{k-1}) \, \mathrm{d} \left( \sum_{j=0}^{k-1} \widehat{P}_m^j P^{d-j}(y_1, ..., y_{k-1}) \right), \end{aligned}$$

where $\mathrm{A}(x, y_1, ..., y_{k-1}) := \{y \in \mathbb{R}^d : x \in \mathfrak{L}(y_1, ..., y_{k-1}, y)\}$. Then,

$$|D_{\text{mod}}(x, \widehat{P}_m) - D(x, P)| \leqslant k \cdot \sup_{I \in \mathcal{I}} |\widehat{P}_m(I) - P(I)|,$$

for the class of sets $\mathcal{I} := \{\mathrm{A}(x, y_1, ..., y_{k-1}) : x, y_1, ..., y_{k-1} \in \mathbb{R}^d\}$. Consequently,

$$\begin{aligned} &\mathbb{P}\left( \sqrt{m} \cdot \sup_{x \in \mathbb{R}^d} |D_{\text{mod}}(x, \widehat{P}_m) - D_(x, P)| > M \right) \\ &\qquad\qquad \leqslant \mathbb{P}\left( \sqrt{m} \cdot \sup_{I \in \mathcal{I}} |\widehat{P}_m(I) - P(I)| > \frac{M}{k} \right). \end{aligned}$$

Moreover, the class of sets $\mathcal{I}$ has finite VC-dimension $V(\mathcal{I})$ (respectively belongs to a larger class of sets $\mathcal{J}$ with finite VC-dimension $V(\mathcal{J})$). Then, as we are under assumption (M), by the Dvoretzky-Kiefer-Wolfowitz inequality [2, Theorem 3.1], we have that for any fixed $\upsilon \geqslant V(\mathcal{J})$ and $\epsilon \in (0, 1]$, there exist constants $m_0(\upsilon, \epsilon)$, $M_0(\upsilon, \epsilon)$ and $\mathcal{K}(\upsilon)$ such that

$$\mathbb{P}\left( \sqrt{m} \sup_{J \in \mathcal{J}} |\widehat{P}_m(J) - P(J)| > \frac{M}{k} \right) \leqslant \mathbb{1}_{\{\frac{M}{k} < M_0\}} + \mathbb{1}_{\{\frac{M}{k} \geqslant M_0\}} \mathcal{K}(\upsilon) \cdot e^{-(2-2\epsilon)(\frac{M}{k})^2}$$

for all $m \geqslant m_0(\upsilon, \epsilon)$ and $M \geqslant M_0(\upsilon, \epsilon)$. Then, by means of the layer cake representation we obtain

$$\mathbb{E}\left( \sup_{x \in \mathbb{R}^d} |D_{\text{mod}}(x, \widehat{P}_m) - D(x, P)|^2 \right) = O(m^{-1}). \tag{71}$$

On the other hand,

$$\sup_{x \in \mathbb{R}^d} |D(x, \widehat{P}_m) - D_{\text{mod}}(x, \widehat{P}_m)|^2 \leqslant \frac{C_{det}}{m^2} \tag{72}$$

as, by Serfling [73, Section 5.7.3], the kernels of both the U-statistic $D(x, \widehat{P}_m)$ and the V-statistic $D_{\text{mod}}(x, \widehat{P}_m)$ are identical and bounded. Then, making use of (71) and (72), we obtain

$$\begin{aligned} &\mathbb{E}\left( \sup_{x \in \mathbb{R}^d} |D(x, \widehat{P}_m) - D(x, P)|^2 \right) \\ &\leqslant 2 \left[ \mathbb{E}\left( \sup_{x \in \mathbb{R}^d} |D(x, \widehat{P}_m) - D_{\text{mod}}(x, \widehat{P}_m)|^2 \right) \right. \end{aligned}$$

$$+\mathbb{E}\left(\sup_{x\in\mathbb{R}^d}|D_{\mathrm{mod}}(x,\widehat{P}_m)-D(x,P)|^2\right)\Bigg]$$
$$=O(m^{-1}).$$

Now, applying Hölder's inequality delivers for $\beta\in(1/2,1]$

$$\mathbb{E}\left(\sup_{x\in\mathbb{R}^d}|D(x,\widehat{P}_m)-D(x,P)|^{2\beta}\right)=O\left(m^{-\beta}\right).$$

□

*Proof of Proposition 3.5.*

Tukey depth. Let us first prove the statement for $\mathrm{Asm}_2$. From Wainwright [86, Propositions 4.19 and 4.20], we have that $\mathcal{H}$ is a system of sets with finite VC-dimension $V(\mathcal{H})$. Thus, as we are under assumption (M), the multivariate Dvoretzky-Kiefer-Wolfowitz inequaltiy holds. Denoting

$$\mathcal{M}_m:=\sqrt{m}\sup_{H\in\mathcal{H}}|\widehat{P}_m(H)-P(H)|,$$

this inequality states that for any fixed $\upsilon\geqslant V(\mathcal{H})$ and $\epsilon\in(0,1]$, there exist constants $m_0(\upsilon,\epsilon), M_0(\upsilon,\epsilon)$ and $\mathcal{K}(\upsilon)$ such that

$$\mathbb{P}(\mathcal{M}_m>M)\leqslant\mathcal{K}(\upsilon)e^{-(2-2\epsilon)M^2}$$

for each $m\geqslant m_0(\upsilon,\epsilon)$ and $M\geqslant M_0(\upsilon,\epsilon)$. The layer cake representation of the expected value delivers for $m$ large enough ($m\geqslant m_0(\upsilon,\epsilon)$)

$$\mathbb{E}\left(\mathcal{M}_m^2\right)=\int_0^\infty 2\cdot M\cdot\mathbb{P}\left(\mathcal{M}_m>M\right)\,\mathrm{d}M\leqslant 2\cdot\int_0^\infty M\cdot\mathcal{K}(\upsilon)e^{-(2-2\epsilon)M^2}\,\mathrm{d}M<\infty.$$

Consequently, $\mathbb{E}\left(\mathcal{M}_m^2\right)=O\left(1\right)$ and so,

$$\mathbb{E}\left(\sup_{H\in\mathcal{H}}|\widehat{P}_m(H)-P(H)|^2\right)=O\left(m^{-1}\right).\tag{73}$$

By Lemma B.1, we have that

$$\mathbb{E}\left(\sup_{x\in\mathbb{R}^d}|D_T(x,\widehat{P}_m)-D_T(x,P)|^2\right)\leqslant\mathbb{E}\left(\sup_{H\in\mathcal{H}}|\widehat{P}_m(H)-P(H)|^2\right)$$

which together with (73) implies that

$$\mathbb{E}\left(\sup_{x\in\mathbb{R}^d}|D_T(x,\widehat{P}_m)-D_T(x,P)|^2\right)=O\left(m^{-1}\right).$$

For $\beta\in(1/2,1]$ Hölder's inequality finally gives $\mathrm{Asm}_2$, i.e.

$$\mathbb{E}\left(\sup_{x\in\mathbb{R}^d}|D_T(x,\widehat{P}_m)-D_T(x,P)|^{2\beta}\right)=O\left(m^{-\beta}\right).$$

Let us prove the statement for Asm$_3$. Making use of the empirical distribution, we have that

$$|\widehat{P}_m(H_u(x)) - \widehat{P}_m^{-\{1\}}(H_u(x))| = \left| -\frac{1}{m(m-1)} \sum_{i=2}^{m} \mathbb{1}_{\{X_i \in H_u(x)\}} + \frac{1}{m} \mathbb{1}_{\{X_1 \in H_u(x)\}} \right|,$$

which results in $|\widehat{P}_m(H_u(x)) - \widehat{P}_m^{-\{1\}}(H_u(x))| = O(m^{-1})$. By Lemma B.1, this leads to Asm$_3$, i.e.

$$\sup_{x \in \mathbb{R}^d} |D_T(x, \widehat{P}_m) - D_T(x, \widehat{P}_m^{-\{1\}})| = O\left(m^{-1}\right).$$

*For Asm$_3$ we just observe that the following empirical depth functions can be represented as U-statistics. For Asm$_2$ we prove for the particular depths that the class of sets* $\mathrm{A}(x, y_1, ..., y_{n-1})$ *has finite VC-dimension.*

Simplicial depth. The U-statistic-estimator for the simplicial depth can be defined by using the kernel $K_x(X_{i_1}, ..., X_{i_{d+1}}) = \mathbb{1}_{\{x \in \mathrm{conv}\{X_{i_1}, ..., X_{i_{d+1}}\}\}}$. Consequently, $\mathrm{A}(x, y_1, ..., y_d) := \{y \in \mathbb{R}^d : x \in \mathrm{conv}\{y_1, ..., y_d, y\}\}$. According to Dümbgen [27, Section 1], $\mathcal{I} := \{\mathrm{A}(x, y_1, ..., y_d) : x, y_1, ..., y_d \in \mathbb{R}^d\}$ belongs to a class of sets $\mathcal{J}$ with finite VC-index $V(\mathcal{J})$. Thus, Asm$_2$ and Asm$_3$ follow immediately by Theorem 3.4.

Spherical depth. The U-statistic-estimator for the spherical depth can be defined by using the kernel $K_x(X_{i_1}, X_{i_2}) = \mathbb{1}_{\{(X_{i_1}-x)^T (X_{i_2}-x) \leqslant 0\}}$. Consequently, $\mathrm{A}(x, y_1) = \{y \in \mathbb{R}^d : (y-x)^T (y_1 - x) \leqslant 0\}$. All of these sets $\mathrm{A}(x, y_1)$ are contained in a subgraph class with finite VC-dimension according to Wainwright [86, Proposition 4.20]. Hence, Asm$_2$ and Asm$_3$ directly arise from Theorem 3.4 for the spherical depth.

Lens depth. The U-statistic-estimator for the lens depth with respect to the Euclidean norm can be defined by using the kernel

$$K_x(X_{i_1}, X_{i_2}) = \mathbb{1}_{\{\|X_{i_1} - X_{i_2}\| \geqslant \max(\|x - X_{i_1}\|, \|x - X_{i_2}\|)\}}.$$

Consequently,

$$\begin{aligned} \mathrm{A}(x, y_1) &= \{y \in \mathbb{R}^d : \|y - y_1\| \geqslant \max(\|x-y\|, \|x - y_1\|)\} \\ &= \{y \in \mathbb{R}^d : \|y - y_1\| \geqslant \|x - y\|\} \cap \{y \in \mathbb{R}^d : \|y - y_1\| \geqslant \|x - y_1\|\}. \end{aligned}$$

All of these sets $\mathrm{A}(x, y_1)$ can be generated by intersecting sets from two subgraph classes with finite VC-dimension according to [86, Proposition 4.19 and 4.20]. Thus, Theorem 3.4 entails Asm$_2$ and Asm$_3$ for the lens depth.

Band depth. The U-statistic-estimator for the band depth based on $k$ elements can be defined by using the kernel

$$K_x(X_{i_1}, ..., X_{i_k}) = \prod_{j=1}^{d} \mathbb{1}_{\{\min_{\ell\in\{1,...,k\}} X_{i_\ell}^{(j)} \leqslant x^{(j)} \leqslant \max_{\ell\in\{1,...,k\}} X_{i_\ell}^{(j)}\}}.$$

Consequently,

$$\begin{aligned}
&\mathrm{A}(x, y_1, ..., y_{k-1}) \\
&= \bigcap_{j=1}^{d} \{y_0 \in \mathbb{R}^d : \min_{\ell\in\{0,...,k-1\}} y_\ell^{(j)} \leqslant x^{(j)} \leqslant \max_{\ell\in\{0,...,k-1\}} y_\ell^{(j)}\} \\
&= \bigcap_{j=1}^{d} \Big( \{y_0 \in \mathbb{R}^d : \min_{\ell\in\{0,...,k-1\}} y_\ell^{(j)} \leqslant x^{(j)}\} \\
&\qquad\qquad \cap \{y_0 \in \mathbb{R}^d : x^{(j)} \leqslant \max_{\ell\in\{0,...,k-1\}} y_\ell^{(j)}\} \Big) \\
&= \bigcap_{j=1}^{d} \left( \left( \cup_{\ell=0}^{k-1} \{y_0 \in \mathbb{R}^d : y_\ell^{(j)} \leqslant x^{(j)}\} \right) \cap \left( \cup_{\ell=0}^{k-1} \{y_0 \in \mathbb{R}^d : x^{(j)} \leqslant y_\ell^{(j)}\} \right) \right).
\end{aligned}$$

Each set $\{y_0 \in \mathbb{R}^d : y_\ell^{(j)} \leqslant x^{(j)}\}$, respectively $\{y_0 \in \mathbb{R}^d : x^{(j)} \leqslant y_\ell^{(j)}\}$, in the expression on the right hand side of the above equality belongs to a certain subgraph class with finite VC-dimension according to Wainwright [86, Proposition 4.19 and 4.20]. Since taking the union respectively the intersection of a finite amount of set classes always generates a set class with finite VC-dimension, we can conclude by Theorem 3.4 that $\mathrm{Asm}_2$ holds for the band depth ($\mathrm{Asm}_3$ holds trivially). □

## Appendix C: Simulation results

In Section 4 of the main paper we summarised the findings from the simulation study. In this section, we report the corresponding full simulation study.

Section C.1 is dedicated to compare the performance of different tests based on the $\mathcal{LS}$-tuple on functional and multivariate data, to confirm the intuition obtained in Section 2.4 when analysing univariate data. Real world functional data typically consist on non-smoothed observations given on a grid. Historically, this raw data are pre-processed into smooth curves [69]. Although preprocessing (including such based on statistical data depths) is very common nowadays, we do know that it influences the data analysis outcome but not necessarily in a good way; see e.g. [10, 75]. To mimic the raw nature of real functional data, we compare the different $\mathcal{LS}$-based tests making use of Brownian motions. In addition, we make use of the integrated Tukey depth that has well-known properties.

Due to the results of Section C.1, the properties of the Joint-TP test for different functional depths are compared in Section C.2. There, in addition to the Brownian motion, we use smooth, simulated, functional data. This is of particular importance as the user may wish to apply the proposed test to smooth functional data. Additionally, the behaviour of some functional depths under raw data is not optimal as it was defined with smooth data in mind.

Then, Section C.5 compares the proposed test with existing functional two-sample tests in the literature.

The simulations in Sections C.1–C.4 were run in MATLAB R2019b (or newer) with our own implementation of the different depth functions used. For the the simulations in Section C.5, R (version 4.2.2) was used, as most implementations of other functional two-sample tests are available in this programming language. Our simulation results are based on 1000 trials and the level of the tests is given by $\alpha = 0.05$. All simulation results are based on 1000 trials and the level of the tests is given by $\alpha = 0.05$. The functional observations are evaluated at 1001 equidistant discretisation points on the interval $[0, 1]$. For generating the plots we made use of [9] and [45].

### C.1. Tests based on the joint $\mathcal{LS}$-tuple for the functional integrated Tukey depth and the multivariate Tukey depth

In this section we show that the findings exibited in Sections 2.2 and 2.4 carry over to functional data. For that, we make use of the $\mathcal{LS}$-tests based on the functional integrated Tukey depth. We concentrate on an alternative falling into one of the blind spots with a huge asymmetry in the two projection statistics. Thus, under $H_0$, both samples consist of standard Brownian motions $(B_t)_{t\in[0,1]}$ (which we will regard as Model 1) while under the alternative we scale and shift the data in one sample (the one drawn from distribution $Q$), considering samples from $(0.8 \cdot B_t + 0.15)_{t\in[0,1]}$. We take $m = n = 100$ observations in each sample. For an illustration, the left panel in Figure 5 displays in grey colour a sample of 10 functional elements drawn from $(B_t)_{t\in[0,1]}$. There, we have also plotted 10 dark blue functional elements drawn from $(0.8 \cdot B_t + 0.15)_{t\in[0,1]}$.

The simulation results are in Table 3, which displays the empirical size (rejection rate under the null hypothesis), power (rejection rate under the alternative hypothesis) and size-corrected power. The size corrected power is computed by replacing the theoretically obtained critical value from the asymptotic result by an empirical quantile of simulated observations under the null hypothesis. Further simulations are displayed in Tables 4 and 5. Table 4 contains the results of contrasting $(B_t)_{t\in[0,1]}$ against the alternative $(B_t + a)_{t\in[0,1]}$ for $a \in \{0.1 + 0.025i : i \in \mathbb{N} \cup \{0\} \text{ with } i \leqslant 10\}$ and Table 5 against the alternative $(b \cdot B_t + 0.15)_{t\in[0,1]}$ for $b \in \{0.85, 0.9, 0.95, 0.99\}$. We also include Tables 6 and 7 with results for multivariate data, which further confirm our observations; with the results in Table 6 corresponding to those in Table 4 and the ones in Table 7 to those in Table 5.

As in Figure 2 for the univariate case, we observe from the first row of Table 3 that in the functional case the two $\mathcal{LS}$-projection tests (first and second columns) and even more so the $\mathcal{LS}$-maximum test (fourth column) do not hold the significance level, with empirical sizes that are far too liberal. Indeed, they are,

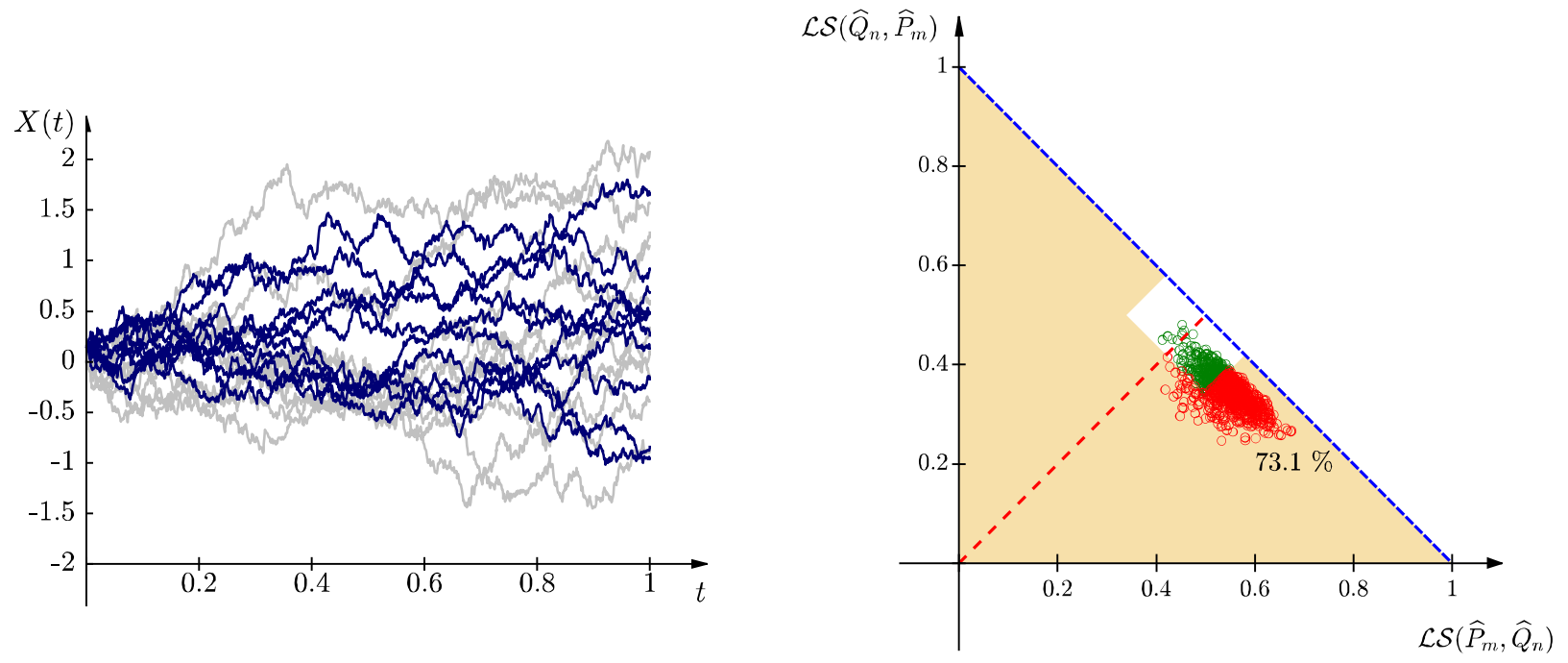

Fig 5. *Left panel: 10 realizations of* $(B_t)_{t\in[0,1]}$ *in grey and 10 of* $(0.8\cdot B_t+0.15)_{t\in[0,1]}$ *in dark blue. In the horizontal axis t is displayed. In the vertical axis* $X(t)$ *is displayed, representing either a realization of* $B_t$ *or a realization* $0.8\cdot B_t+0.15$*. Each of the 20 realizations are drawn independently. Right panel: The circles represent 1000 independent realizations of the* $\mathcal{LS}$*-tuple* $(\mathcal{LS}(\widehat{P}_m,\widehat{Q}_n)$ *in the horizontal axis and* $\mathcal{LS}(\widehat{Q}_n,\widehat{P}_m)$ *in the vertical axis) under* $(B_t)_{t\in[0,1]}$ *(distribution P) against* $(0.8\cdot B_t+0.15)_{t\in[0,1]}$ *(distribution Q). The functional integrated Tukey depth is used and* $m=n=100$ *observations in each sample. The red circles (unlike the green ones) result in rejections made by the Joint-TP statistic and the shaded orange region is its rejection region. For both panels: the functional observations are evaluated at 1001 equidistant discretisation points.*

Table 3

*Size (first row) and power comparison (second row) under* $(B_t)_{t\in[0,1]}$ *(distribution P) against* $(0.8\cdot B_t+0.15)_{t\in[0,1]}$ *(distribution Q) of the statistics* $\mathcal{LS}(\widehat{P}_m,\widehat{Q}_n)$ *(first column),* $\mathcal{LS}(\widehat{Q}_n,\widehat{P}_m)$ *(second column),* $\mathcal{LS}$*-difference (third column),* $\mathcal{LS}$*-maximum (fourth column), Joint-CC (fifth column) and Joint-TP (sixth column). The third column displays the size corrected power for the statistics with size above the significance level. The statistics are based on the functional integrated Tukey depth and* $m=n=100$ *observations in each sample.*

| | $\mathcal{LS}(\widehat{P}_m,\widehat{Q}_n)$ | $\mathcal{LS}(\widehat{Q}_n,\widehat{P}_m)$ | $\mathcal{LS}$-difference | $\mathcal{LS}$-maximum | Joint-CC | Joint-TP |
|---|---|---|---|---|---|---|
| Size | 7.4% | 8.2% | 4.5% | 14.9% | 4.6% | 4.6% |
| Power | 20.3% | 97.1% | 69.7% | 97.2% | 73.6% | 73.1% |
| Corrected Power | 16.4% | 94.2% | – | 91.6% | – | – |

respectively, 2.4, 3.2 and 9.9 points above the 5% level under the null hypothesis. In contraposition, the $\mathcal{LS}$-difference test (third column), the Joint-CC (fifth column) and Joint-TP test (sixth column) exhibit very good size behaviour.

Moreover, the second row of Table 3 shows that the $\mathcal{LS}$-difference test (third column), the Joint-CC (fifth column) and Joint-TP test (sixth column) exhibit an acceptable power behaviour, under the alternative. The power is not as good as for the $\mathcal{LS}$-maximum test (fourth column) and one of the two projection tests ($\mathcal{LS}(\widehat{Q}_n,\widehat{P}_m)$, second column), even after size-correcting them. This is due to the specific location of the $\mathcal{LS}$-tuple under this alternative, which we have depicted in the right panel of Figure 5. In fact, this is also what makes the power between the two projection tests (first and second columns) differ greatly in the

Table 4

*Rejection rates of $(B_t)_{t\in[0,1]}$ (distribution $P$) against $(B_t+a)_{t\in[0,1]}$ (distribution $Q$) for the different $a \in \{0.1+0.025\cdot i : i\in\mathbb{N}_0$ with $i\leqslant 10\}\cup\{0\}$ and the statistics $\mathcal{LS}(\widehat{P}_m,\widehat{Q}_n)$ (rows (a)), $\mathcal{LS}(\widehat{Q}_n,\widehat{P}_m)$ (rows (b)), $\mathcal{LS}$-difference (rows (c)), $\mathcal{LS}$-maximum (rows (d)), Joint-CC (rows (e)) and Joint-TP (rows (f)). In rows (a), (b) and (d) the size corrected power for the corresponding statistic is also included (the respective rejection rate under the null hypothesis is that given in the first column, under $a=0$). The statistics are based on the functional integrated Tukey depth and $m=n=100$ observations in each sample.*

| | ***a*** | | | | | |
|---|---|---|---|---|---|---|
| | 0 | 0.1 | 0.125 | 0.15 | 0.175 | 0.2 |
| (a) $\mathcal{LS}(\widehat{P}_m,\widehat{Q}_n)$ | 7.4% | 13.5% | 19.0% | 25.1% | 28.0% | 35.0% |
| (a) Corrected | 5.0% | 9.9% | 15.3% | 20.3% | 23.5% | 29.8% |
| (b) $\mathcal{LS}(\widehat{Q}_n,\widehat{P}_m)$ | 8.2% | 18.3% | 20.4% | 23.4% | 32.6% | 40.0% |
| (b) Corrected | 5.0% | 9.6% | 12.4% | 14.9% | 20.7% | 26.4% |
| (c) $\mathcal{LS}$-difference | 4.5% | 3.5% | 3.3% | 3.5% | 3.6% | 3.0% |
| (d) $\mathcal{LS}$-maximum | 14.9% | 31.5% | 38.7% | 47.2% | 56.5% | 68.4% |
| (d) Corrected | 5.0% | 11.8% | 17.6% | 21.5% | 29.5% | 38.9% |
| (e) Joint-CC | 4.6% | 5.0% | 6.4% | 9.3% | 19.2% | 29.9% |
| (f) Joint-TP | 4.6% | 4.7% | 5.5% | 7.5% | 15.5% | 24.2% |
| | ***a*** | | | | | |
| | 0.225 | 0.25 | 0.275 | 0.3 | 0.325 | 0.35 |
| (a) $\mathcal{LS}(\widehat{P}_m,\widehat{Q}_n)$ | 44.2% | 55.6% | 62.3% | 70.5% | 79.2% | 84.1% |
| (a) Corrected | 39.5% | 48.3% | 56.1% | 65.0% | 75.0% | 80.4% |
| (b) $\mathcal{LS}(\widehat{Q}_n,\widehat{P}_m)$ | 46.0% | 53.7% | 63.8% | 73.8% | 79.1% | 86.6% |
| (b) Corrected | 34.8% | 41.9% | 51.4% | 61.8% | 68.3% | 77.3% |
| (c) $\mathcal{LS}$-difference | 4.0% | 2.1% | 1.3% | 1.4% | 1.1% | 1.3% |
| (d) $\mathcal{LS}$-maximum | 77.6% | 86.9% | 92.8% | 96.5% | 98.9% | 99.7% |
| (d) Corrected | 51.5% | 61.1% | 71.9% | 81.0% | 90.5% | 95.7% |
| (e) Joint-CC | 44.5% | 63.5% | 78.0% | 89.7% | 95.5% | 98.8% |
| (f) Joint-TP | 38.5% | 58.0% | 71.8% | 84.6% | 92.2% | 97.4% |

functional case. The same effect has been observed in Figure 1 for the univariate case.

Furthermore, while the power of the $\mathcal{LS}$-difference statistic is satisfactory in this example, the power against a mere location shift $(B_t+0.15)_{t\in[0,1]}$ breaks down with the empirical power below 5% (falling again into the blind spot of the statistic), see Table 4. Also, from Table 4, we can observe that the Joint-CC and the Joint-TP tests exhibit a high power for large location shifts, being even higher than the $\mathcal{LS}$-maximum size corrected test.

Consequently, the two joint tests are the preferable ones, due to their size and power behaviour; and the huge blind spots the two projection statistics exhibit. In fact, the Joint-TP and Joint-CC test have similar size and power behaviour with a slightly higher power of the Joint-CC test which is due to the larger rejection region, as it has a faster contraction. In the rest of the simulation study, we restrict our attention to the Joint-TP test due to its conservative contraction rate.

Table 5

*Rejection rates for $(B_t)_{t\in[0,1]}$ (distribution $P$) against $(b\cdot B_t+0.15)_{t\in[0,1]}$ (distribution $Q$) for $b\in\{0.85,0.9,0.95,0.99\}$ and the statistics $\mathcal{LS}(\widehat{P}_m,\widehat{Q}_n)$ (rows (a)), $\mathcal{LS}(\widehat{Q}_n,\widehat{P}_m)$ (rows (b)), $\mathcal{LS}$-difference (row (c)), $\mathcal{LS}$-maximum (rows (d)), Joint-CC (row (e)) and Joint-TP (row (f)). In rows (a), (b) and (d) the size corrected power for the corresponding statistic is also included (the respective rejection rate under the null hypothesis is that given in the first column of Table 4, under $a=0$). The statistics are based on the functional integrated Tukey depth and $m=n=100$ observations in each sample.*

| | $b$ | | | |
|---|---|---|---|---|
| | 0.99 | 0.95 | 0.9 | 0.85 |
| (a) $\mathcal{LS}(\widehat{P}_m,\widehat{Q}_n)$ | 19.5% | 10.8% | 6.1% | 8.6% |
| (a) Corrected | 15.5% | 8.4% | 4.2% | 6.8% |
| (b) $\mathcal{LS}(\widehat{Q}_n,\widehat{P}_m)$ | 32.5% | 46.8% | 67.8% | 89.1% |
| (b) Corrected | 20.4% | 32.7% | 55.3% | 82.1% |
| (c) $\mathcal{LS}$-difference | 3.6% | 6.7% | 17.7% | 41.6% |
| (d) $\mathcal{LS}$-maximum | 50.0% | 55.9% | 71.2% | 89.7% |
| (d) Corrected | 22.9% | 31.1% | 50.0% | 76.6% |
| (e) Joint-CC | 10.5% | 15.4% | 24.6% | 48.5% |
| (f) Joint-TP | 9.2% | 13.3% | 22.4% | 46.6% |

Table 6

*Rejection rates for $\mathbf{N}((\mathbf{0},\mathbf{0})^{\mathbf{T}},\mathbf{I_2})$ (distribution $P$) against $\mathbf{N}((\mathbf{c},\mathbf{c})^{\mathbf{T}},\mathbf{I_2})$ (distribution $Q$) for $c\in\{0,0.2,0.4,0.6,0.8\}$ and the statistics $\mathcal{LS}(\widehat{P}_m,\widehat{Q}_n)$ (rows (a)), $\mathcal{LS}(\widehat{Q}_n,\widehat{P}_m)$ (rows (b)), $\mathcal{LS}$-difference (row (c)), $\mathcal{LS}$-maximum (rows (d)), Joint-CC (row (e)) and Joint-TP (row (f)). In rows (a), (b) and (d) the size corrected power for the corresponding statistic is also included (the respective rejection rate under the null hypothesis is that given in the first column, under $c=0$). All statistics are computed with the Tukey depth in $\mathbb{R}^2$ and $m=n=100$ observations in each sample.*

| | $c$ | | | | |
|---|---|---|---|---|---|
| | 0 | 0.2 | 0.4 | 0.6 | 0.8 |
| (a) $\mathcal{LS}(\widehat{P}_m,\widehat{Q}_n)$ | 10.2% | 16.1% | 41.7% | 78.3% | 96.7% |
| (a) Corrected | 5.0% | 9.1% | 27.7% | 64.8% | 94.0% |
| (b) $\mathcal{LS}(\widehat{Q}_n,\widehat{P}_m)$ | 14.0% | 20.1% | 41.0% | 76.2% | 96.8% |
| (b) Corrected | 5.0% | 8.7% | 24.5% | 60.9% | 90.9% |
| (c) $\mathcal{LS}$-difference | 4.8% | 4.9% | 3.2% | 2.7% | 2.5% |
| (d) $\mathcal{LS}$-maximum | 23.7% | 35.7% | 74.6% | 98.2% | 100.0% |
| (d) Corrected | 5.0% | 9.4% | 31.6% | 83.2% | 99.3% |
| (e) Joint-CC | 4.8% | 6.0% | 37.4% | 93.2% | 99.9% |
| (f) Joint-TP | 4.8% | 5.6% | 31.0% | 89.7% | 99.8% |

### *C.2. Joint-TP test for different functional depths*

In order to further investigate the small sample behaviour of the Joint-TP test, we consider two different models for generating functional data: Model 2 represents smooth functional data and Model 3 smooth fluctuating functional data. This is in addition to Model 1 (standard Brownian motion), which represents a way to generate noisy non-smooth functional data. For the smooth models, let

Table 7

*Rejection rates for* $\mathbf{N}((\mathbf{0},\mathbf{0})^{\mathbf{T}},\mathbf{I_2})$ *(distribution* $P$*) against* $\mathbf{N}((\mathbf{0.2},\mathbf{0.2})^{\mathbf{T}},\widetilde{\mathbf{c}}\cdot\mathbf{I_2})$ *(distribution* $Q$*) for* $\widetilde{c}\in\{1.01\}\cup\{1+0.05i : i\in\mathbb{N}_0$ *with* $i\leqslant 8\}$ *and the statistics* $\mathcal{LS}(\widehat{P}_m,\widehat{Q}_n)$ *(rows (a)),* $\mathcal{LS}(\widehat{Q}_n,\widehat{P}_m)$ *(rows (b)),* $\mathcal{LS}$*-difference (row (c)),* $\mathcal{LS}$*-maximum (rows (d)), Joint-CC (row (e)) and Joint-TP (row (f)). In rows (a), (b) and (d) the size corrected power for the corresponding statistic is also included (the respective rejection rate under the null hypothesis is that given in the first column of Table 6, under* $c=0$*). All statistics are computed with the Tukey depth in* $\mathbb{R}^2$ *and* $m=n=100$ *observations in each sample.*

| | $\widetilde{c}$ | | | | | | | | | |
|---|---|---|---|---|---|---|---|---|---|---|
| | 1 | 1.01 | 1.05 | 1.1 | 1.15 | 1.2 | 1.25 | 1.3 | 1.35 | 1.4 |
| (a) $\mathcal{LS}(\widehat{P}_m,\widehat{Q}_n)$ | 17.1% | 17.7% | 28.4% | 36.6% | 45.5% | 58.5% | 61.7% | 71.9% | 78.2% | 82.2% |
| (a) Corrected | 9.1% | 8.2% | 13.5% | 23.1% | 28.8% | 40.7% | 46.1% | 55.2% | 65.6% | 68.4% |
| (b) $\mathcal{LS}(\widehat{Q}_n,\widehat{P}_m)$ | 20.0% | 19.1% | 12.1% | 9.6% | 6.8% | 6.1% | 7.3% | 10.7% | 15.4% | 18.6% |
| (b) Corrected | 8.7% | 8.1% | 5.4% | 2.6% | 1.3% | 1.6% | 1.8% | 3.6% | 5.6% | 6.7% |
| (c) $\mathcal{LS}$-difference | 5.0% | 4.8% | 4.8% | 9.8% | 12.2% | 19.5% | 26.3% | 33.9% | 44.1% | 51.1% |
| (d) $\mathcal{LS}$-maximum | 36.6% | 36.3% | 39.7% | 44.7% | 49.6% | 60.6% | 62.9% | 73.2% | 78.6% | 82.3% |
| (d) Corrected | 9.4% | 8.7% | 9.5% | 15.1% | 18.2% | 26.6% | 32.1% | 40.9% | 51.7% | 57.0% |
| (e) Joint-CC | 6.4% | 6.0% | 5.6% | 10.9% | 12.5% | 20.0% | 27.3% | 35.6% | 44.6% | 51.5% |
| (f) Joint-TP | 6.1% | 5.5% | 5.5% | 10.5% | 12.4% | 20.0% | 27.1% | 34.6% | 44.6% | 51.4% |

us consider the Fourier basis functions

$$e_j(t)=\begin{cases}\sin(j\pi t), & j \text{ odd},\\ \cos(j\pi t), & j \text{ even}\end{cases}$$

of $L^2[0,1]$, where the constant function is excluded in order to achieve a better comparability of the models. Model 2 is $(C_t)_{t\in[0,1]}$ with

$$C_t=\sum_{l=1}^{20} W_l\cdot e_l(t),\qquad W_l \overset{\text{ind.}}{\sim} N\left(0,3^{-l}\right).$$

Model 3 is $(D_t)_{t\in[0,1]}$ with

$$D_t=\sum_{l=1}^{20} Z_l\cdot e_l(t),\qquad Z_l \overset{\text{ind.}}{\sim} N\left(0,\frac{J}{l}\right) \text{ with } J=\sum_{j=1}^{20}\frac{1}{j}=\frac{55835135}{15519504}\approx 3.6.$$

Here, we follow [6], who use the variance $3^{-l}$ to mimic a fast decay of the eigenvalues of the covariance operator and the variance $J/l$ for a slow decay.

Figure 6 shows some realizations to illustrate Model 2 (central panels) and Model 3 (right panels), in addition to Model 1 (left panels). In each of the panels we have represented 10 functional observations of the respective null distribution in grey: In each of the three columns corresponding to the three models, only the coloured sample (under the alternative) varies while the grey observations (corresponding to the null distribution) are the same in each panel. The plots clearly show that only Model 2 and 3 result in smooth curves, with the curves from Model 2 being smoother.

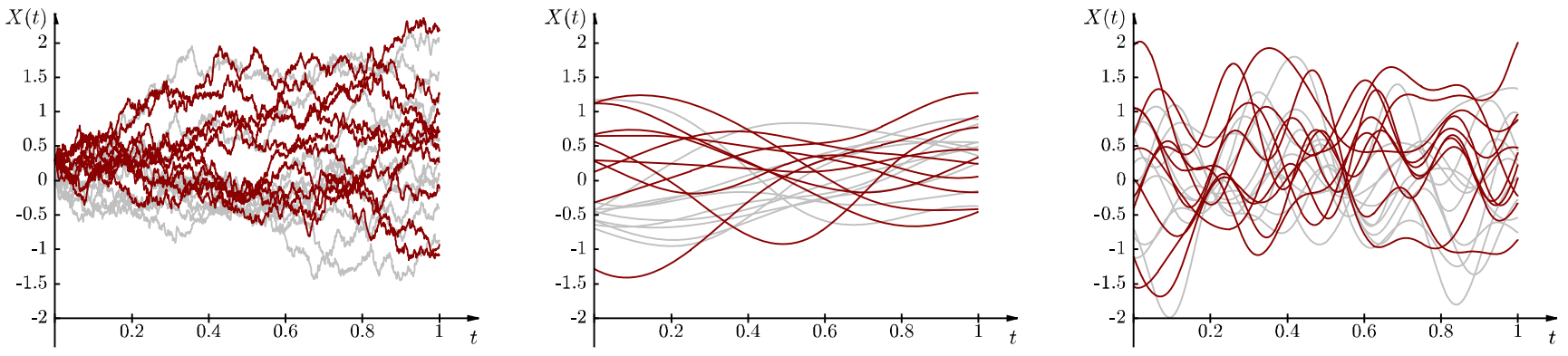

(a) $(B_t)_{t\in[0,1]}$, $(B_t+0.3)_{t\in[0,1]}$ (b) $(C_t)_{t\in[0,1]}$, $(C_t+0.3)_{t\in[0,1]}$ (c) $(D_t)_{t\in[0,1]}$, $(D_t+0.3)_{t\in[0,1]}$

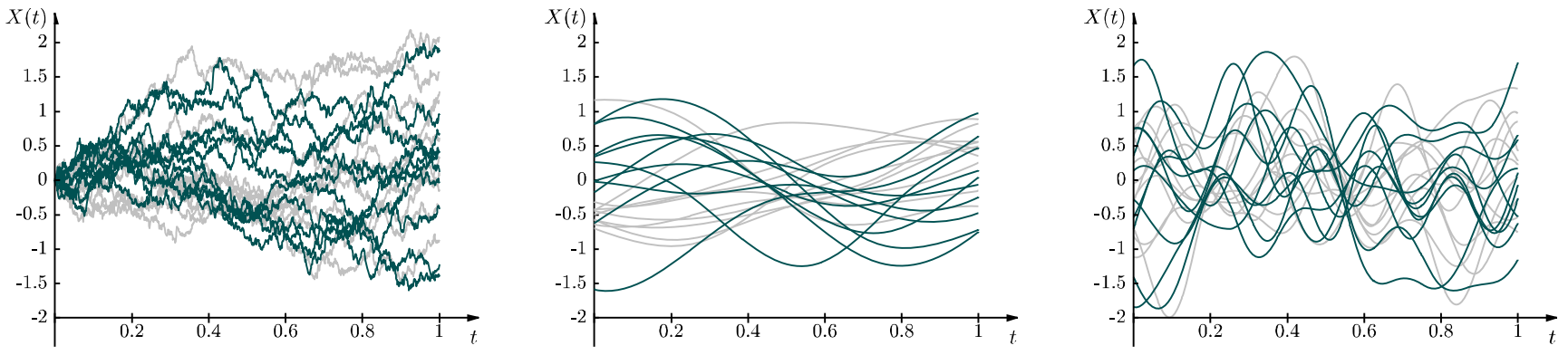

(d) $(B_t)_{t\in[0,1]}$, $(B_t+0.3\cdot\sin(2\pi t))_{t\in[0,1]}$ (e) $(C_t)_{t\in[0,1]}$, $(C_t+0.3\cdot\sin(2\pi t))_{t\in[0,1]}$ (f) $(D_t)_{t\in[0,1]}$, $(D_t+0.3\cdot\sin(2\pi t))_{t\in[0,1]}$

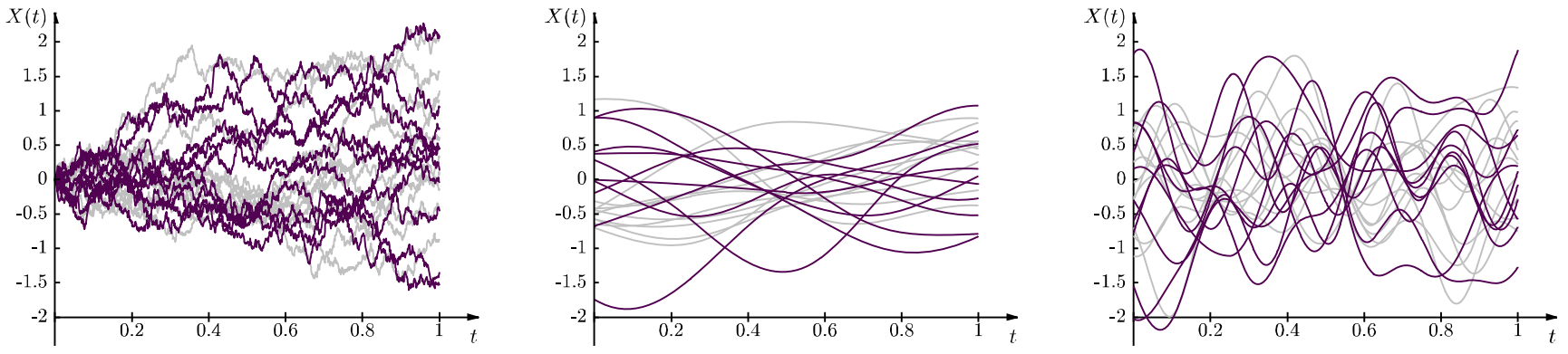

(g) $(B_t)_{t\in[0,1]}$, $(1.1\cdot B_t)_{t\in[0,1]}$ (h) $(C_t)_{t\in[0,1]}$, $(1.1\cdot C_t)_{t\in[0,1]}$ (i) $(D_t)_{t\in[0,1]}$, $(1.1\cdot D_t)_{t\in[0,1]}$

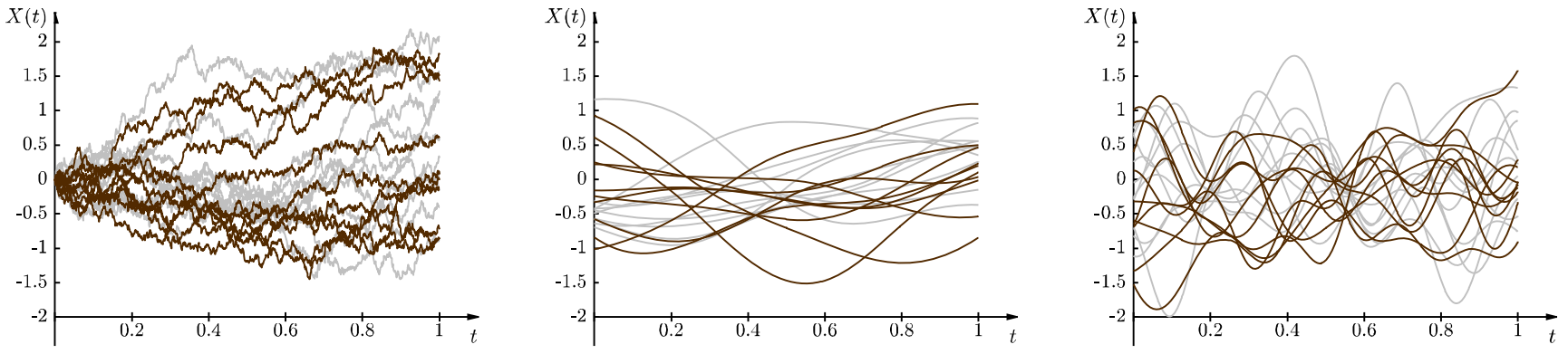

(j) $(B_t)_{t\in[0,1]}$, $\left(\frac{B_t+\sqrt{1-2\cdot 0.3}\widetilde{B}_t}{\sqrt{2}}+\frac{\sqrt{32}\cdot 0.3\cdot t(t-1)}{\sqrt{2}}\right)_{t\in[0,1]}$ (k) $(C_t)_{t\in[0,1]}$, $\left(\frac{C_t+\sqrt{1-2\cdot 0.3}\widetilde{C}_t}{\sqrt{2}}+\frac{\sqrt{32}\cdot 0.3\cdot t(t-1)}{\sqrt{2}}\right)_{t\in[0,1]}$ (l) $(D_t)_{t\in[0,1]}$, $\left(\frac{D_t+\sqrt{1-2\cdot 0.3}\widetilde{D}_t}{\sqrt{2}}+\frac{\sqrt{32}\cdot 0.3\cdot t(t-1)}{\sqrt{2}}\right)_{t\in[0,1]}$

FIG 6. 10+10 *functional observations generated with Models 1 (left panels), 2 (central panels) and 3 (right panels). Each row of panels corresponds to a class of alternatives: the grey curves are generated from the first distribution in each panel caption – the coloured curves from the second one. The first row depicts a simple location shift. In the second row there is a more advanced location shift by adding a scaled sine curve to the data. For the third row, a constant multiplicative scaling factor was added. The last row contains a simultaneous (more complex) location-scale difference, where* $\widetilde{B}_t$ *(respectively* $\widetilde{C}_t$, $\widetilde{D}_t$*) is an independent copy of* $B_t$ *(respectively* $C_t$, $D_t$*).*

In this section, we investigate the impact of the selection of depth function on the size and power of our proposed testing procedure. We consider the integrated depth, the $h$-depth, the spatial depth [11], the lens metric depth [39], the random Tukey depth [20] and the random projection depth [22]. As it can be observed from (14), the integrated depth requires of a multivariate depth for its computation. We make use of two popular options: Tukey and simplicial depths which are here reduced to univariate depths as in this section we simulate univariate functional data. Nevertheless, the Tukey and the simplicial depth with respect to the same finite probability measure in $\mathbb{R}$ do not necessarily lead to the same ordering as illustrated by the following example.

**Example C.1.** *Let $\{x_1, ..., x_{10}\}$ be a set of real numbers with $x_i < x_{i+1}$ for $i = 1, ..., 9$ and $P$ the discrete uniform distribution on the set. The Tukey depth leads to an inside-out ordering, i.e. $D_T(x, P) \leqslant D_T(y, P)$ for any $x, y \in \mathbb{R}$ with $|x - (x_5 + x_6)/2| \geqslant |y - (x_5 + x_6)/2|$. In contrast, for the simplicial depth (16) we have $D_S(x_3, P) = 24/45 > D_S(y, P) = 21/45$ for any $y \in (x_3, x_4)$. So, the simplicial depth with respect to finite probability measures does not necessarily entail an inside-out ordering.*

The simplicial depth does have a very good behaviour when computed with respect to continuous distributions [48]. To benefit from it, we also apply the integrated depth on a modified version of the simplicial depth, where we add some noise on each vertex of the simplices. In particular, we add a zero mean normal distribution with variance $10^{-16}$.

The $h$-depth has been extensively studied from a theoretical point for a fix $h$ [21, 88], and it has very good properties [61, 62]. However, it has been seen that its behaviour is better when $h$ depends on the distribution [22, 59, 63]. Thus, in our simulations we make use of $h = 1$ and $h$ equal to the 0.15-quantile of the $L^2$-distance among the curves of the underlying sample. Moreover, for the lens metric depth we make use of the $L^2([0, 1])$-metric. Also, we consider the random Tukey depth based on 2 projections and the random projection depth based on 10 projections.

The rejection rates under the null hypothesis for different sample sizes are displayed in Table 8 for Model 1, 2 and 3. The size behaviour displayed in the tables seems appropriate, since our proposed procedure is based on asymptotics; with the exception of the integrated simplicial depth in Model 3, which has an obvious size problem in small sample sizes. This size problem is related to the simplicial depth ordering explained in Example C.1. Table 8 clearly indicates that the proposed modification of the simplicial depth solves the problem resulting in a reasonable size behaviour.

To study the power of the test, we consider the following alternatives:

a simple location shift

$$(M_t + a)_{t \in [0,1]}, \tag{74}$$

a more complex location difference

$$(M_t + a \cdot \sin(2\pi t))_{t \in [0,1]}, \tag{75}$$

TABLE 8

*Rejection rates for Model 1-3 under $H_0$ for different pairs of sample sizes with several depth functions: the integrated Tukey depth (IntegrTukey), the integrated simplicial depth (IntegrSimplicial), the integrated modified simplicial depth (IntegrSimplicialMod), the $h$-depth with $h = 1$ ($h$-const), the $h$-depth with data-adaptively chosen $h$ ($h$-adaptive), the spatial depth (Spatial), the lens metric depth (LensMetric), the random Tukey depth based on two random projections (RandomTukey) and the random projection depth based on ten random projections (RandomProjection). For smaller samples the deviation from 5% is slightly larger because of the decision rule based on asymptotics. Unbalanced and small sample sizes can lead to rejection rates above 5%, e.g. with the lens metric depth.*

| **Model 1** | | | | | | |
|---|---|---|---|---|---|---|
| | ***m/n*** | | | | | |
| | 50/50 | 75/75 | 100/100 | 75/150 | 500/500 | 1000/1000 |
| IntegrTukey | 4.6% | 4.3% | 4.6% | 4.1% | 5.2% | 4.4% |
| IntegrSimplicial | 4.7% | 3.4% | 4.2% | 4.4% | 4.4% | 4.8% |
| IntegrSimplicialMod | 4.4% | 4.8% | 4.8% | 4.4% | 4.6% | 4.8% |
| $h$-const | 4.2% | 4.3% | 4.6% | 4.1% | 5.0% | 4.2% |
| $h$-adaptive | 3.5% | 3.7% | 3.9% | 3.9% | 4.6% | 4.7% |
| Spatial | 4.7% | 4.0% | 4.5% | 5.0% | 4.8% | 4.5% |
| LensMetric | 5.9% | 5.3% | 5.1% | 5.6% | 5.2% | 4.2% |
| RandomTukey | 4.7% | 5.0% | 4.6% | 4.4% | 5.3% | 4.4% |
| RandomProjection | 3.5% | 4.3% | 4.4% | 4.1% | 4.9% | 4.2% |
| **Model 2** | | | | | | |
| | ***m/n*** | | | | | |
| | 50/50 | 75/75 | 100/100 | 75/150 | 500/500 | 1000/1000 |
| IntegrTukey | 4.1% | 5.4% | 4.3% | 4.8% | 4.9% | 5.0% |
| IntegrSimplicial | 4.3% | 5.0% | 4.4% | 5.3% | 4.7% | 5.0% |
| IntegrSimplicialMod | 4.2% | 5.7% | 4.7% | 3.8% | 4.7% | 5.1% |
| $h$-const | 4.3% | 4.9% | 4.4% | 3.8% | 4.3% | 4.7% |
| $h$-adaptive | 4.6% | 4.6% | 4.3% | 4.7% | 4.4% | 4.7% |
| Spatial | 4.5% | 4.9% | 4.3% | 4.5% | 4.6% | 4.9% |
| LensMetric | 4.7% | 5.2% | 4.6% | 4.7% | 4.4% | 4.8% |
| RandomTukey | 5.1% | 5.6% | 5.2% | 4.1% | 4.5% | 5.0% |
| RandomProjection | 5.0% | 4.1% | 3.9% | 4.5% | 4.1% | 5.0% |
| **Model 3** | | | | | | |
| | ***m/n*** | | | | | |
| | 50/50 | 75/75 | 100/100 | 75/150 | 500/500 | 1000/1000 |
| IntegrTukey | 6.3% | 5.4% | 5.2% | 6.9% | 4.4% | 4.3% |
| IntegrSimplicial | 84.1% | 24.2% | 6.5% | 10.5% | 4.0% | 4.5% |
| IntegrSimplicialMod | 5.1% | 5.3% | 5.0% | 5.7% | 4.0% | 4.6% |
| $h$-const | 5.0% | 5.0% | 5.5% | 5.1% | 4.6% | 5.0% |
| $h$-adaptive | 4.7% | 4.6% | 5.3% | 5.2% | 4.3% | 4.8% |
| Spatial | 5.0% | 5.3% | 5.7% | 5.5% | 4.5% | 4.8% |
| LensMetric | 5.0% | 5.7% | 5.7% | 5.4% | 5.0% | 4.9% |
| RandomTukey | 5.8% | 7.6% | 4.9% | 6.1% | 4.1% | 5.2% |
| RandomProjection | 4.4% | 4.9% | 5.6% | 4.0% | 4.9% | 6.1% |

a scale difference

$$(a \cdot M_t)_{t\in[0,1]} \tag{76}$$

and a simultaneous difference in the location as well as in the second order

structure

$$\left(\frac{M_t + \sqrt{1-2a}\widetilde{M}_t + \sqrt{32}at(t-1)}{\sqrt{2}}\right)_{t\in[0,1]}, \tag{77}$$

with $\widetilde{M}_t$ an independent copy of $M_t$. $M$ stands for $B$ when making use of Model 1, for $C$ when making use of Model 2 and $D$ when using Model 3. These alternatives depend on a constant $a \in \mathbb{R}\backslash\{0\}$ for (74), (75) and (77) and $a \in \mathbb{R}\backslash\{1\}$ for (76). To illustrate them, each panel of Figure 6 displays 10 functional observations in colour that were generated under an alternative hypothesis. The first row panels display in red the alternative in (74), the second (75) in green, the third (76) in blue and the fourth (77) in orange. For the illustration we have used $a = 0.3$ in (74), (75) and (77) and $a = 1.1$ in (76).

Figure 7 shows the power curves for the studied alternatives: (74) in the first row (plots (a)–(c)), (75) in the second (plots (d)–(f)), (76) in the third (plots (g)–(i)) and (77) in the fourth (plots (j)–(l)). The plots in the $i$-th column regard Model $i$, for $i = 1, 2, 3$.

When making use of a simple location shift in Models 1 (plot (a)) and 3 (plot (c)), the curves showing most power are those of the integrated simplicial depth (orange) followed by the integrated Tukey depth (navy blue). These two models are related in that there are multiple crossings among each two curves in a sample. They are different in that Model 3 results in smooth curves while Model 1 does not. However, the smoothness of the data does not affect the integrated depths as, in practice, the sample curves are observed over a grid of equally spaced points and the integrated depths are computed over the sum of the corresponding univariate depth at each grid point. For Model 1, these power curves are also just followed by the integrated modified simplicial depth power curve (red). Thus, the integrated depths are the ones resulting in most power for simple location shifts of the Brownian motion. We can argue that the location of the image of the curves corresponding to the alternative sample gets higher with higher values of $a$. Thus, for each fixed grid point, as $a$ increases, the depth value of a point corresponding to one of the samples has low univariate depth when computed with respect to the other sample. For Model 3 (plot (c)), the results of the integrated simplicial depth (orange) are affected by the high rejection rates under the null hypothesis, which are observable from Table 8. In this case we can consider the integrated Tukey depth (navy blue) a preferable choice, as the integrated modified simplicial depth (red) results in a power curve that is much lower than that of other depth functions.

The high rejection rates of the integrated simplicial depth (orange) under the null hypothesis in Model 3 are also observable in the plots of the other three alternatives (plots (f), (i) and (l)). The behaviour under the complex location difference (plot (f)) is very similar to that of the simple location shift (plot (c)) analysed above, resulting in the integrated Tukey depth (navy blue) being a preferable choice. For a simultaneous difference in the location as well as in the second order structure (plot (l)), the integrated simplicial depth power curve (orange) is only just above those of the adaptive $h$-depth (green), the spatial depth (cyan), the $h$-depth (yellow) and the lens metric depth (pink). These four power curves are indeed the ones that are higher for the scale difference (plot (i)). It is worth saying that the tests based on these depth functions are the ones that make use of the distance between the curves. Thus, in these two last cases we recommend the use of any of the four distance based depth functions – in particular, the adaptive $h$-depth.

For Model 2, the adaptive $h$-depth power curve (green) is generally higher irrespective of the studied alternative (plots (b), (e), (h) and (k)). The main difference between Model 2 and the other two models is that in Model 2 there are typically just two crossings between each two curves. This, results in that the norm between each two curves is a relevant measure, which for the adaptive $h$-depth is used both in the computation of the $h$ as well as in the depth itself, being amplified by the Gaussian kernel used.

The adaptive $h$-depth (green) is also a good choice for Model 1 for the alternatives given in (75), (76) and (77) (plots (d), (g) and (j)). In fact, it is the highest power curve for (76) and (77) (plots (g) and (j)) and just below the integrated simplicial depth power curve (orange) for (75) (plot (d)).

Summarizing, the best results are obtained by the distance based depth functions and the integrated depth functions. Among the distance based depth functions, the clear winner is the adaptive $h$-depth (green), which makes a double use of the distance between curves. For the integrated depth functions, the best results are obtained by the integrated simplicial depth (orange), despite the high rejection rates under the null hypothesis in Model 3. Some more conservative results are given by the integrated Tukey depth (navy blue). Furthermore, the random depth functions, random projection depth and random Tukey depth, provide non-relevant results.

#### *C.2.1. Special alternative: fixed location and second order structure, but different shape*

We investigate the ability of the Joint-TP test to detect differences in the shape of the data, i.e. finding differences in the third- or higher order structures. Therefore, we consider the models

$$\begin{aligned} S_t &= U, \text{ where } U \sim U(0,1), t \in [0,1] \text{ and} \\ T_t &= 5\cdot(t+V) - \lfloor 5\cdot(t+V)\rfloor, \text{ where } V \sim U(0,1), t \in [0,1]. \end{aligned}$$

An illustration is given in Figure 8, where 10 curves of each sample are depicted. The rejection rates in Table 9 indicate that the tests using integrated depths have less power than tests using depths from other families. The reason for this phenomenon is that integrated depths are computed by averaging pointwise computed values over the time axis, but the shape difference does not occur pointwise.

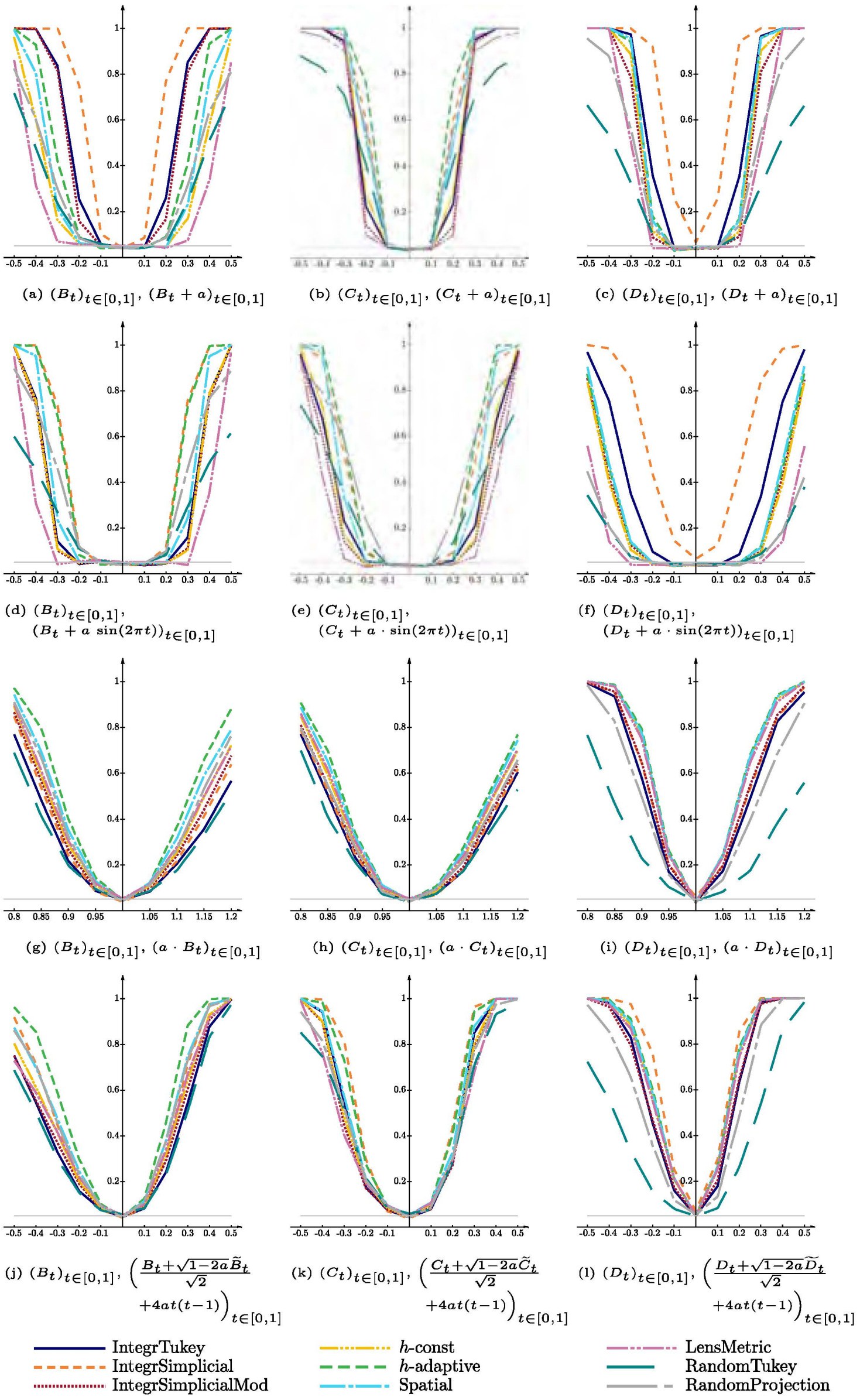

FIG 7. *Power-plots (sample sizes: $m = n = 100$) for several alternatives depending on the parameter $a$ (on the horizontal axis). Each curve represents the behaviour of the Joint-TP test with a particular depth function: the integrated Tukey depth (IntegrTukey), the integrated simplicial depth (IntegrSimplicial), the integrated modified simplicial depth (IntegrSimplicialMod), the h-depth with $h = 1$ (h-const), the h-depth with data-adaptively chosen h (h-adaptive), the spatial depth (Spatial), the lens metric depth (LensMetric), the random Tukey depth based on two random projections (RandomTukey) and the random projection depth based on ten random projections (RandomProjection).*

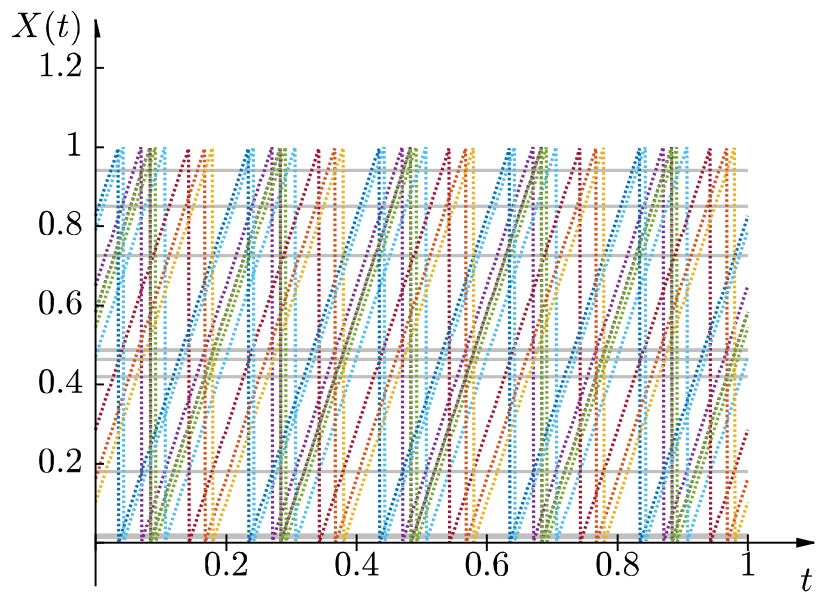

Fig 8. *Shape difference alternative: 10 horizontal curves following the distribution of $S_t$ (grey) and 10 dotted zigzag curves following the distribution of $T_t$ (coloured).*

Table 9

*Rejection rates for the shape difference alternative with sample sizes $m = n = 100$ and several depth functions: the integrated Tukey depth (IntegrTukey), the integrated simplicial depth (IntegrSimplicial), the integrated modified simplicial depth (IntegrSimplicialMod), the $h$-depth with $h = 1$ ($h$-const), the $h$-depth with data-adaptively chosen $h$ ($h$-adaptive), the spatial depth (Spatial), the lens metric depth (LensMetric), the random Tukey depth based on two random projections (RandomTukey) and the random projection depth based on ten random projections (RandomProjection).*

| **Model:** $(S_t)_{t\in[0,1]}$ against $(T_t)_{t\in[0,1]}$, $m = n = 100$ | |
|---|---|
| IntegrTukey | 2.9% |
| IntegrSimplicial | 46.2% |
| IntegrSimplicialMod | 46.4% |
| $h$-const | 56.8% |
| $h$-adaptive | 100.0% |
| Spatial | 100.0% |
| LensMetric | 99.6% |
| RandomTukey | 94.4% |
| RandomProjection | 100.0% |

### C.3. $\mathcal{LS}$-ellipsoidal statistic

Shi, Zhang and Fu [79] consider the following $\mathcal{LS}$-ellipsoidal test

$$\frac{12 \cdot m \cdot n}{m+n} \cdot \left( w \left( \mathcal{LS}(\widehat{P}_m, \widehat{Q}_n) - \frac{1}{2} \right)^2 + (1-w) \cdot \left( \mathcal{LS}(\widehat{Q}_n, \widehat{P}_m) - \frac{1}{2} \right)^2 \right) \tag{78}$$

with nonnegative weight $0 < w < 1$. The null asymptotic of this class of test statistic also follows immediately from Theorem 2.3. Figure 9 shows that the corresponding non-rejection regions are ellipsoidal and reports the empirical size in the situation as of Figure 2.

Clearly, the $\mathcal{LS}$-ellipsoidal tests eliminate some of the blind spots of the asymmetric $\mathcal{LS}$-statistics given by the projections. The best elimination is obtained for $w = 1/2$. In the other cases, the $\mathcal{LS}$-ellipsoidal test has unnecessary blind spots and a size problem due to its rejection region (see in particular (c) in

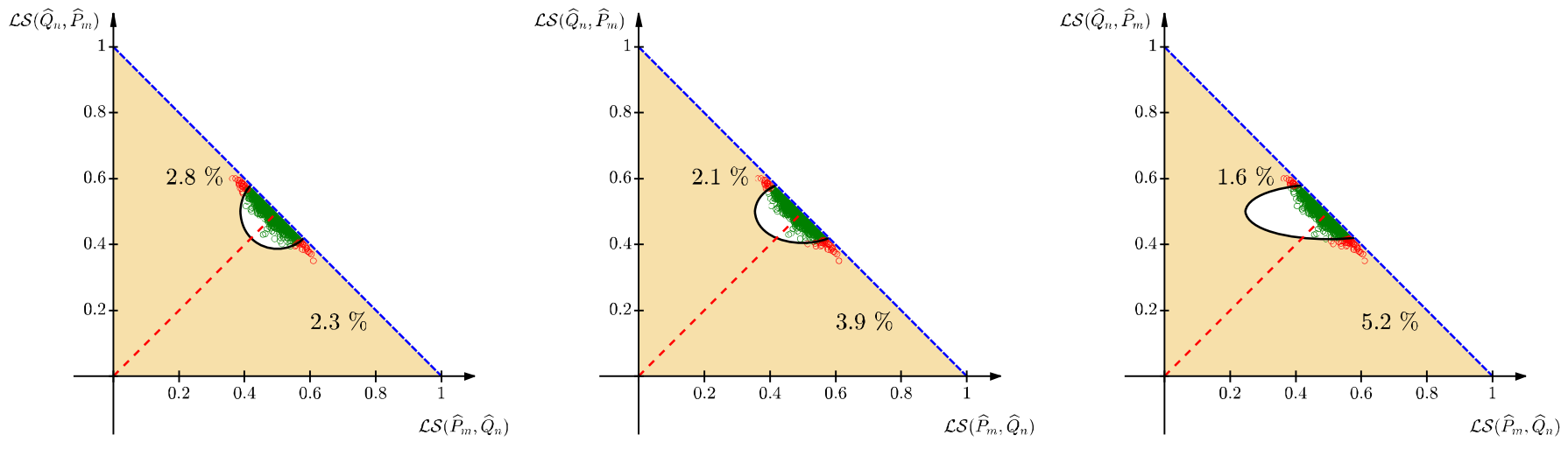

(a) $w_1 = w_2 = \frac{1}{2}$ (b) $w_1 = \frac{3}{10}, w_2 = \frac{7}{10}$ (c) $w_1 = \frac{1}{10}, w_2 = \frac{9}{10}$

Fig 9. *The circles are* 1000 *realisations of the* $\mathcal{LS}$*-tuple obtained with the univariate Tukey depth and* $U(0,1)$*-distributed observations, 100 in each sample, where the red dots (unlike the green ones) are rejected by the respective variation of* (78)*, where* $w_1, w_2$ *denote the chosen weights. The numbers in the plots give the empirical rejection rates to both sides of the red dashed main diagonal, while the shaded orange region shows the corresponding rejection region of the corresponding variation of the* $\mathcal{LS}$*-statistic.*

Figure 2). For $w = 1/2$ the non-rejection region is the inscribed circle of the non-rejection region of the Joint-TP test, eliminating some of the blind spots closer to the triangle point but at the cost of a possibly too large rejection region closer to the dark blue dash-dotted secondary diagonal. While this is not a problem for the example in Figure 2, in terms of the empirical size, in the situation of Model 3 (integrated Tukey depth, $m = n = 100$ observations in each sample) this leads to an enlarged empirical size of 8.8% as opposed to 5.2% for the Joint-TP test (see Table 8).

Shi, Zhang and Fu [79] give extensive simulation studies concerned with the empirical power of the $\mathcal{LS}$-maximum test and the $\mathcal{LS}$-ellipsoidal tests but do not take the empirical size into account. As such it is not surprising that they come to the conclusion that the $\mathcal{LS}$-maximum test is preferable over the other tests despite the fact that it does not hold the size to such an extend that it should not be used in practice.

### *C.4. On the location of the $\mathcal{LS}$-tuple*

Under the null hypothesis, it holds $\mathcal{LS}(\widehat{P}_m, \widehat{P}_m) = 1/2$, and thus, we expect the $\mathcal{LS}$-tuple to take its values close to the mid point of the unit square. Indeed, Theorem 2.3 entails that the $\mathcal{LS}$-tuple contracts asymptotically to the mid point of the unit square under the null hypothesis. With most depth functions, our simulations entailed that the $\mathcal{LS}$-tuple takes only values below the secondary diagonal of the unit square, but this is not always the case under both the null and the alternative hypothesis.

Under the null hypothesis this is illustrated using the lens metric depth in Figure 10. As the decision rule of the Joint-TP test holds for points outside the lower triangle (the lighter coloured areas are part of the rejection region), this effect does not affect the size of the Joint-TP test. In fact, Corollary 2.5 also allows to cut off the non-rejection region of the Joint-TP test above the secondary diagonal in the same way as proposed below the secondary diagonal. Under the alternative hypothesis, a similar effect occurs in Figure 11 (orange crosses) with the integrated Tukey depth. Nevertheless, in our extensive simulations we could not observe that the $\mathcal{LS}$-tuple takes any values in the unit square above the secondary diagonal that are far from this secondary diagonal. It is worth mentioning that there exist special cases, for which the $\mathcal{LS}$-tuple is guaranteed to take values below the secondary diagonal under the null hypothesis, see Theorem 2.6.

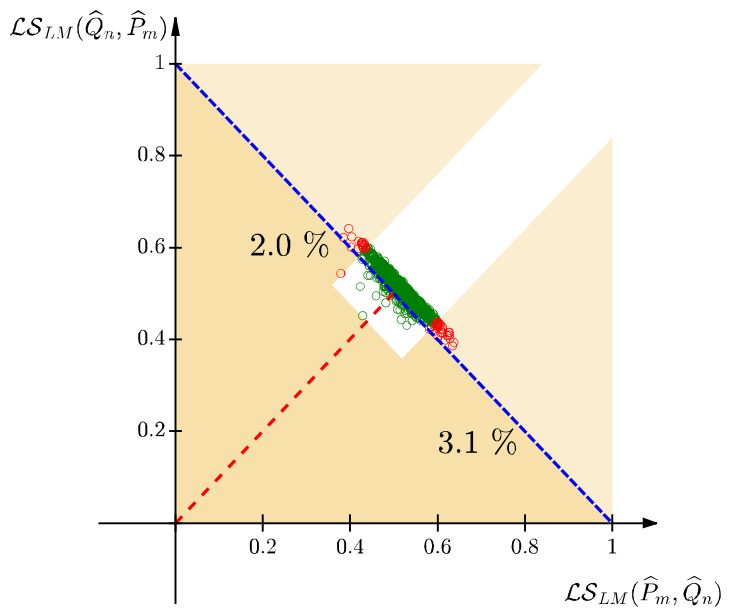

FIG 10. *For 1000 trials, it is plotted for $m, n = 100$ the $\mathcal{LS}$-tuple $(\mathcal{LS}(\widehat{P}_m, \widehat{Q}_n), \mathcal{LS}(\widehat{Q}_n, \widehat{P}_m))$ with the lens metric depth. Green circles: $P = Q = (B_t)_{t\in[0,1]}$.*

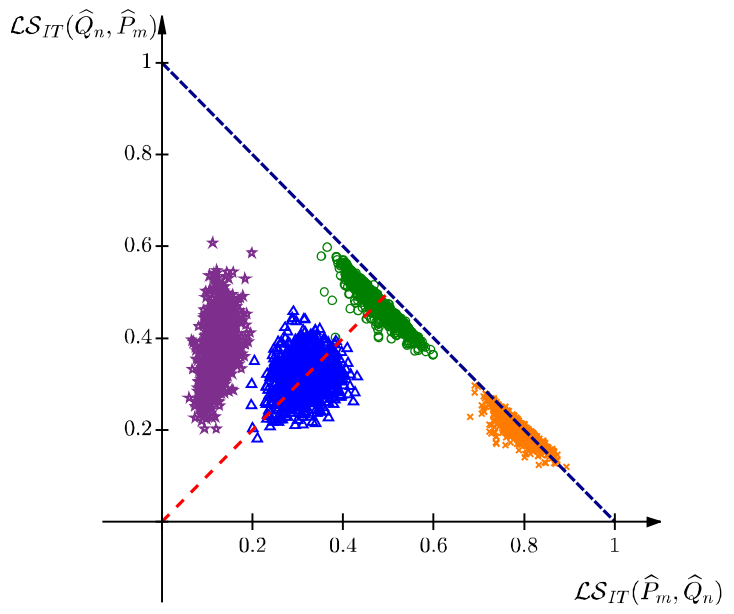

FIG 11. *For 1000 trials, it is plotted for $m, n = 100$ the $\mathcal{LS}$-tuple $(\mathcal{LS}(\widehat{P}_m, \widehat{Q}_n), \mathcal{LS}(\widehat{Q}_n, \widehat{P}_m))$ with the integrated Tukey depth. Green circles: $P = Q = (B_t)_{t\in[0,1]}$. Blue triangles: $P = (B_t)_{t\in[0,1]}, Q = (B_t + 0.5)_{t\in[0,1]}$. Orange crosses: $P = (B_t)_{t\in[0,1]}, Q = (B_t + 0.5)_{t\in[0,1]}$. Purple stars: $P = (B_t)_{t\in[0,1]}, Q = (2 \cdot B_t + 1)_{t\in[0,1]}$.*

Moreover, we want to point out that the value of the $\mathcal{LS}$-tuple may help to draw conclusions on the type of detected difference between the underlying probability distributions $P$ and $Q$ without applying further testing procedures, as illustrated in Figure 11: A location difference (blue triangles) leads to a shift of the points alongside the red dashed main diagonal. A scale decrease (orange crosses) causes a shift downwards the dark blue dash-dotted secondary diagonal, while a scale increase causes a shift upwards this diagonal. A simultaneous difference in location and scale (purple stars) leads to a shift alongside both diagonals. Moreover, the exact value of the $\mathcal{LS}$-tuple gives a hint on which sample lies deeper within the other one. Liu and Singh [52] originally proposed to make use of their version of the test for scale increases (under simultaneous location shifts). In future work, the above observations could also be used to construct unbiased one-sided tests for scale differences (under simultaneous location shifts) with better size behaviour by adapting the rejection regions correspondingly.

### C.5. Comparison with the state-of-the-art in functional two-sample tests

We compare three $\mathcal{LS}$-tests based on the adaptive $h$-depth, the integrated simplicial depth and the integrated Tukey depth with some other methods available in the literature, namely:

(𝔞) Two bootstrap-free ANOVA tests for detecting differences in the mean of functional data. One of them is an F-type test based on a naive and bias-reduced method of estimation [78, 89] and the other one is a globalised version of the pointwise F-test [90]. Both are implemented in the R-package ‘fdANOVA’ [40] and denoted as ‘FB’ respectively ‘GPF’ within the function *fanova.tests*.

(𝔟) A nonparametric two-sample test for functional data for detecting differences in the covariance operator based on PCA [36] and implemented in the R-package ‘fChange’ [84] which not maintained on CRAN at the moment (name of the function: *Cov_test*). For our purposes, we consider this test with five respectively ten principal components.

(𝔠) A nonparametric two-sample test based on the energy distance for multivariate data [83] from the R-package ‘energy’ [70] (name of the function: *eqdist.etest*). Note, that discretised functional data can be regarded as multivariate data.

(𝔡) A nonparametric ball divergence two-sample test for functional data [65] implemented in the R-package ‘Ball’ [92]. The decision rule of this test is based on the asymptotic behaviour of the test statistic. The name of the function in the R-package is *bd.test* and the there included method ‘limit’ is used.

The R-code for the integrated depths is our own implementation in order to unify the versions used for the simulations in this paper (for example, the R-package ‘fda.usc’ makes use of other estimators for these depths and this has an impact on the performance of the $\mathcal{LS}$-test). For computing the depth functions, we used own implementations that correspond to the definitions given in this paper, as there exist different estimators for most depth functions with respect to empirical probability measures which have an impact on the performance of the Joint-TP test. The $h$-depth coincedes with the implementation from the R-package ‘fda.usc’ [32].

We run a simulation study with the set of alternatives provided in Section C.2. The results are displayed in Figure 12. Obviously, the ANOVA-procedures ($\mathfrak{a}$) – in particular the global one – that are designed for detecting differences in the location generally have higher power than the other tests for this kind of alternative (Figure 12, (a)–(f)), but they completely fail to detect scale differences (Figure 12, (g)–(i)). Even in the case of a simultaneous location-scale difference, they are a good choice as long as the location difference dominates (Figure 12, (j)–(l)). The PCA-based tests for differences in the covariance operator ($\mathfrak{b}$) cannot detect any location differences (Figure 12, (a)–(f)). For pure scale differences, the $\mathcal{LS}$-test with the data-adaptive $h$-depth shows similar or even better performance than all the other tests. This is remarkable, as the PCA-tests are particularly designed for detecting these kind of alternatives (Figure 12, (g)–(i)). The nonparametric two-sample test for multivariate data ($\mathfrak{c}$) has high power against differences in the location (Figure 12, (a)–(f)), but difficulties to discriminate between alternatives in which only the scaling is different (Figure 12, (g)–(i)). The ball divergence statistic ($\mathfrak{d}$) has less power than the three versions of the Joint-TP test when it comes to pure second order differences in non-smooth or smooth data with less fluctuation (Figure 12, (g)–(i)).

In presence of location differences, it seems to perform better than the Joint-TP tests with smooth data (Figure 12, (b)–(c) and (e)–(f)), but in contrast, with non-smooth data the Joint-TP test with the integrated simplicial depth is more sensitive (Figure 12, (a) and (d)). All in all, there exist alternatives in which some of these test classes dominate all the other tests, in particular, if they were designed for that particular type of alternatives. Overall, the Joint-TP test has a competitive performance in particular with non-smooth raw data or for shifts that are dominated by a difference in the second order structure.

## *C.6. Performance in the presence of an additive outlier*

Depth functions are known for being robust and, indeed, the notion of functional depth [61] does imply qualitative robustness. Functional data may contain different types of outlying trajectories, including magnitude and/or shape outliers [3]. While, generally, lower depth values correspond to potential outliers, not all types of outliers will get low depth values and which types these are depends on the depth function. As the $\mathcal{LS}$-test is based on depth functions it is of interest to investigate its performance in the presence of outliers. We do this in comparison to the two ANOVA-methods that were also considered in the previous section. Here, we concentrate on magnitude outliers generated by adding a constant function of value +50 to a given curve; and refer them as additive outliers. These additive outliers can be expected to have a low depth value with respect to all considered depth functions and will also be outlying for the ANOVA procedures.

More precisely, we consider the following scenarios with 100 observations in each of the two samples in each case.

(a) *Null hypothesis without outliers:* $(B_t)_{t\in[0,1]}$ against $(B_t)_{t\in[0,1]}$.

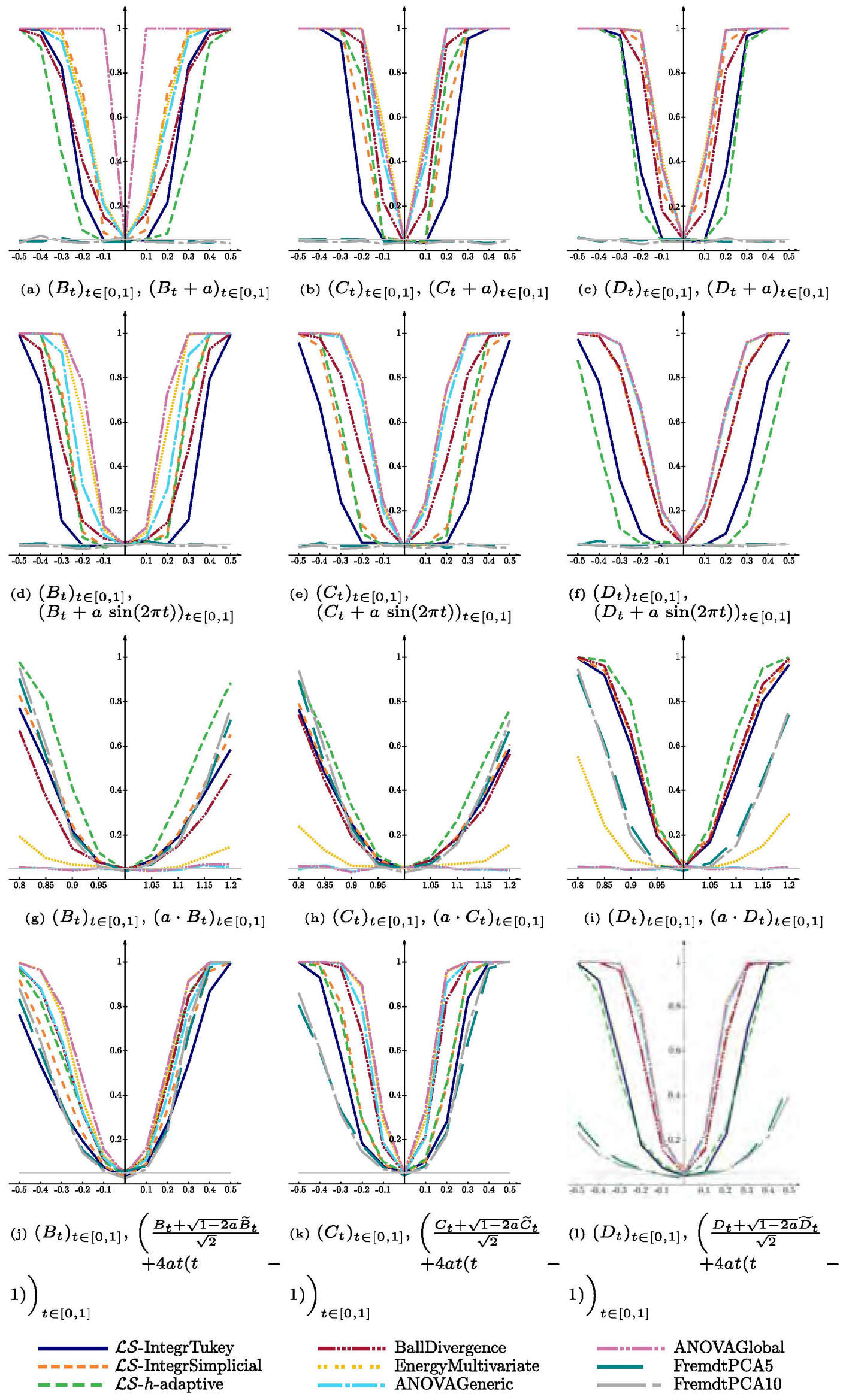

Fig 12. *Power-plots (sample sizes: $m = n = 100$) for several alternatives depending on the parameter $a$ (on the horizontal axis). Each curve represents the behaviour of a particular test. $\mathcal{LS}$-IntegrTukey, $\mathcal{LS}$-IntegrSimplicial and $\mathcal{LS}$-$h$-adaptive stand for the Joint-TP test with the integrated Tukey, the integrated simplicial and the data-adaptively $h$-depth. BallDivergence is a functional nonparametric two-sample test, EnergyMultivariate a multivariate nonparametric two-sample test (that can be also applied to functional data on a time grid), ANOVAGeneric as well as ANOVAGlobal are ANOVA-tests for functional data, and FremdtPCA5 respectively FremdtPCA10 are two versions of a PCA-based test (with 5 respectively 10 principal components) for detecting scale differences.*

(b) *Null hypothesis with additive outlier in both samples:* $(B_t)_{t\in[0,1]}$ against $(B_t)_{t\in[0,1]}$, where the first observation in each of the samples was shifted by $+50$.

(c) *Null hypothesis with additive outlier in one sample:* $(B_t)_{t\in[0,1]}$ against $(B_t)_{t\in[0,1]}$, where the first observation in the first sample was shifted by $+50$.

(d) *Alternative hypothesis without outliers:* $(B_t)_{t\in[0,1]}$ against $(B_t+0.25)_{t\in[0,1]}$.

(e) *Alternative hypothesis with additive outlier in both samples:* $(B_t)_{t\in[0,1]}$ against $(B_t+0.25)_{t\in[0,1]}$, where the first observation in each of the samples was shifted by $+50$.

(f) *Alternative hypothesis with additive outlier in the sample with lower mean:* $(B_t)_{t\in[0,1]}$ against $(B_t+0.25)_{t\in[0,1]}$, where the first observation in the first sample was shifted by $+50$.

(g) *Alternative hypothesis with additive outlier in the sample with higher mean:* $(B_t)_{t\in[0,1]}$ against $(B_t+0.25)_{t\in[0,1]}$, where the first observation in the second sample was shifted by $+50$.

Table 10

*Empirical rejection rates for the scenarios (a)–(g), where (a) and (d) do not contain any outliers.*

| | Null hypotheses | | | Alternative hypotheses | | | |
|---|---|---|---|---|---|---|---|
| | (a) | (b) | (c) | (d) | (e) | (f) | (g) |
| $\mathcal{LS}$-IntegrTukey | 4.9% | 4.1% | 4.4% | 55.3% | 51.9% | 48.6% | 56.9% |
| $\mathcal{LS}$-IntegrSimplicial | 4.3% | 3.5% | 3.7% | 95.6% | 95.1% | 93.9% | 95.9% |
| $\mathcal{LS}$-$h$-adaptive | 4.5% | 3.4% | 3.6% | 18.3% | 16.6% | 17.7% | 17.8% |
| ANOVAGeneric | 4.8% | 0.0% | 0.0% | 84.9% | 0.0% | 0.0% | 100.0% |
| ANOVAGlobal | 5.2% | 0.0% | 0.0% | 100.0% | 0.0% | 0.0% | 100.0% |

The empirical rejection rates at nominal level $\alpha = 0.05$ based on 1000 repetitions of each scenario are given in Table 10. The empirical size of the ANOVA procedures breaks down to 0, if either both or only one of the samples contain such an outlier, while the $\mathcal{LS}$-tests become only slightly more conservative. Similarly, under the alternatives, the ANOVA procedures are not robust towards outliers, in the sense that their empirical power breaks down to 0, unless the additive outlier works in the same direction as the level shift of the alternative (Scenario (g)). On the other hand, the $\mathcal{LS}$-tests have a comparable empirical power through all alternatives.

If not only one but five outliers are added to both samples in an analogous fashion, the $\mathcal{LS}$-tests become more conservative with rejection rates of 2.4% for $\mathcal{LS}$-Tukey, 2.0% for $\mathcal{LS}$-Simplicial, 2.1% for $\mathcal{LS}$-$h$-adaptive and – as before – 0 for both of the ANOVA-procedures. The picture changes if only one of the two samples contains five outliers of this type. In that case, the ANOVA-based tests reject in 100% of the cases, while the $\mathcal{LS}$-tests become unbiased with rejection rates of 6.7% for $\mathcal{LS}$-Simplicial, 7.2% for $\mathcal{LS}$-$h$-adaptive and 8.1% for $\mathcal{LS}$-Tukey. This effect is not surprising given that all tests search for distribu-tional differences in the two samples: If only one of the two samples contains several additive outliers, the distributions of both samples do indeed differ. Nevertheless, the $\mathcal{LS}$-tests are much less affected by this than the ANOVA-based tests.

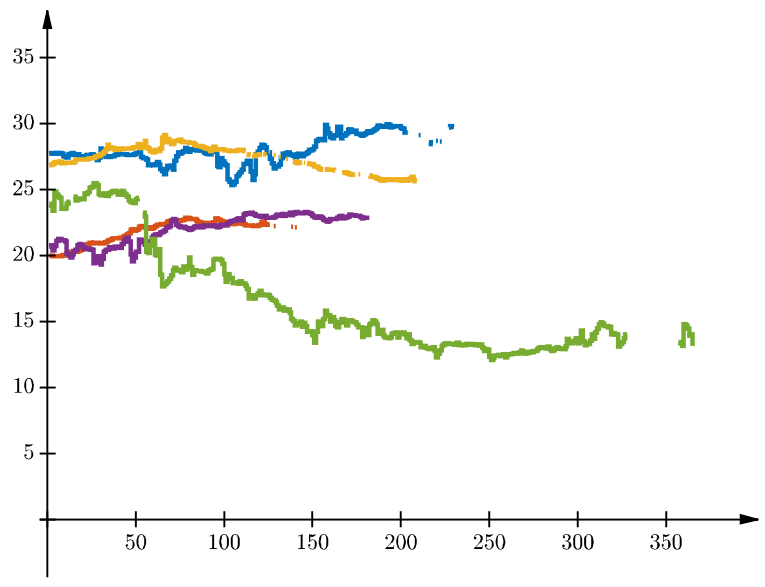

FIG 13. *Five temperature curves from the dataset of the year 2021.*

### C.7. Analysis of ocean drifter temperature data containing curves not observed over the whole grid domain

Commonly, functional data are not observed over an equally spaced dense grid. For instance, a curve can be observed over a grid domain with sparse and dense regions or it could be observed only on a few grid points, which may or may not be concentrated on certain regions of the domain.

As a proof of concept, we apply our methodology again to the ocean drifter data as in Section 4.2 of the main paper but extended by additional curves that are not observed over the whole grid domain. Following Elías et al. [31], we make use of empirical integrated depth functions that restrict the integration to the observed domains. However, we do not check whether such empirical depth functions fulfil our assumptions.

More precisely, in addition to all fully observed curves used in Section 4.2 of the main paper, we include all data from drifters that were in the relevant region at January 1st of the respective year and label all data from time points where the drifter was not in the relevant region as missing. Typically, these curves are densely observed for certain periods of the year with other periods where the observations are missing. In a few instances, there are only isolated observations at certain times. Figure 13 plots five of these types of curves to give an impression of the type of missingness.

The applied smoothing is as in Section 4.2 of the main paper, which results in assigning NA to a whole day within the smoothed curve when the raw data have a NA value that day. Smoothed curves that only contain NA-values are removed. Depending on the year, incorporating the curves with missing observations significantly increases the number of included drifters by a factor of at least two but sometimes even more than four.

TABLE 11
*p-values of the Joint-TP statistic with the integrated Tukey depth (IT) and the integrated simplicial depth (IS). 2015 is considered as very strong El Niño year, 2018 is considered a weak El Niño period and 2020–2021 is considered a moderate La Niña period.*

| **Sample 1** | **Sample 2** 2018 (470 curves) | 2020 (389 curves) | 2021 (309 curves) |
|---|---|---|---|
| 2015 (470 curves) | $0.227 \cdot 10^{-4}$ (IT)<br>$0.205 \cdot 10^{-4}$ (IS) | $0.332 \cdot 10^{-3}$ (IT)<br>$0.297 \cdot 10^{-3}$ (IS) | $0.446 \cdot 10^{-8}$ (IT)<br>$0.416 \cdot 10^{-8}$ (IS) |
| 2018 (470 curves) | | 0.880 (IT)<br>0.787 (IS) | 0.030 (IT)<br>0.039 (IS) |
| 2020 (389 curves) | | | 0.063 (IT)<br>0.087 (IS) |

TABLE 12
*p-values of the Joint-TP statistic with the integrated Tukey depth (IT) and the integrated simplicial depth (IS): 2011–2012 as well as 2020–2021 are considered as moderate La Niña periods. In 2002 a moderate El Niño occured.*

| **Sample 1** | **Sample 2** 2002 (325 curves) | Pooled 2011–2012 (640 curves) |
|---|---|---|
| Pooled 2020–2021 (698 curves) | 0.015 (IT)<br>0.011 (IS) | 0.085 (IT)<br>0.047 (IS) |

Table 11 is the analogue to Table 1 in the main paper, while Table 12 is the analogue to Table 2. All previously significant tests are still significant with decreased $p$-value, sometimes by magnitudes of order which we believe to be due to the largely increased sample size despite the fact that the additional data are only partially observed. Additionally, it appears, that the weak El Niño year of 2018 may be closer to the moderate La Niña period of 2020 than that of 2021, which was not indicated by the analysis with only the curves observed over the whole grid domain. This is consistent with the fact that 2020 was still relatively warm at the beginning of the year where most of the partial curves are observed, compare with the Oceanic Niño index as given in [64].

## Acknowledgments

We would like to thank Adam Sykulski (Imperial College London) for bringing the ocean drifter data to our attention. Felix Gnettner wants to thank Norbert Gaffke (Otto-von-Guericke-Universität Magdeburg) for a helpful discussion on the proof of Theorem 2.6 and Łukasz Smaga (Uniwersytet im. Adama Mickiewicza, Poznań) for sending some R-code on his tests for functional data.

## Funding

This work was supported by the Deutsche Forschungsgemeinschaft (DFG, German Research Foundation) - 314838170, GRK 2297 MathCoRe, as well as KI

1443/6-1. Additionally, it was supported by grant PID2022-139237NB-I00 funded by MCIN/AEI/10.13039/501100011033 and “ERDF A way of making Europe”.